\documentclass{article}
\usepackage[margin=3cm]{geometry}
\usepackage{mathtools, authblk, booktabs, xcolor}
\usepackage{amsmath, amssymb, amsthm}
\usepackage{subfig, caption, graphicx, enumitem, tikz}
\usepackage[bookmarks]{hyperref}
\hypersetup{
	colorlinks = true,
	linkcolor = blue,
	citecolor = red
}
\usepackage[maxbibnames=1000,backend=bibtex]{biblatex}
\newcommand{\uu}{{\mathbf u}}

\newtheorem{definition}{Definition}[section]
\newtheorem{theorem}{Theorem}[section]

\newtheorem{remark}[theorem]{Remark}
\newtheorem{example}{Example}

\begin{document}

\title{Spectral analysis of the stiffness matrix sequence in the approximated Stokes equation}

\author[1]{Samuele Ferri\thanks{sferri1@uninsubria.it}}
\author[2]{Chiara Giraudo\thanks{chiarg@math.uio.no}}
\author[1]{Valerio Loi\thanks{vloi@uninsubria.it}}
\author[3]{Miroslav Kuchta\thanks{miroslav@simula.no}}
\author[456]{Stefano Serra-Capizzano\thanks{s.serracapizzano@uninsubria.it}}
\affil[1]{\footnotesize Department of Theoretical and Applied Sciences, University of Insubria, Varese, Italy}
\affil[2]{Department of Mathematics, University of Oslo, Oslo, Norway}
\affil[3]{Department of Numerical Analysis and Scientific Computing, Simula Research Laboratory, Oslo, Norway}
\affil[4]{Department of Science and High Technology, University of Insubria, Como, Italy}
\affil[5]{Department of Information Technology, Uppsala University, Uppsala, Sweden}
\affil[6]{CLAP Research Center, University of Insubria, Como, Italy}

\date{}

\maketitle

\begin{abstract}
	In the present paper, we analyze in detail the spectral features of the matrix sequences arising from the Taylor-Hood $\mathbb{P}_2$-$\mathbb{P}_1$ approximation of variable viscosity for $2d$ Stokes problem under weak assumptions on the regularity of the diffusion. Localization and distributional spectral results are provided, accompanied by numerical tests and visualizations. A preliminary study of the impact of our findings on the preconditioning problem is also presented. A final section with concluding remarks and open problems ends the current work.
\end{abstract}

\section{Introduction}
The Stokes model describes the flow of a viscous, incompressible fluid in the absence of inertial effects. Such systems occur e.g. in geodynamics \cite{van2008community}, ice sheet modeling \cite{pattyn2007benchmark}, planetology \cite{kuchta2015despinning}, two-phase flows \cite{he2012preconditioning}, study of non-Newtonian fluids \cite{rajagopal2006implicit} or fluid-structure interaction problems \cite{wichrowski2023exploiting}. In these applications, the fluid viscosity typically depends on a different unknown model, for example, temperature \cite{farrell2022finite}, pressure, or shear rate \cite{hirn2012finite}. However, the resulting spatial variations of viscosity, in particular steep gradients, present a challenge for the construction of efficient iterative methods.\\
In this study, we analyze the effect of spatially varying viscosity on spectral features of the block
matrices that constitute discretization of the Stokes problem
\begin{equation}\label{eq: Stokes equations}
	\begin{aligned}
		-\nabla\cdot(\mu \nabla u) + \nabla p &= g \quad&\text{ in }\Omega,\\
		-\nabla\cdot u                        &=0\quad&\text{ in }\Omega,\\
		u &= 0 \quad&\text{ on }\partial\Omega.
	\end{aligned}
\end{equation}
Here $\Omega\subset\mathbb{R}^d$, $d=2$, $\mu:\Omega\rightarrow\mathbb{R}$ such that $\mu>0$ in $\Omega$, $g:\Omega\rightarrow\mathbb{R}^d$ is a given body force and the model seeks velocity $u:\Omega\rightarrow\mathbb{R}^d$ and pressure $p:\Omega\rightarrow\mathbb{R}$. The discretization of the equations in \eqref{eq: Stokes equations} then leads to a symmetric block-structured saddle-point system
\begin{equation}\label{eq:stokesblocks}
	\mathcal{A}_h
	\begin{bmatrix}
		\uu_h\\p_h
	\end{bmatrix} = \begin{bmatrix}
		{\mathbf g}_h\\0
	\end{bmatrix},\text{ where }
	\mathcal{A}_h=
	\begin{bmatrix}
		A_h & B_h^T\\
		B_h & 0
	\end{bmatrix}
\end{equation}
and $A_h$, $B_h$ discretize\footnote{We refer to Example \ref{ex:example} for details of the discretization.} respectively the diffusion operator $A=-\nabla(\mu\nabla\cdot)$ and the divergence operator $B=-\nabla\cdot$. Importantly, we assume that the discretization is stable so that \eqref{eq:stokesblocks} yields convergent approximations of \eqref{eq: Stokes equations}. Examples of such discretizations are staggered finite differences \cite{armfield1991finite} or the finite element method using Taylor-Hood elements \cite{brezzi1991stability} which shall be employed below.\\
In real-life applications, the system \eqref{eq:stokesblocks} needs to be solved by preconditioned Krylov or other iterative methods. Here, the canonical (block-diagonal) preconditioner is the operator $\mathcal{B}_h=\text{diag}(A_h, S_h)^{-1}$ with $S_h=B_hA^{-1}_hB^T_h$ being the Schur complement. As shown in \cite{murphy2000note} the spectrum of $\mathcal{B}_h\mathcal{A}_h$ contains only 3 distinct eigenvalues leading to rapid convergence of iterative solvers. However, due to the exact Schur complement, this preconditioner is prohibitively expensive. A practical preconditioner, based on the approximation of $S_h$, is the operator
\begin{equation}\label{eq:wathen}
	\mathcal{B}_h = \begin{bmatrix}
		A_h & 0\\
		0   & M_h\\
	\end{bmatrix}^{-1}
\end{equation}
with $M_h$ the discretization of the $\mu^{-1}$-weighted $L^2$-inner product, i.e. $\mu^{-1}$-weighted mass matrix.
In particular, \cite{grinevich2009iterative} shows that the spectrum of $M^{-1}_h S_h$ is bounded independent of $h$ (i.e. preconditioner \eqref{eq:wathen} is independent of the discretization parameter) but, critically, the spectral bounds depend on the smoothness of $\mu$. In turn, \eqref{eq:wathen} performs poorly in the presence of sharp gradients as observed, e.g. in \cite{rudi2017weighted} in benchmark problems inspired by magma dynamics. We demonstrate this behavior further in a simple setting in Example \ref{ex:example}.

\begin{example}[Variable viscosity and standard Stokes preconditioner]\label{ex:example}
	We consider \eqref{eq: Stokes equations} on domain $\Omega=(-1, 1)\times (-1, 1)$ and with the piecewise
	linear viscosity field
	\begin{equation}\label{eq:viscosity}
		\mu(x, y) = \begin{cases}
			\mu_1 & \lvert x \rvert < w,\\
			\mu_0 & \lvert x \rvert > w + \delta,\\
			\mu_0 + \frac{1}{\delta}\frac{w+\delta-\lvert x \rvert}{\mu_1 - \mu_0} & \text{otherwise}.
		\end{cases}
	\end{equation}
	Here $\mu_0>0$, $\mu_1>0$ and the parameters $w>0$ and $\delta\geq 0 $ control the width of the strips in which $\mu$ is constant and takes value $\mu_0$ or $\mu_1$, or respectively varies linearly between $\mu_0$ and $\mu_1$, For $\delta=0$ the viscosity is discontinuous.\\
	To discretize the model, we consider the finite element method with lowest order Taylor-Hood elements, i.e.	the approximate velocity and pressure $(u_h, p_h)\in V_h\times Q_h$ where
	\begin{displaymath}
		\begin{aligned}
			V_h&=\left\{v\in (H^1_0(\Omega_h))^2,\,v|_K \in (\mathbb{P}_2(K))^2\, \forall K\in \Omega_h \right\},\\
			Q_h&=\left\{q\in H^1(\Omega_h),\,q|_K \in \mathbb{P}_1(K)\,\forall K\in \Omega_h \right\}.
		\end{aligned}
	\end{displaymath}
	Here $\Omega_h$ is a triangular mesh\footnote{To resolve the sharp viscosity interface, the element size in the transition region is always smaller than $\delta$. Furthermore, the meshes always conform to the lines $x=\pm w$ and $x=\pm(w+\delta)$.} of $\Omega$ with characteristic mesh size $h>0$ and $\mathbb{P}_k(K)$ is the space of polynomials of degree $k$ on the element $K$. Denoting by $\uu_h$, $\mathbf{p}_h$ the coefficient vectors of $u_h\in V_h$, respectively $p_h\in Q_h$ and $\langle\cdot, \cdot\rangle$ the $l^2$-inner product, the operators in \eqref{eq:stokesblocks} and \eqref{eq:wathen} are defined as
	\begin{displaymath}
		\begin{aligned}
			&\langle A_h \uu_h, \mathbf{v}_h \rangle = \int_{\Omega_h} \mu \nabla u_h\cdot \nabla v_h,\quad
			\langle B_h \uu_h, \mathbf{p}_h \rangle = \int_{\Omega_h} p_h \nabla\cdot u_h,\\
			&\langle M_h \mathbf{p}_h, \mathbf{q}_h \rangle = \int_{\Omega_h} \mu^{-1} p_h q_h.
		\end{aligned}
	\end{displaymath}
	In the above integrals we approximate $\mu$ within each element by its midpoint value.\\	
	To investigate the effect of viscosity on \eqref{eq:wathen}-preconditioned Stokes problem we set $w=0.1$, $\mu_0=1$ and vary $1\leq \mu_1 \leq 10^8$, $0\leq \delta \leq 0.2$. In Figure \ref{fig:example} we then report the spectral condition number\footnote{The spectral condition number is defined in terms of the eigenvalues of $\mathcal{B}_h\mathcal{A}_h\uu_h=\lambda\uu_h$ as the ratio of the largest in magnitude and the smallest in magnitude non-zero eigenvalues. We recall that $\mathcal{A}_h$ is singular since functions in $Q_h$ do not necessarily have zero mean.} of the operators $\mathcal{B}_h\mathcal{A}_h$ when meshes $\Omega_h$ are refined. Therein we observe that the condition number blows up as the ratio $\mu_1/\mu_0$ increases. In comparison, there appears to be little effect of the parameter $\delta$. As expected the condition numbers are independent of the mesh size.\\
	We also evaluate {the} performance of the preconditioner in terms {of the} number of MinRes iterations\footnote{To deal with the singularity of $\mathcal{A}_h$ we pass the nullspace vector $(\mathbf{0}_h, \mathbf{1}_h)$ to the Krylov solver.} required to solve the Stokes model. Here, the iterations are started with 0 initial guess and use relative error tolerance of $10^{-12}$ as convergence criteria. In Figure \ref{fig:example}, we observe that the iteration count grows as $\mu_1/\mu_0$, however, the increase is not as dramatic as for the condition number, according to the potential theory and distribution results in \cite{beckermann2001superlinear,kuijlaars2006convergence}.
	
	\begin{figure}
		\centering
		\includegraphics[width=\textwidth]{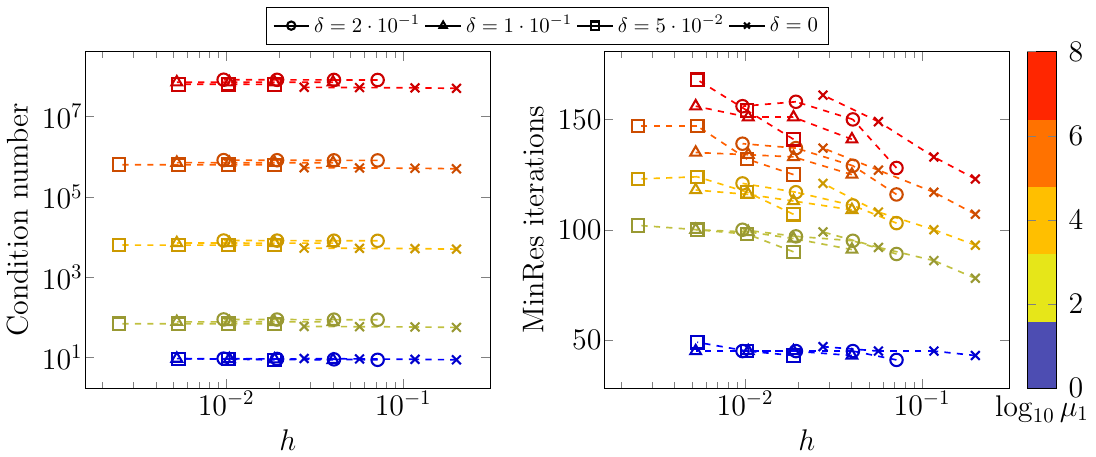}
		\vspace{-10pt}
		\caption{
			Effect of variable viscosity \eqref{eq:viscosity} on conditioning of the Stokes model
			\eqref{eq: Stokes equations} using the preconditioner \eqref{eq:wathen}.
		}
		\label{fig:example}
	\end{figure}
\end{example}

In recent years, several methods have been developed to mitigate the poor performance of \eqref{eq:wathen} for Stokes problems with high viscosity contrast. These include monolithic multigrid methods, e.g. \cite{borzacchiello2017box, wichrowski2022matrix}, where tailored block smoothers are a key ingredient for $\mu$-robustness. Another class of methods are block diagonal preconditioners, where we distinguish Schur complement and augmented Lagrangian (AL) preconditioners. The Schur complement-based preconditioners have been developed e.g. in \cite{rudi2015extreme, rudi2017weighted, burkhart2026robust}. In these works the inverse of the Schur complement is approximated using BFBT approach (derived from a Least-Squares Commutator \cite{elman2006block}) as
\begin{equation}
	S_h^{-1} \approx (B_hD_h^{-1}B_h^T)^{-1} (B_hD_h^{-1} A_hD_h^{-1}B_h^T)(B_hD_h^{-1}B_h^T)^{-1}.
\end{equation}
Here $D_h$ is a suitably chosen weighting matrix. Specifically, $D_h=\text{diag}A_h$ in \cite{rudi2015extreme}, while \cite{rudi2017weighted} set $D_h$ as the lumped matrix of the $\sqrt{\mu}$-weighted $L^2$-inner product and in \cite{burkhart2026robust} multiplication with $D_h$ amounts to applications of a multigrid cycle with matrix $A_h$. In general, these preconditioners yield faster convergence than \eqref{eq:wathen}, however, their performance degenerates for very strong viscosity variations. Finally, in AL approach the Stokes problem \eqref{eq:stokesblocks} is replaced by an equivalent parameter-dependent linear system for which the inverse Schur complement is easier to approximate robustly in $\mu$. In particular, \cite{he2012preconditioning,he2011preconditioning,shih2022robust} show that linear combinations of inverses of specifically $\mu$-weighted mass matrices are sufficient. However, the velocity block of the AL preconditioner, being now different than $A_h$, might require specialized approximate solvers, e.g. tailored multigrid methods as in \cite{shih2022robust}.\\
The above issues with convergence of Schur complement preconditioners motivate the need to study the effect of viscosity on the spectra of preconditioned Stokes problem. In particular, the condition number estimates from the continuous analysis, e.g. \cite{grinevich2009iterative} for the preconditioner \eqref{eq:wathen} or \cite{rudi2017weighted} for BFBT preconditioner, only provide (pessimistic, cf. Example \ref{ex:example}) bounds for convergence rates of Krylov solvers, while sharp estimates can be obtained from Weyl spectral distributions \cite{FEM-distr-conv,gergelits2019laplacian}, using tools from potential analysis \cite{beckermann2001superlinear,kuijlaars2006convergence}. We briefly review these tools in the following section.

\subsection*{Weyl distributions}
The study of Weyl spectral distribution results goes back to the work of Szeg\"o on Toeplitz matrix sequences $\{T_{n}(f)\}_n$ with $f \in L^\infty([-\pi, \pi])$ and real-valued; see e.g. \cite[Section 6.5]{GS-I} (and also Definition \ref{def-distribution} below) for the formal definitions.\\
After several decades, major advancements were made by Tyrtyshnikov in 1996 \cite{tyrtyshnikov1996unifying}, who developed a unifying approach for the spectral and singular value distribution of Toeplitz and multilevel Toeplitz matrices generated by symbols in {$L^2([-\pi, \pi]^d),  d\geq1$}, with the eigenvalue distribution proved under the restriction that $f$ is essentially real-valued. By relaxing the assumption on the generating function to the minimal requirement $f\in L^1([-\pi,\pi]^{d})$, the singular value distribution was extended in \cite{tyrtyshnikov1998spectra}, while, for the eigenvalue distribution, the function $f$ has to be essentially real-valued.\\
In 1998 Tyrtyshnikov and Zamarashkin \cite{tyrtyshnikov1998spectra} introduced matrix-theoretical approximation tools which anticipate the notion of approximating classes of sequences \cite{algebra}; see also Definition \ref{def:rectangular_acs}. A further significant contribution which is mathematically very elegant was provided by Tilli in \cite{TilliNota}, where he extended Szegő-type distribution in the case where the generating function is Lebesgue integrable and matrix-valued: the eigenvalue distribution is proven for Hermitian $d$-level, $s$-block Toeplitz matrix sequences, $d,s \geq 1$.
Conversely, when the operator involves variable coefficients as in the present work, the shift invariance is lost and the resulting matrices are no longer of Toeplitz type. This limitation naturally motivated the introduction of a more general framework capable of capturing the asymptotic spectral behavior of matrices arising from the discretization of variable-coefficient partial differential equations (PDEs). The challenge was addressed by Tilli with the theory of Locally Toeplitz (LT) sequences in \cite{TilliLT}.\\

Since these foundational contributions by Tyrtyshnikov and Tilli, there has been a surge of interest in the spectral analysis of structured matrix sequences, leading to the development of the comprehensive theory of Generalized Locally Toeplitz (GLT) sequences \cite{SerraCapizzano2003, SerraCapizzano2006}. This theoretical framework has rapidly become central in both pure and applied mathematics, due to its remarkable ability to capture the asymptotic singular value and eigenvalue distributions of matrices arising from virtually any meaningful approximation of (fractional) partial differential equations ((F)PDEs); see, e.g., the books and review papers~\cite{,block-glt-dD,GS-I,GS-II,garoni-tutorial} and references therein. In fact, virtually any practical discretization scheme of local type, such as finite differences and finite elements of any precision order, discontinuous Galerkin methods, finite volumes, isogeometric analysis, produces GLT matrix sequences.\\
Mathematically, for fixed positive integers $d$ and $r$, the set of $d$-level, $r$-block GLT matrix sequences forms a maximal $*$-algebra of matrix sequences, which is isometrically equivalent to the space of $d$-variate, $r \times r$ matrix-valued measurable functions defined on $[0,1]^d \times [-\pi,\pi]^d$. Each such GLT sequence $\{A_n\}_n$ is uniquely associated with a measurable, matrix-valued function $\kappa$ the GLT symbol on the domain $D = [0,1]^d \times [-\pi, \pi]^d$.
Notice that the set $[0, 1]^d$ can be replaced by any bounded Peano-Jordan measurable subset of $\mathbb{R}^d$ as occurring with the notion of reduced GLT $*$-algebras, see \cite{SerraCapizzano2003}[pp. 398-399, formula (59)] for the first occurrence with applications of approximated PDEs on general non-Cartesian domains in $d$ dimensions, \cite{SerraCapizzano2006}[Section 3.1.4] for the first formal proposal, and \cite{reducedGLT} for an exhaustive treatment, containing both the $*$-algebra theoretical results and some applications. This symbol provides a powerful tool for analyzing the singular value and eigenvalue distributions when the matrices $A_n$ are part of a matrix sequence of increasing size. The notation $\{A_n\}_n \sim_{\mathrm{GLT}} \kappa$ indicates that $\{A_n\}_n$ is a GLT sequence with symbol $\kappa$. Notably, the symbol of a GLT sequence is unique in the sense that if $\{A_n\}_n \sim_{\mathrm{GLT}} \kappa$ and $\{A_n\}_n \sim_{\mathrm{GLT}} \xi$, then $\kappa = \xi$ almost everywhere in $[0, 1]^d \times [-\pi, \pi]^d$. Furthermore, by the $*$-algebra structure, $\{A_n\}_n \sim_{\mathrm{GLT}} \kappa$ and $\{B_n\}_n \sim_{\mathrm{GLT}} \kappa$ implies that $\{A_n-B_n\}_n \sim_{\mathrm{GLT}} 0$, i.e., the matrix sequence $\{A_n-B_n\}_n$ is zero-distributed and the latter is very important for building explicit matrix sequences approximating a GLT matrix sequence and whose inversion is computationally cheap, in the context of preconditioning of large linear systems. It must be noticed that the $d$-level, $r$-block GLT $\ast$-algebra includes all Toeplitz sequences $\{T_n(f)\}_n$ with $f$ matrix-valued and $f \in L^{1}([-\pi,\pi]^{d})$ \cite{block-glt-1D,block-glt-dD}. Hence, already for $s=d=1$, the GLT class with $r=d=1$ represents a substantial generalization of the LT matrix sequences.\\

In this work, as a first step towards the spectral analysis, we employ the Toeplitz technology \cite{serra1999rate} and the theory of generalized locally Toeplitz (GLT) matrix sequences \cite{block-glt-dD} to analyse the spectral features of blocks $A_h$ and $B_h$. In particular, we study the distributions in terms of eigenvalues and singular values of the matrix sequences associated to $A_h$ and $B_h$ and the asymptotic behaviour of extreme eigenvalues of $A_h$, also in the variable coefficient setting. The study is critically discussed, with preliminary preconditioning proposals and with related numerical tests.\\
The current work is organized as follows. In Section \ref{sec:notation} we set the methodology while in Section \ref{sec:spectral} the main spectral tools are presented. Section \ref{sec:A and B} contains the main results while in Section \ref{sec:precond} we discuss a matrix approximation strategy based on the spectral analysis and we report few numerical experiments. Section \ref{sec:end} contains conclusions, open problems, and future directions of investigation.

\section{Methodology and notation}\label{sec:notation}
The present work contains a detailed study of the spectral (and singular value) distribution of the matrix sequences $\{A_n\}_n$ and $\{B_n\}_n$ {obtained} from a Taylor-Hood discretization of a variable-viscosity Stokes problem \eqref{eq:stokesblocks} with $n\sim 1/h$. Here $A_n$ identifies the symmetric stiffness block associated with the second derivatives in the space variables: in our setting, the contributions in $x$ and $y$ give rise to a structure {such} that the {entire} eigenvalue distribution is {determined} by one of them, since $\{A_n\}_n$ is block diagonal with equal blocks. The matrix $B_n=[B_{x,n}\;B_{y,n}]$ denotes the (rectangular) discrete divergence block: differently from $A_n$, it is independent of the viscosity coefficient $\mu$ and we will study the singular value distribution associated with its matrix sequence. To this end we heavily use the GLT theory (see \cite{GS-I,GS-II,block-glt-1D,block-glt-dD,rect-glt,reducedGLT,barbarino2020non} and references therein for the theoretical apparatus, and \cite{garoni-tutorial,garoni2019} for reviews and tutorials), especially in the block multilevel setting \cite{block-glt-dD} and in its rectangular extension \cite{rect-glt}. The aims of this work are as follows:
\begin{enumerate}
    \item understand if there is a useful multilevel/block structure in the matrices $A_n$ and $B_n$ that we can exploit for our spectral analysis (after suitable permutations and $O(n)$-rank corrections);
    \item modify the problem to apply the GLT spectral theory directly and without unnatural assumptions, both in the square/Hermitian case ($\{A_n\}_n$) and in the rectangular case ($\{B_n\}_n$).
    \item use the GLT theory to compute the GLT symbols: for $\{A_n\}_n$, a Hermitian matrix-valued symbol that, given the real symmetric nature of $\{A_n\}_n$, coincides with the spectral symbol; for $\{B_n\}_n$, a rectangular matrix-valued symbol describing its singular value distribution; 
    \item if feasible, combine the block information to analyse the saddle-point structure to design a useful preconditioner.
\end{enumerate}
To apply the block multilevel GLT machinery to our problem, we make use of suitable ad hoc projections and similarity transformations by permutation matrices. In particular, we appeal to \cite[Eq.\,(2.3) and Th.\,2.4.7 at page 136]{block-glt-dD}. The tools used are summarized in the subsequent list:
\begin{enumerate}
    \item spectral (resp.\ singular value) distribution invariance under (relatively) small dimensional compressions: the Hermitian case \cite[Theorem 4.3]{curl-div} and the rectangular case \cite[Theorem 2.4]{furci2024block};
    \item special similarity permutations for $d$-level $s$-block Toeplitz matrices \cite[Eq.\,(2.3) and Th.\,2.4.7 at page 136]{block-glt-dD} with $d\le 2$, $s>1$;
    \item the $d$-level $s$-block GLT $*$-algebra with $d\le 2$, $s>1$, together with its rectangular extension \cite{block-glt-dD,rect-glt};
    \item basic approximation theory when the viscosity coefficient $\mu$ is not constant.
\end{enumerate}
We describe the steps in each section.
The notation used follows the conventional standards taken from the books \cite{GS-I,GS-II} and from the papers \cite{block-glt-1D,block-glt-dD}, with $V^*$ denoting the conjugate transpose of $V$ and $V^{\dag}$ denoting the pseudo-inverse of $V$, $V$ any matrix with entries in the complex field.
Regarding our specific problem, a further customary and useful notation is reported below:
\begin{equation}\label{tridiag}
	\text{diag}(X) =
	\begin{bmatrix}
		X &  & \\
		 & \ddots & \\
		 &  & X\\
	\end{bmatrix},
	\qquad
	\text{tridiag}(X,D,Y) =
	\begin{bmatrix}
		D & Y & & \\
		X & \ddots & \ddots & \\
		 & \ddots & \ddots & Y\\
		 &  & X & D
	\end{bmatrix}.
\end{equation}

\section{Spectral analysis tools}\label{sec:spectral}
In order to be self-contained, in this section we review the essential technical tools for the spectral analysis of our saddle point problem blocks. Central to our approach is the concept of the singular value or spectral distribution, described by a matrix-valued function known as the spectral symbol. We use a formalism aligned with the well established theory and formalism for Toeplitz and GLT matrix sequences \cite{GS-I,GS-II,block-glt-1D,block-glt-dD}. \\
Our strategy uses the GLT algebra and its spectral theory as a foundation for analysing more complex matrix sequences. Indeed, using permutation similarities and semiorthogonal compressions, we recast our problem into a suitable form for GLT analysis. In the following subsections, we recall the definitions and spectral properties of Toeplitz and GLT sequences and review methods for handling block-structured and compressed sequences.

\subsection{Asymptotic spectral symbol and approximation tools}\label{ssec:approximation}
We begin by defining the concept of asymptotic (spectral) singular value distribution for matrix sequences
\begin{definition}\label{def-distribution}
	Let $f \colon D \,\to\, \mathbb{C}^{s\times s}$	be a measurable matrix-valued function whose eigenvalues and singular values are denoted by
	\[
	\{\lambda_i(f)\}_{i=1}^s
	\quad\text{and}\quad
	\{\sigma_i(f)\}_{i=1}^s,
	\]
	respectively. Assume that \(D \subset \mathbb{R}^d\) is Lebesgue measurable with
	\[
	0 \;<\; \mu_d(D) \;<\; \infty,
	\]
	where \(\mu_d\) is the Lebesgue measure in \(\mathbb{R}^d\). Let \(\{A_n\}_{n}\) be a sequence of matrices such that
	\[
	\dim(A_n) \;=\; d_n \;\xrightarrow[n\to\infty]{}\; \infty,
	\]
	and let the eigenvalues and singular values of \(A_n\) be
	\[
	\{\lambda_j(A_n)\}_{j=1}^{d_n}
	\quad\text{and}\quad
	\{\sigma_j(A_n)\}_{j=1}^{d_n},
	\]
	respectively. We introduce the following notions:
	\begin{itemize}
		\item \emph{Eigenvalue distribution.}
		We say that \(\{A_n\}_{n}\) is \emph{distributed as \(f\) over \(D\) in the sense of the eigenvalues}, and write
		\[
		\{A_n\}_{n}\;\sim_{\lambda}\;(f,D),
		\]
		if for every continuous function \(F\) with compact support,
		\begin{equation}\label{distribution:eig}
			\lim_{n\to\infty}\frac{1}{d_n}\sum_{j=1}^{d_n}F\bigl(\lambda_j(A_n)\bigr)
			\;=\;
			\frac{1}{\mu_d(D)} \int_{D}
			\frac{1}{s}\sum_{i=1}^s F\bigl(\lambda_i(f(t))\bigr)
			\,\mathrm{d}t.
		\end{equation}
		In this case, we say that \(f\) is the \emph{spectral symbol} of \(\{A_n\}_{n}\).
		
		\item \emph{Singular value distribution.}
		We say that \(\{A_n\}_{n}\) is \emph{distributed as \(f\) over \(D\) in the sense of the singular values}, and write
		\[
		\{A_n\}_{n} \;\sim_{\sigma}\; (f,D),
		\]
		if for every continuous function \(F\) with compact support,
		\begin{equation}\label{distribution:sv}
			\lim_{n\to\infty}\frac{1}{d_n}\sum_{j=1}^{d_n}F\bigl(\sigma_j(A_n)\bigr)
			\;=\;
			\frac{1}{\mu_d(D)} \int_{D}
			\frac{1}{s}\sum_{i=1}^s F\bigl(\sigma_i(f(t))\bigr)
			\,\mathrm{d}t.
		\end{equation}
		In this case, \(f\) is referred to as the \emph{singular value symbol} of \(\{A_n\}_{n}\).
		
		\item \emph{Rectangular case.}
		The notion \(\{A_n\}_{n}\sim_{\sigma}(f,D)\) also applies to the rectangular setting where \(f \colon D \to \mathbb{C}^{s_1\times s_2}\). In that scenario, the index \(s\) in \eqref{distribution:sv} is replaced by \(\min(s_1,s_2)\), and \(A_n \in \mathbb{C}^{d_n^{(1)}\times d_n^{(2)}}\) with \(d_n\) in \eqref{distribution:sv} given by \(\min(d_n^{(1)}, d_n^{(2)})\). Naturally, the eigenvalue distribution notion does not extend to rectangular matrices.
	\end{itemize}
\end{definition}

Throughout the paper, when the domain can be easily inferred from the context, we replace the notation $\{A_n\}_n\sim_{\lambda,\sigma}(f,D)$ with $\{A_n\}_n\sim_{\lambda,\sigma} f$.
\begin{remark}\label{simply}
	{W}e assume that every matrix sequence we consider has a measurable symbol $f$, in terms of singular or eigenvalues. This property helps in algebraic approximation techniques, as they are sparsely unbounded (s.u.) (see \cite[Definition 5.2]{GS-I} for the definition of s.u. sequences and \cite[Proposition 5.4]{GS-I} for its relation with the asymptotic singular value distribution notion). Unless stated otherwise, in all the reasoning we make in this work, we assume that the matrix sequences are distributed according to some function $f$, either in the sense of the singular {values} or in the sense of the eigenvalues.
\end{remark}

\begin{remark}\label{rem:approx}
	If $f$ is smooth enough, an informal interpretation of the limit relation \eqref{distribution:eig} (resp. \eqref{distribution:sv}) is that when $n$ is sufficiently large, the eigenvalues (resp. singular values) of $A_{n}$ can be subdivided into $s$ different subsets of the same cardinality. Then $d_n/s$ eigenvalues (resp. singular values) of $A_{n}$ can	be approximated by a sampling of $\lambda_1(f)$ (resp. $\sigma_1(f)$) on a uniform equispaced grid of the domain $D$, and so on until the	last $d_n/s$ eigenvalues (resp. singular values), which can be approximated by an equispaced sampling of $\lambda_s(f)$ (resp. $\sigma_s(f)$) in the domain $D$.
\end{remark}

\begin{remark}
	We say that $\{A_n\}_n$ is \textit{zero-distributed} in the sense of the eigenvalues if $\{A_n\}_n\sim_\lambda 0$ and all the $A_n$ are Hermitian (in general only the quasi-Hermitian character of the matrix sequence $\{A_n\}_n$ is required; see \cite{barbarino2020non}).
\end{remark}

A key tool to compute singular and eigenvalue distributions is to approximate matrix sequences with more manageable ones. This approximation is considered in the so-called approximating class of sequences (a.c.s.) topology, which is described concretely by the a.c.s. limit notion we describe below. We limit ourselves to describing it in the elementary matrix sequence approximation perspective.

\begin{definition}[{\rm \cite[Definition 2.6]{rect-glt}}]
	\label{def:rectangular_acs}
	Let $\{A_n\}_{n\in\mathbb{N}}$ be a sequence of matrices, where $A_n$ has size $d_n \times e_n$.
	Let $\{\{B_{n,m}\}_{n\in\mathbb{N}}\}_{m\in\mathbb{N}}$ be a double sequence of matrices, with
	$B_{n,m}$ also of size $d_n \times e_n$.
	
	We say that $\{\{B_{n,m}\}_n\}_m$ is an approximating class of sequences (a.c.s.) for
	$\{A_n\}_n$, and we write
	\begin{equation}\label{eq:acs_convergence}
		\bigl\{B_{n,m}\bigr\}_n \;\xrightarrow{\text{a.c.s.}}\; \bigl\{A_n\bigr\}_n,
	\end{equation}
	if there exist functions $c : \mathbb{N} \to \mathbb{R}_{\ge 0}$ and
	$\omega : \mathbb{N} \to \mathbb{R}_{\ge 0}$ such that
	\[
	\lim_{m \to \infty} c(m) \;=\; 0
	\quad\text{and}\quad
	\lim_{m \to \infty} \omega(m) \;=\; 0,
	\]
	and for every $m \in \mathbb{N}$ there is an $n_m \in \mathbb{N}$ with the property that for all $n \ge n_m$
	there exist matrices $R_{n,m}$ and $N_{n,m}$ (both of size $d_n \times e_n$) satisfying:
	\begin{enumerate}
		\item $A_n \;=\; B_{n,m} \;+\; R_{n,m} \;+\; N_{n,m}$,
		\item $\mathrm{rank}\!\bigl(R_{n,m}\bigr) \;\le\; c(m)\,\bigl(d_n \wedge e_n\bigr)$,
		\item $\|N_{n,m}\| \;\le\; \omega(m)$,
	\end{enumerate}
	where $d_n \wedge e_n := \min\{\,d_n,\, e_n\,\}$, and $\|\cdot\|$ denotes the spectral norm.
\end{definition}

We observe that the notion of convergence described in Definition \ref{def:rectangular_acs} identifies matrix sequences that are equal to each other up to small perturbations. More precisely, we say that two matrix sequences $\{A_n\}_n$ and $\{B_n\}_n$ are a.c.s. equivalent, in symbol
\begin{equation}\label{eq:acs_equivalence}
	\bigl\{
	A_n
	\bigr\}_{n}
	\;\equiv_{\mathrm{a.c.s.}}\;
	\bigl\{
	B_n
	\bigr\}_{n},
\end{equation}
if
\begin{enumerate}
	\item $A_n \;=\; B_{n} \;+\; R_{n} \;+\; N_{n}$,
	\item $\mathrm{rank}\,\!\bigl(R_{n}\bigr) \; = \; o(d_n \wedge e_n)$,
	\item $\|N_{n}\| \;\to \; 0$.
\end{enumerate}
(\ref{eq:acs_equivalence}) is clearly an equivalence relation. The reason a.c.s convergence is a tool to compute spectral and singular value distributions is the following:
\begin{theorem}[{\rm \cite[Theorem 5.6]{GS-I}}]\label{thm:acs_conv}
	Let \(\{A_n\}_{n}\) and \(\{B_{n,m}\}_{n}\) be sequences of matrices, and let
	$f, f_{m}: D \subset \mathbb{R}^{k} \rightarrow$ $\mathbb{C}$ be measurable functions.
	Suppose the following conditions hold:
	\begin{enumerate}
		\item \(\{B_{n,m}\}_{n} \sim_{\sigma} f_{m}\) for every \(m\).
		\item \(\{B_{n,m}\}_{n} \xrightarrow{\mathrm{a.c.s.}} \{A_n\}_{n}\).
		\item \(f_{m} \to f\) in measure.
	\end{enumerate}
	Then it follows that \(\{A_n\}_{n} \sim_{\sigma} f\). If all the sequences are Hermitian, the distributional results also hold for the eigenvalue symbol.
\end{theorem}

As a consequence of Theorem \ref{thm:acs_conv}, we can state the following result.
\begin{theorem}\label{thm:acs_equiv_distributions}
	Let \(\{A_n\}_{n}\) and \(\{B_n\}_{n}\) be sequences of matrices such that
	\[
	\bigl\{A_n\bigr\}_n
	\;\equiv_{\mathrm{a.c.s.}}\;
	\bigl\{B_n\bigr\}_n.
	\]
	Then, for any function \(f\), we have
	\[
	\big\{A_n\bigr\}_n \sim_{\sigma} f
	\quad \Longleftrightarrow \quad
	\bigl\{B_n\bigr\}_n \sim_{\sigma} f.
	\]
	Moreover, if \(\{A_n\}_n\) and \(\{B_n\}_n\) are Hermitian, the same equivalence holds for their eigenvalue distributions.
\end{theorem}

Given that we are dealing with distributed matrix sequences and adhering to the convention outlined in Remark \ref{simply}, nice algebraic properties for a.c.s. convergence hold. It is important to note that the assertions of the subsequent theorem trivially apply to a.c.s. equivalence as well.

\begin{theorem}[{\rm \cite[Proposition 5.2 and Proposition 5.5]{GS-I}}]\label{thm:acs-properties}
	
	Let \(\{A_n\}_n\), \(\{A'_n\}_n\), \(\{B_{n,m}\}_{n,m}\), and \(\{B'_{n,m}\}_{n,m}\) be matrix sequences such that
	\[
	\bigl\{B_{n,m}\bigr\}_{m} \; \xrightarrow{\mathrm{a.c.s.}} \;\bigl\{A_n\bigr\}_n
	\quad\text{and}\quad
	\bigr\{B'_{n,m}\bigl\}_{m} \; \xrightarrow{\mathrm{a.c.s.}} \; \bigl\{A'_n\bigl\}_n.
	\]
	Consequently, assuming matrix sizes permit these operations, the following algebraic properties are valid:
	\begin{enumerate}
		\item \[
		\bigl\{ B_{n,m}^*\bigr\}_m
		\; \xrightarrow{\mathrm{a.c.s.}} \;
		\bigl\{ A_n^*\bigr\}_n.
		\]
		\item For all \(\alpha, \beta \in \mathbb{C}\),
		\[
		\bigl\{\alpha B_{n,m} + \beta B'_{n,m}\bigr\}_m
		\; \xrightarrow{\mathrm{a.c.s.}} \;
		\bigl\{\alpha A_n + \beta A'_n\bigr\}_n.
		\]
		\item
		\[
		\bigl\{B_{n,m} B'_{n,m}\bigr\}_m
		\; \xrightarrow{\mathrm{a.c.s.}} \;
		\bigl\{A_n A'_n\bigr\}_n.
		\]
		\item For all \(\{C_n\}_n\)
		\[
		\bigl\{B_{n,m} C_n\bigr\}_m
		\; \xrightarrow{\mathrm{a.c.s.}} \;
		\bigl\{A_n C_n\bigr\}_n.
		\]
	\end{enumerate}
\end{theorem}

\subparagraph{Toeplitz structures and GLT tools}\label{ssec:GLTbackground}
We begin by formalizing the notion of block multilevel Toeplitz and block multilevel circulant matrix sequences associated with a matrix-valued Lebesgue integrable function. We then introduce the concepts of eigenvalue (spectral) and singular value distribution and recall the basic tools from the GLT theory.

\begin{definition}[Toeplitz, block-Toeplitz, block multilevel Toeplitz matrices]\label{def_not1}
	A (finite- or infinite-dimensional) Toeplitz matrix $T$ has constant entries along each descending diagonal from left to right. Concretely, we have
	\begin{equation}\label{eq_toeplits}
		T \;=\;
		\begin{bmatrix}
			a_0 & a_{-1} & a_{-2} & \cdots\\
			a_{1} & a_0 & a_{-1} & \cdots\\
			a_{2} & a_{1} & \ddots & \ddots\\
			\vdots & \vdots & \ddots & \ddots
		\end{bmatrix},
		\quad
		T_n \;=\;
		\begin{bmatrix}
			a_0 & a_{-1} & a_{-2} & \cdots & a_{-n+1}\\
			a_{1} & a_0   & a_{-1} & \ddots & \vdots\\
			\vdots & \ddots & \ddots & \ddots & a_{-1}\\
			a_{n-1} & \cdots & a_2 & a_1 & a_0
		\end{bmatrix},
	\end{equation}
	where $T$ can be infinite-dimensional, while $T_n$ is $n\times n$. We consider matrix sequences $\{T_n\}_n$ where $n \to \infty$.
	
	In general, each $a_k$ may be an $s_1 \times s_2$ matrix (a \emph{block}), giving a \emph{block Toeplitz} matrix. Then $T_n$ is partitioned into $n\times n$ blocks, each of size $s_1\times s_2$, so $T_n$ has total size $N_1\times N_2$, with $N_1 = n\,s_1$ and $N_2 = n\,s_2$. In this notation, we may also write $X_{N_1,N_2} = T_n$.
	
	A block \emph{multilevel} Toeplitz matrix is obtained recursively by letting each block itself be a unilevel or multilevel Toeplitz matrix. Equivalently, a standard (unilevel) Toeplitz matrix is a special case of this construction.
\end{definition}

\begin{definition}[Toeplitz sequences and generating functions]\label{def_not3}
	Let $f$ be a $d$-variate, complex-valued, Lebesgue integrable function defined on $Q^d = [-\pi,\pi]^d$, where $d \ge 1$ and $\mu_d(Q^d)=(2\pi)^d$. Denote by $f_k$ its Fourier coefficients,
	\[
	f_k
	\;=\; \frac{1}{(2\pi)^d} \int_{Q^d} f(\theta)\,e^{-\,i\,(k,\theta)}\,d\theta,
	\quad
	k = (k_1,\dots,k_d)\,\in \mathbb{Z}^d,
	\]
	with $i^2=-1$, $\theta=(\theta_1,\dots,\theta_d)$, and $(k,\theta)=\sum_{j=1}^d k_j\,\theta_j$.
	Define
	\[
	T_n \;=\; \{\,f_{\,k-\ell}\}_{k,\ell=\mathbf{e}^T}^{\,n}
	\;\in\; \mathbb{C}^{\,N(n)\times N(n)},
	\]
	where $\mathbf{e}=[1,1,\dots,1]\in \mathbb{N}^d$ and $N(n)=n_1\cdots n_d$. Then $\{T_n\}_n$ is a \emph{Toeplitz matrix sequence generated by $f$}, denoted $T_n(f)$.
	
	For $d=1$, the matrix
	\[
	T_n(f) \;=\;
	\begin{bmatrix}
		f_0 & f_{-1} & \cdots & f_{-n+1}\\
		f_1 & f_0 & \ddots & \vdots\\
		\vdots & \ddots & \ddots & f_{-1}\\
		f_{n-1} & \cdots & f_{1} & f_{0}
	\end{bmatrix}.
	\]
	For $d=2$, say $n=(2,3)$, one obtains a block structure:
	\[
	T_n(f)
	\;=\;
	\begin{bmatrix}
		F_0 & F_{-1}\\
		F_1 & F_0
	\end{bmatrix},
	\quad
	F_k
	\;=\;
	\begin{bmatrix}
		f_{(k,\,0)} & f_{(k,\, -1)} & f_{(k,\, -2)}\\
		f_{(k,\,1)} & f_{(k,\, 0)} & f_{(k,\, -1)}\\
		f_{(k,\,2)} & f_{(k,\, 1)} & f_{(k,\,0)}
	\end{bmatrix},
	\quad k=0,\pm1.
	\]
	We call $f$ the \emph{generating function} (or \emph{symbol}) of $T_n(f)$ and of the sequence $\{T_n(f)\}_n$.
	
	If $f$ takes values in $\mathbb{C}^{s_1\times s_2}$ (matrix-valued) and is still integrable over $Q^d$, then $T_n(f)$ is a block $d$-level Toeplitz matrix of size $s_1 N(n)\times s_2 N(n)$, and $f$ is again its generating function. If $s_1=s_2=s$, we write $f\in L^1(Q^d,s)$ and $T_n(f)$ has square size $s\,N(n)\times s\,N(n)$.
\end{definition}

The following theorem by Tilli generalizes a century of results by Szeg\H{o}, Widom, Avram, Parter, Tyrtyshnikov, and Zamarashkin (see \cite{GS-II,block-glt-dD} and references therein).
\begin{theorem}[{\cite{TilliNota}}]\label{szego-herm}
	Let $f \in L^1(Q^d,\,s_1 \times s_2)$. Then the sequence $\{T_{n}(f)\}_{{n}}$ is distributed in the singular value sense as $f$, i.e.\ $\{T_{n}(f)\}_n \sim_\sigma(f, Q^d)$. If $s_1=s_2=s$ and $f$ is Hermitian matrix-valued, then $\{T_{n}(f)\}_{{n}}\sim_\lambda(f,Q^d)$.
\end{theorem}

\begin{remark}
	The general tools described here can be applied in a variety of contexts, including linear systems arising from isogeometric analysis or finite-volume schemes; see the applications described in \cite{GS-II,block-glt-dD,garoni2019}.
\end{remark}

\subsubsection{Definition and properties of the block multilevel GLT class}\label{GLT}
The multilevel Toeplitz, diagonal sampling and zero distributed matrix sequences, along with the notion of approximating class of sequences (a.c.s.) (see \cite{algebra} Definition 2.1) are the main building blocks of the GLT class. Instead of a rigorous definition and further technical details we here give six basic axioms of the class, proven to be sufficient for studying the spectral and singular value distribution of matrix sequences of interest. It should be noted that the whole set of axioms characterizes uniquely the GLT class and they can be considered as an alternative more friendly definition (see \cite{block-glt-dD}).
\begin{description}[leftmargin=1.75em] 
	\item[\textbf{GLT 1}]
	A $d$-level GLT sequence $\{A_n\}_n$ of $d_n\times d_n$ matrices is associated with a unique Lebesgue measurable function
	\begin{displaymath}
		\kappa : [0,1]^d \times [-\pi,\pi]^d \;\longrightarrow\; \mathbb{C}^{s\times s},
	\end{displaymath}
	called its \emph{symbol}, and we write
	\begin{displaymath}
		\bigl\{A_n\bigr\}_n \,\sim_{\textsc{glt}}\, \kappa.
	\end{displaymath}
	The singular values of $\{A_n\}_n$ are distributed as $\kappa$. If each $A_n$ is Hermitian, then the eigenvalues are distributed as $\kappa$.
	
	\item[\textbf{GLT 2}]
	Every $d$-level Toeplitz sequence whose generating function
	\begin{displaymath}
		 f\in L^1\bigl([-\pi,\pi]^d, \mathbb{C}^{s \times s} \bigr)
	\end{displaymath}
	is a GLT sequence with symbol $\kappa(x,\theta) = f(\theta)$.
	
	\item[\textbf{GLT 3}]
	Every diagonal matrix sequence whose entries are the uniform sampling of an a.e.\ continuous function
	\begin{displaymath}
		\alpha \; : \; [0,1]^d \;\to\; \mathbb{C}^{s\times s}
	\end{displaymath}
	is a GLT sequence with symbol $\kappa(x,\theta) = \alpha(x)$.
	
	\item[\textbf{GLT 4}]
	Every zero-distributed matrix sequence is a GLT sequence with symbol $\kappa(x,\theta)=0$.
	
	\item[\textbf{GLT 5}]
	The collection of all GLT sequences forms a unital $*$-algebra. In particular, if $\{A_n\}_n \sim_{\textsc{glt}} \kappa$ and $\{B_n\}_n \sim_{\textsc{glt}} \xi$, then:
	\begin{enumerate}[label=(\alph*), topsep=0.3em, itemsep=0.3em] 
		\item $\{A_n^*\}_n \sim_{\textsc{glt}} \kappa^* $
		\item $\{\alpha A_n + \beta B_n\}_n \sim_{\textsc{glt}} \alpha\,\kappa + \beta\,\xi$ for any $\alpha,\beta \in \mathbb{C}$.
		\item $\{A_n B_n\}_n \sim_{\textsc{glt}} \kappa\,\xi$, whenever $\kappa$ and $\xi$ are block-dimension compatible.
		\item If $\kappa$ is invertible a.e., then $\{A_n^{\dag}\}_n \sim_{\textsc{glt}} \kappa^{-1}$. 
	\end{enumerate}
	
	\item[\textbf{GLT 6}]
	$\{A_n\}_n \sim_{\textsc{glt}} \kappa$ if and only if there exist GLT sequences
	\begin{displaymath}
		\{B_{n,m}\}_n \;\sim_{\textsc{glt}}\; \kappa_m
	\end{displaymath}
	such that $\kappa_m \to \kappa$ in measure and $\{\{B_{n,m}\}_n\}_m$ is an a.c.s.\ for $\{A_n\}_n$.
\end{description}

\noindent
\textbf{Rectangular GLT sequences.}
In our applications, we need to extend the notion of GLT sequences in the rectangular case accordingly. We refer to \cite{rect-glt} for the technical details and relationships with the classical GLT algebra. We outline the main building blocks:
\begin{enumerate}[label=\textbf{(\arabic*)}, topsep=0.3em, itemsep=0.3em] 
	\item \textbf{$(r,s)$-Toeplitz sequences.}
	Let $f : [-\pi, \pi]^d \to \mathbb{C}^{r \times s}$ belong to the space $L^1\bigl([-\pi,\pi]^d, \mathbb{C}^{r\times s}\bigr)$.
	The sequence $\{T_n(f)\}_n$ of rectangular Toeplitz matrices generated by $f$ is a rectangular GLT sequence with symbol $\kappa(\theta) = f(\theta)$.
	
	\item \textbf{$(r,s)$-Diagonal sequences.}
	Let $a : [0,1]^d \to \mathbb{C}^{r \times s}$ be a.e.\ continuous.
	The family $\{D_n(a)\}$ of diagonal sampling matrices is a rectangular GLT sequence with symbol $\kappa(x) = a(x)$.
	
	\item \textbf{$(r,s)$-Zero-distributed sequences.}
	A matrix sequence $\{Z_n\}$ is $(r \times s)$-zero-distributed if
	\[
	Z_n = N_n + E_n,
	\quad \text{where}\quad \lim_{n\to\infty}\!\Bigl(\|E_n\|_2 + \frac{\operatorname{rank}(N_n)}{n}\Bigr) = 0.
	\]
	Such a sequence has symbol $\kappa(x,\theta)=0$.
\end{enumerate}

\subsection{Properties}
For rectangular GLT sequences, the analogue of \textbf{GLT 1--6} holds, provided block dimensions are compatible for sums or products. Specifically, if $\{A_n\}_n \sim_{\textsc{GLT}} \kappa$ and $\{B_n\}_n \sim_{\textsc{GLT}} \xi$, then
\begin{equation}
    \bigl\{\alpha A_n + \beta B_n\bigr\}_n \;\sim_{\textsc{glt}}\; \alpha\,\kappa + \beta\,\xi, \quad\text{and}\quad \bigl\{A_n\,B_n\bigr\}_n \;\sim_{\textsc{glt}}\; \kappa\,\xi,
\end{equation}
if the matrix operations are well defined. If the singular values of $\kappa$ are non-zero a.e., then
\begin{equation}
    \bigl\{A_n^\dagger\bigr\}_n \;\sim_{\textsc{glt}}\; \kappa^\dagger.
\end{equation}

\subsection{Block-structured sequences and the extradimensional approach}
In our applications, matrices and matrix sequences arise naturally from functional calculus and algebra of matrix sequences and compressed block structure algebras, where varying block sizes are managed via partial isometries (compressions, unitary matrices, etc.). The key objects remain GLT sequences; here, we clarify how GLT theory applies in our context.\\
We handle block structures through permutation similarities and semi-orthogonal compressions, translating the original spectral analysis into an equivalent GLT-based problem. To this end, we briefly recall how permutations and embeddings preserve the spectral distribution of GLT sequences. We also highlight precisely when and how GLT algebra and its functional calculus transfer to these compressed or similar algebras. We remark that the approximation theory for compressed, unitarily equivalent sequences with varying spectral support is summarized in \cite[Definition 2.1]{gltmoving} via the notion of \emph{generalized approximated class of sequences}; but this additional formalism is not required here.\\
We begin by introducing the basic permutation we use:
\begin{equation}
	P_{k_1, k_2} \;= \;
	\begin{bmatrix}
		I_{k_1} \otimes (\mathbf{e}_1^{(k_2)})^T \\
		I_{k_1} \otimes (\mathbf{e}_2^{(k_2)})^T \\
		\vdots \\
		I_{k_1} \otimes (\mathbf{e}_{k_2}^{(k_2)})^T
	\end{bmatrix}
	\; = \;
	\sum_{i=1}^{k_2} \mathbf{e}_i^{(k_2)} \otimes I_{k_1} \otimes (\mathbf{e}_i^{(k_2)})^T
\end{equation}
and we further define
\begin{equation}
	\Pi_{\mathbf{n},s,r}=P_{s, N(n)} \otimes I_r.
\end{equation}
These are the matrices that we use for the permutation theorems of blocking structures. We begin with the very general theorem for handling multilevel block structures with rectangular GLT symbols.

\begin{theorem}[{\rm \cite[Theorem 4.2]{rect-glt}}]\label{thm:permutation_rect_eq}
	Let \(i,j = 1,\dots,s\), and suppose that \(\{A_{n,ij}\}_{n}\) is a
	\(d\)-level \((r,t)\)-block GLT sequence with symbol
	\(\kappa_{ij} : [0,1]^d \times [-\pi,\pi]^d \to \mathbb{C}^{r \times t}\).
	Define
	\[
	A_n \;=\;
	\begin{bmatrix}
		A_{n,11} & \cdots & A_{n,1s} \\[6pt]
		\vdots   & \ddots & \vdots   \\[6pt]
		A_{n,s1} & \cdots & A_{n,ss}
	\end{bmatrix}
	\quad\text{and}\quad
	\kappa \;=\;
	\begin{bmatrix}
		\kappa_{11} & \cdots & \kappa_{1s} \\[4pt]
		\vdots      & \ddots & \vdots      \\[4pt]
		\kappa_{s1} & \cdots & \kappa_{ss}
	\end{bmatrix}.
	\]
	Then the sequence
	\(\{\Pi_{n,s,r}\,A_n\,\Pi_{n,s,t}^{T}\}_{n}\)
	is a \(d\)-level $(sr,st)$-block GLT sequence with symbol \(\kappa\).
\end{theorem}

We approximate matrix sequences of possibly different (yet asymptotically equal) sizes via semi-orthogonal compressions. The following theorems provide conditions that ensure that the compressed sequences retain the asymptotic eigenvalue or singular value distribution of the original sequences (see also \cite{Blo_I, Blo_II} for more refined results).
\begin{theorem}[{\rm \cite[Theorem 2.4]{furci2024block}}]\label{thm:compression-equivalence}
	Let \(\{A_n\}_n\) be a (rectangular) matrix sequence with each
	\(A_n \in \mathbb{C}^{n^{(1)}\times n}\) and \(n^{(1)} \ge n\).
	Suppose \(P_n \in \mathbb{C}^{n\times n'}\) and
	\(P_{n^{(1)}} \in \mathbb{C}^{n^{(1)}\times n^{(1)'}}\) are compression matrices
	such that \(n'<n\), \(n^{(1)'} < n^{(1)}\), and
	\[
	P_n^*\,P_n \;=\; I_{n'} \quad\text{and}\quad
	P_{n^{(1)}}^*\,P_{n^{(1)}} \;=\; I_{n^{(1)'}}.
	\]
	Define
	\[
	Y_{n'} \;=\; P_{n^{(1)}}^*\,A_n\,P_n.
	\]
	Under the assumption that
	\[
	\lim_{n,n'\to\infty}\,\frac{n'}{n}
	\;=\;
	\lim_{n^{(1)},\,n^{(1)'}\to\infty}\,\frac{n^{(1)'}}{n^{(1)}}
	\;=\; 1
	\]
	we have
	\[
	\{A_n\}_n \;\sim_{\sigma}\; f
	\quad\Longleftrightarrow\quad
	\{Y_{n'}\}_{n'} \;\sim_{\sigma}\; f.
	\]
\end{theorem}
\begin{theorem}[{\rm \cite[Theorem 4.3]{curl-div}}]
	\label{thm:hermitian-compression}
	Let \(\{X_n\}_n\) be a Hermitian matrix sequence with each
	\(X_n \in \mathbb{C}^{n\times n}\) and \(X_n = X_n^*\).
	Suppose \(P_n \in \mathbb{C}^{n\times n'}\) is a compression matrix
	satisfying \(n' < n\) and \(P_n^*\,P_n = I_{n'}\).
	Define
	\[
	Y_{n'} \;=\; P_n^*\,X_n\,P_n.
	\]
	Under the assumption
	\[
	\lim_{n,\,n'\to\infty}\,\frac{n'}{n} \;=\; 1
	\quad\bigl(\text{i.e., }n \;=\; n' + o(n)\bigr),
	\]
	we have
	\[
	\{X_n\}_n \;\sim_{\lambda}\; f
	\quad\Longleftrightarrow\quad
	\{Y_{n'}\}_{n'} \;\sim_{\lambda}\; f.
	\]
\end{theorem}

\section{The block multilevel structures}\label{sec:A and B}
\begin{figure}
	\centering
	\includegraphics[width=0.32\textwidth]{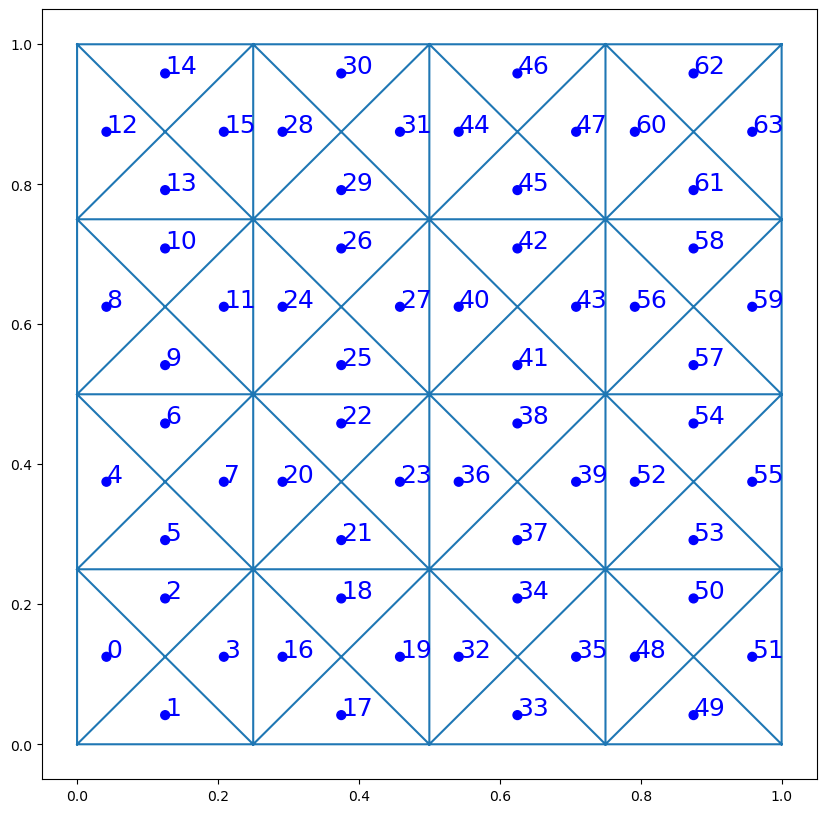} 
	\includegraphics[width=0.32\textwidth]{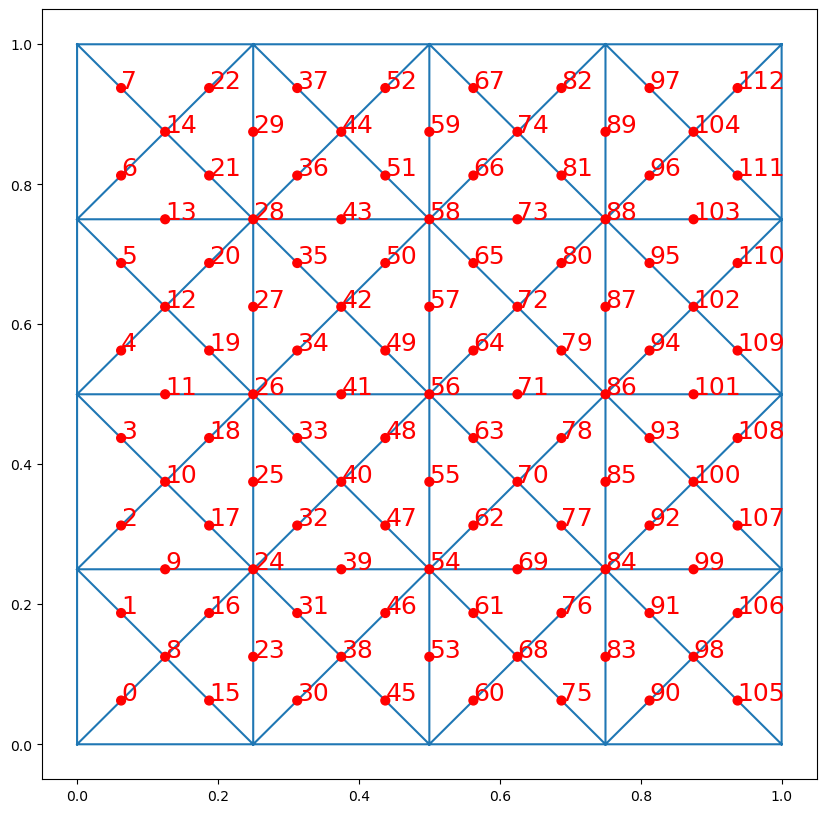} 
	\includegraphics[width=0.32\textwidth]{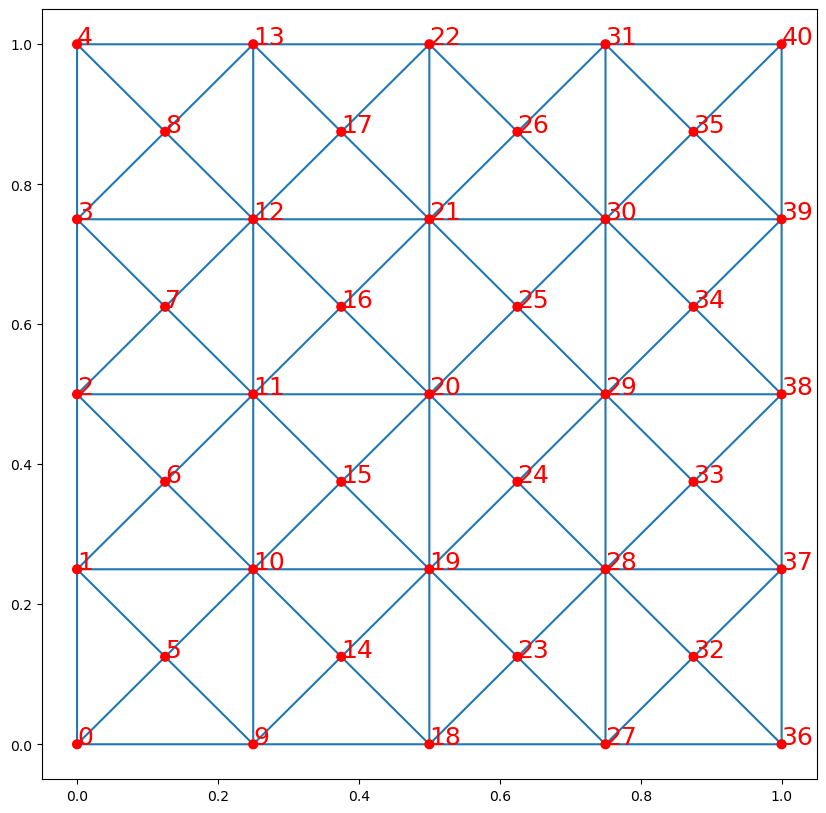} 
	\caption{Ordering for the cells/viscosity (left), 
		$\mathbb{P}_2$ degrees of freedom for velocity components (center) and $\mathbb{P}_1$ pressures (right) for case 
		$n=4$. Note that only the interior velocity degrees of freedom, i.e. the basis functions in $(H^1_0(\Omega))^2$, enter into the linear system.}
	\label{fig:ordering}
\end{figure}
In this section we apply the GLT methodology to matrix sequences stemming from the Stokes problem \eqref{eq:stokesblocks}, when employing the discretization by Taylor-Hood $\mathbb{P}_2$-$\mathbb{P}_1$ elements. The finite element discretization is setup over uniform triangulation of domain $\Omega=(0, 1)^2$ formed by isosceles triangles, see Figure \ref{fig:ordering}. In the following, we assume ordering of the unknowns by the lexicographic ordering. In Appendix \ref{app:A} and \ref{app:B} we report the matrix entries of both $A_n$ and $B_n^T$.

\subsection{The block multilevel structure of \texorpdfstring{$A_n$}{A\_n}}\label{sec:Ablock-structure}
Firstly, we consider the matrix sequence \(\bigl\{A_n\bigr\}_n\), where each \(A_n\) is a diagonal block matrix of the form
\begin{equation}
    A_n \;=\; \begin{bmatrix}
	A_{x,n} & 0 \\
	0         & A_{y,n}
    \end{bmatrix},
\end{equation}
with blocks \(A_{x,n}\) and \(A_{y,n}\) each of dimension \(8n^2 - O(n)\). Since $A_{x,n} = A_{y,n}$, we focus our attention on \(A_{x,n}\).\label{subsec:regular-structure-axx}\\
By construction, the symmetric matrix \(A_{x,n}\) has the following block diagonal components:
\begin{equation}
A_{x,n} \;=\;
\left[
\begin{array}{ccc|c|ccc|c|c}
	D_1 & B_1^T & 0 & 0 & 0 & 0 & 0	& 0 & \\
	B_1 & D_2 & B_2^T & C_1^T & 0 & 0 & 0 & 0 & \\
	0 & B_2 & D_3 & B_3^T & F_1^T & 0 & 0 & 0 & \\
	\hline
	0 & C_1 & B_3 & D_4 & B_4^T & C_2^T & 0 & 0 & \\
	\hline
	0 & 0 & F_1 & B_4 & D_5 & B_5^T & 0 & 0 & \\
	0 & 0 & 0 & C_2 & B_5 & D_6 & B_6^T & C_3^T & \\
	0 & 0 & 0 & 0 & 0 & B_6 & D_7 & B_7^T & \\
	\hline
	0 & 0 & 0 & 0 & 0 & C_3 & B_7 & D_8 & \\
	\hline
	&&&&&&&& \ddots \\
\end{array}
\right],
\end{equation}
where
\begin{itemize}
	\item $D_i$ are square matrices of sizes
	\begin{itemize}
		\item $2(n-2)\times 2(n-2)$ if $i \equiv 0 \pmod{4}$,
		\item $2n\times 2n$ otherwise.
	\end{itemize}
	\item $B_i$ are matrices of sizes
	\begin{itemize}
		\item $(2n-2)\times 2n$ if $i\equiv 3 \pmod 4$,
		\item $2n\times (2n-2)$ if $i\equiv 0 \pmod4$,
		\item $2n\times 2n$ otherwise.
	\end{itemize}
	\item $C_i$ are rectangular matrices of size
	\begin{itemize}
		\item $(2n-2)\times 2n$ if $i$ is odd,
		\item $2n\times (2n-2)$ if $i$ is even.
	\end{itemize}
    \item $F_i$ are square matrices of size $2n \times 2n$ and rank 1.
\end{itemize}
Indeed the collection of $\{F_i\}$ represents a low-rank perturbation, since they appear $O(n)$-times in the block structure. Therefore, we can split $A_{x,n}$ in
\begin{equation}
	A_{x,n} = \hat{A}_{x,n} + E_n,
\end{equation}
where $E_n$ is the matrix that gathers all the perturbations $F_i$. Thus
\begin{equation}
    \lim_{x\rightarrow +\infty} \frac{\operatorname{rk}(E_n)}{n^2} = 0.
\end{equation}
Since \(\{E_{n}\}_n\) is zero-distributed and symmetric, it does not affect the limiting spectral distribution by Theorem \ref{thm:acs_equiv_distributions} and in this section we will focus on the more regular matrix sequence $\{\hat{A}_n\}_n$.

\subparagraph{Simplifying the block structure by using projections}
The structure of the matrix is regular and periodic, but we see irregularities in block sizes.
To simplify the analysis, we format each block as \(2n \times 2n\), adding a minimal number of rows and columns proportional to the overall dimensions of the matrix. Specifically, we insert a row at all the positions \( i \equiv 6n \pmod{8n}\) and \(  j \equiv 0 \pmod {8n} \), and symmetrically add columns at the corresponding positions. Thus, we can consider this new matrix $\hat{A}_{x,n}$ structured with an $8n$-block structure
\begin{equation}
	\hat{A}_{x,n}=\text{tridiag}(\mathcal{A}_1,\mathcal{A}_0,\mathcal{A}_1^T)
\end{equation}
where
\begin{equation*}
	\mathcal{A}_0 = \begin{bmatrix}
		D_1 	& B_1^T & 0 	& 0 \\
		B_1 	& D_2 	& B_2^T & C_1^T	 \\
		0 		& B_2 	& D_3 	& B_3^T \\
		0 		& C_1 	& B_3 	& D_4 	\\
	\end{bmatrix},
    \qquad
	\mathcal{A}_1= \begin{bmatrix}
		0		& 0 	& 0		& B_4 \\
		0		& 0		& 0		& C_2 \\
		0		& 0		& 0		& 0	\\
		0		& 0		& 0		& 0	\\
	\end{bmatrix}.
\end{equation*}
We can then arrange these blocks into an \((8n^2 - 2n)\times (8n^2 - 2n)\) matrix. By suitably adding (or removing) a small number of rows and columns (on the order of $n$), we preserve the overall structure. Hence, from Theorem~\ref{thm:hermitian-compression} and the Definition \ref{def_not1}, the sequence \(\{{A}_n\}_n\) can be embedded into a larger sequence \(\{\widetilde{A}_n\}_n\), which is a block two-level Toeplitz matrix of dimension $8n^2$. This embedding does not alter the spectral distribution. Ultimately, after a permutation similarity transformation, we obtain a GLT matrix sequence whose symbol is an $8\times 8$ Hermitian matrix-valued function, which we describe in detail. Using Theorem \ref{szego-herm}, it follows that its eigenvalues distribute according to this matrix-valued function.\\
The final goal is to show that, up to a similarity transformation by permutations, the enlarged sequence is a block two-level GLT sequence generated by a Hermitian matrix-valued function (in fact a $8\times 8$ Hermitian matrix-valued function, whose entries are trigonometric polynomials times a viscosity function), as in Definition \ref{def_not3}, so that its order is $N(n)=8n_1n_2=8n^2$. In this way, the distribution Theorem \ref{szego-herm} can be invoked and we obtain that the GLT and spectral symbol is the function $\hat{G}(\theta_1,\theta_2)$ in (\ref{mu m hat g}) and subsequent expression. Summarising we derive the spectral analysis of the larger sequence $\{ \widetilde{A}_n \}$ through compressions described just before represented by a semi-orthogonal compression denoted by $\{P_n\}_n$.\\

Now, we want to show the structure of each block of $A_{x,n}$. In our problem, each block of $A_{x,n}$ can be viewed as
\begin{equation}\label{eq:hadamard}
	X \;=\; T \;\diamond\; \mathbb{M},
\end{equation}
where \(T\) is a block matrix with Toeplitz blocks, and \(\mathbb{M}\) carries the discretized viscosity according to our discretization method.
Here, ``\(\diamond\)'' denotes the Hadamard (elementwise) product:
\begin{displaymath}
	(X \diamond Y)_{ij} \;=\; X_{ij}\,Y_{ij}.
\end{displaymath}
The samples in $\mathbb{M}$ appear as sums of some adjacent samples $\mu_i$, as defined in the discretization. We see four types of these coefficients: singletons, the sum of two elements, the sum of four elements and the sum of eight elements.\\
In the following, according to the extradimensional idea, for studying the distributions via the multilevel block Toeplitz structure, we eliminate few rows or columns in a number which do not impact on the resulting distributions, as stated in Theorem \ref{thm:compression-equivalence} and Theorem \ref{thm:hermitian-compression}. More specifically, for $i=1,\ldots,4$, the matrices $D_i$, $B_i$ are replaced by $\hat{D}_i$, $\hat{B}_i$, respectively, while for $i=1,2$, the matrix $C_i$ is replaced by $\hat{C}_i$.

\subparagraph{Diagonal blocks \texorpdfstring{$D_i$}{Di}.}\label{subsubsec:diagonal-blocks}
For \(i=1,\ldots,4\), each diagonal block \(\hat{D}_i\) has one of the following forms:
\begin{equation*}
	\hat{D}_i \;=\;
	\begin{cases}
		\mathrm{diag}\left(\tfrac{8}{3}\right),
		& \text{if } i \equiv 1 \text{ or } 3 \pmod{4},\\[6pt]
		\mathrm{diag}\!\left(
		\begin{bmatrix}
			1 & 0\\
			0 & \tfrac{8}{3}
		\end{bmatrix}
		\right),
		& \text{if } i \equiv 2 \pmod{4},\\[13pt]
		\mathrm{diag}\!\left(
		\begin{bmatrix}
			\tfrac{8}{3} & 0\\
			0 & \tfrac12
		\end{bmatrix}
		\right),
		& \text{if } i \equiv 0 \pmod{4}.
	\end{cases}
\end{equation*}
In addition, each \(\tfrac{8}{3}\) is multiplied by the sum of 2-samples from \(\mathbb{M}\), each \(1\) by the sum of 4-samples, and each \(\tfrac12\) by the sum of 8-samples.

\subparagraph{Off-diagonal blocks \texorpdfstring{$B_i$}{B\_i}.}\label{subsubsec:off-diagonal-blocks-B}
For \(i=1,\ldots,4\), the blocks \(\hat{B}_i\) appear along the first subdiagonal and are defined by
\begin{equation*}
	\hat{B}_i \;=\;
	\begin{cases}
		\text{tridiag}\!\left(
		0,\,
		\begin{bmatrix}
			-\tfrac{2}{3} & -\tfrac{2}{3}\\[2pt]
			0 & -\tfrac{4}{3}
		\end{bmatrix},\,
		\begin{bmatrix}
			0 & 0\\[2pt]
			-\tfrac{4}{3} & 0
		\end{bmatrix}
		\right),
		& \text{if } i \equiv 1 \pmod{4},\\[12pt]
		\text{tridiag}\!\left(
		0,\,
		\begin{bmatrix}
			-\tfrac{2}{3} & -\tfrac{4}{3}\\[2pt]
			0 & -\tfrac{4}{3}
		\end{bmatrix},\,
		\begin{bmatrix}
			0 & 0\\[2pt]
			-\tfrac{2}{3} & 0
		\end{bmatrix}
		\right),
		& \text{if } i \equiv 2 \pmod{4},\\[12pt]
		\text{tridiag}\!\left(
		\begin{bmatrix}
			0 & -\tfrac{4}{3}\\[2pt]
			0 & 0
		\end{bmatrix},\,
		\begin{bmatrix}
			-\tfrac{4}{3} & 0\\[2pt]
			-\tfrac{2}{3} & -\tfrac{2}{3}
		\end{bmatrix},\,
		0
		\right),
		& \text{if } i \equiv 3 \pmod{4},\\[12pt]
		\text{tridiag}\!\left(
		\begin{bmatrix}
			0 & -\tfrac{2}{3}\\[2pt]
			0 & 0
		\end{bmatrix},\,
		\begin{bmatrix}
			-\tfrac{4}{3} & 0\\[2pt]
			-\tfrac{4}{3} & -\tfrac{2}{3}
		\end{bmatrix},\,
		0
		\right),
		& \text{if } i \equiv 0 \pmod{4}.
	\end{cases}
\end{equation*}
Every occurrence of \(-\tfrac{2}{3}\) is multiplied by a sum of 2-samples from \(\mathbb{M}\), and every occurrence of \(-\tfrac{4}{3}\) by a single sample.

\subparagraph{Blocks \texorpdfstring{$C_i$}{C\_i}.}\label{subsubsec:blocks-C}
For \(i=1,2\), the blocks \(\hat{C}_i\) have a single pattern,
\begin{equation*}
	\hat{C}_i \;=\;
	\begin{cases}
		\text{tridiag}\!\left(
		0,\,
		\begin{bmatrix}
			0 & 0\\[2pt]
			\tfrac{1}{6} & 0
		\end{bmatrix},\,
		\begin{bmatrix}
			0 & 0\\[2pt]
			\tfrac{1}{6} & 0
		\end{bmatrix}
		\right),
		& \text{if $i$ is odd,}\\[12pt]
		\text{tridiag}\!\left(
		\begin{bmatrix}
			0 & \tfrac{1}{6}\\[2pt]
			0 & 0
		\end{bmatrix},\,
		\begin{bmatrix}
			0 & \tfrac{1}{6}\\[2pt]
			0 & 0
		\end{bmatrix},\,
		0
		\right),
		& \text{if $i$ is even.}
	\end{cases}
\end{equation*}
In all nonzero entries, \(\tfrac{1}{6}\) is multiplied by a sum of 2-samples from \(\mathbb{M}\).\\

We note that, except for $O(n)$ irregularities, each sampling component linked to every Toeplitz block samples the function based on some lexicographic ordering pattern beginning with the $y$ coordinate. We have not analyzed these patterns yet, as we analyze them when the problem is adjusted using a similarity operator later.

\subsection{Spectral analysis of \texorpdfstring{$A_n$}{A\_n}}
As introduced before, this matrix sequence shows a regular structure but is not of GLT type. Having a GLT structure of the sequence is needed to use important algebraic properties of spectral symbols. We need some similarity transformation to turn the sequence into a GLT. To understand the spectral symbol, we are going to:
\begin{enumerate}
	\item find a suitable {similarity permutation} that makes our problem solvable with the GLT formalism;
	\item study the Toeplitz structure of the permuted coefficient matrix and find its spectral symbol;
	\item study its symbol, which depends on the viscosity function used;
	\item finally, find the final spectral symbol.
\end{enumerate}
To apply the permutation across all blocks, we tensorize once more, resulting in our final similarity operator
\begin{equation}
	\mathbf{\Pi}_{4,n}^{\prime}=I_n \otimes \Pi_{2n,4,2}.
\end{equation}

\subparagraph{Toeplitz symbol}
We start applying the permutation operator to $\hat{A}_{x,n}$
\begin{equation}
	\mathbf{\Pi}_{4,n} \, \hat{A}_{x,n} \, \mathbf{\Pi}_{4,n}^*=\text{tridiag}(T_n(g_1),T_n(g_0),T_n(g_1)^T) \diamond \mathbb{M}.
\end{equation}
Now each block is a $8\times 8$ block 1-level Toeplitz matrix with symbols
\begin{displaymath}
	\resizebox{\linewidth}{!}{%
	$g_0(\theta_2) = \left[
	\begin{array}{cc|cc|cc|cc}
		\frac{8}{3} & 0 & -\frac{2}{3} & -\frac{4}{3}e^{i\theta_2} & 0 & 0 & 0 & 0 \\
		0 & \frac{8}{3} & -\frac{2}{3} & -\frac{4}{3} & 0 & 0 & 0 & 0 \\
		\hline
		-\frac{2}{3} & -\frac{2}{3} & 1 & 0 & -\frac{2}{3} & -\frac{2}{3}e^{i\theta_2} & 0 & h_2(\theta_2)\\
		-\frac{4}{3}e^{-i\theta_2} & -\frac{4}{3} & 0 & \frac{8}{3} & -\frac{4}{3} & -\frac{4}{3} & 0 & 0 \\
		\hline
		0 & 0 & -\frac{2}{3} & -\frac{4}{3} & \frac{8}{3} & 0 & -\frac{4}{3} & -\frac{2}{3} \\
		0 & 0	& -\frac{2}{3}e^{-i\theta_2} & -\frac{4}{3} & 0 & \frac{8}{3} & -\frac{4}{3}e^{-i\theta_2} & -\frac{2}{3} \\
		\hline
		0 & 0 & 0 & 0 & -\frac{4}{3} & -\frac{4}{3}e^{i\theta_2} & \frac{8}{3} & 0 \\
		0 & 0 & \overline{h_2(\theta_2)} & 0 & -\frac{2}{3}	& -\frac{2}{3} & 0 & \frac{1}{2} \\
	\end{array}
	\right],$}
\end{displaymath}
where $h_2(\theta_2)=\frac{1}{6}(1+e^{i\theta_2})$ and
\begin{displaymath}
	{g_1}(\theta_2) = \left[
	\begin{array}{cc|cc|cc|cc}
		0 & 0 & 0 & 0 & 0 & 0 & -\frac{4}{3} & -\frac{2}{3}e^{i\theta_2} \\
		0 & 0 & 0 & 0 & 0 & 0 & -\frac{2}{3} & -\frac{2}{3} \\
		\hline
		0 & 0 & 0 & 0 & 0 & 0 & 0 & h_2(\theta_2) \\
		0 & 0 & 0 & 0 & 0 & 0 & 0 & 0 \\
		\hline
		0 & 0 & 0 & 0 & 0 & 0 & 0 & 0 \\
		0 & 0 & 0 & 0 & 0 & 0 & 0 & 0 \\
		\hline
		0 & 0 & 0 & 0 & 0 & 0 & 0 & 0 \\
		0 & 0 & 0 & 0 & 0 & 0 & 0 & 0 \\
	\end{array}
	\right].
\end{displaymath}
The global matrix $\hat{T}$ is thus a $8\times 8$ block $2$-level Toeplitz matrix with symbol $G(\theta_1,\theta_2)$ equal to
\begin{equation}\label{muConst_symbol}
	\resizebox{\linewidth}{!}{%
	$\hat{G}(\theta_1,\theta_2)= \left[
	\begin{array}{cc|cc|cc|cc}
		\frac{8}{3} & 0 & -\frac{2}{3} & -\frac{4}{3}e^{i\theta_2} & 0 & 0 & -\frac{4}{3}e^{i\theta_1} & -\frac{2}{3}e^{i(\theta_2+\theta_1)} \\
		0 & \frac{8}{3} & -\frac{2}{3} & -\frac{4}{3} & 0 & 0 & -\frac{4}{3}e^{i \theta_1} & -\frac{2}{3}e^{i \theta_1} \\
		\hline
		-\frac{2}{3} & -\frac{2}{3} & 1 & 0 & -\frac{2}{3} & -\frac{2}{3}e^{i\theta_2} & 0 & h_3 \\
		-\frac{4}{3}e^{-i\theta_2} & -\frac{4}{3} & 0 & \frac{8}{3}& -\frac{4}{3} & -\frac{4}{3} & 0 & 0 \\
		\hline
		0 & 0 & -\frac{2}{3} & -\frac{4}{3} & \frac{8}{3}  & 0 & -\frac{4}{3} & -\frac{2}{3} \\
		0 & 0 & -\frac{2}{3}e^{-i\theta_2} & -\frac{4}{3} & 0 & \frac{8}{3} & -\frac{4}{3}e^{-i\theta_2} & -\frac{2}{3} \\
		\hline
		-\frac{4}{3}e^{-i\theta_1} & -\frac{4}{3}e^{-i \theta_1} & 0 & 0 & -\frac{4}{3} & -\frac{4}{3}e^{i\theta_2} & \frac{8}{3} & 0 \\
		-\frac{2}{3}e^{-i(\theta_2+\theta_1)} & -\frac{2}{3}e^{-i \theta_1} & \overline{h_3} & 0 & -\frac{2}{3} & -\frac{2}{3} & 0 & \frac{1}{2} \\
	\end{array}
	\right],$}
\end{equation}
\hspace{1cm}\\
where $h_3(\theta_1,\theta_2)=\frac{1}{6}(1+e^{i\theta_1})(1+e^{i\theta_2})$. In the symbol, we observe two Fourier variables, $\theta_1$ and $\theta_2$, which indicate its $2$-level nature. Additionally, the fact that the symbol is an $8\times 8$ matrix-valued can be understood from the overall matrix's structure, which consists of basic blocks of size $8\times 8$, all in accordance with Definition \ref{def_not3}.

\subparagraph{Analysis of the GLT symbol.} We now examine the spectral symbol that emerges in the distribution of $\left\{A_{x,n}\right\}_n$. To derive the Toeplitz structure from the stiffness matrix, we initially represent $\left\{A_{x,n}\right\}_n$ as a Hadamard product of spatial viscosity coefficients and the Toeplitz matrix.\\
For algebraic purpose, it is more useful to write matrices as
\begin{equation}
	A_{x,\textsc{glt},n} = D_n T_n
\end{equation}
where $D_n$ is some uniform sampling diagonal matrix and $T_n$ is a Toeplitz matrix.\\
Let $\mu(x): [0,1]^2 \to \mathbb{R}$ denote the viscosity function for our scenario. The enumeration of $\mu$ samples is a hard task, but we can handle this by following an approximation approach. The function \(\mu\) is Riemann-integrable and, in particular, belongs to $L^1(\mathbb{R}^2)$. Given the density of continuous functions in $L^1(\mathbb{R}^2)$, we consider a sequence of continuous viscosity functions \(\mu^{[m]}(x,y)\) that converge to \(\mu(x,y)\) in $L^1(\mathbb{R}^2)$. $L^1(\mathbb{R}^2)$ convergence implies convergence in measure, allowing the use of axiom {\bf GLT 6}.

\subparagraph{Approximating sampling structures.}
After the permutation operation, the sample coefficients for each block follow a lexicographic order. Using the continuity of \(\mu^{[m]}(x,y)\), we approximate each sum by multiples of \(\mu^{[m]}(x_i,y_i)\) at specific nodes $(x_i,y_i)$ on an uniform grid of $[0,1]^2$. Up to an approximation with a small error converging to zero, we can approximate the samplings of $\mu_i^{[m]}$ arranged along one of the eight uniform grids in lexicographic coordinate order for each block element. Consequently, the samples now have a simpler pattern for GLT purposes. By adjusting each Toeplitz matrix through redistributing every non-unit coefficient from the sampling matrix arising from the preceding sampling sum ($2$-sum, $4$-sum, etc.) to the relative Toeplitz element, we can express  $\widetilde{A}^{[m]}_{x,\textsc{glt},n}$ from every non-trivial diagonal of this matrix as
\begin{equation}
	\widetilde{A}^{[m]}_{x,\textsc{glt},n} = \sum_{i} D_{i,n}^{(m)} T_{i,n}{,}
\end{equation}
where $D_{i,n}^{(m)}$ is a multilevel uniform diagonal sampling matrix and $T_{i,n}$ is the Toeplitz matrix that represents the specific diagonal. 
The symbol for each \(\{D^{(m)}_{i,n}\}_n\) is expressed as $\mu^{[m]}(x,y)$ thanks to the approximation considered. The reason is standard in the a.c.s convergence of GLT sequences. Indeed, let $\widetilde{A}_{x,\textsc{glt},n}:=\sum_i D_{i,n}T_{i,n}$ denote the analogous matrix obtained from the exact viscosity sampling. Since $\mu^{[m]}\to\mu$ in measure, we can choose $\varepsilon_m\downarrow0$ so that 
\begin{equation}
	E_m:=\{(x,y)\in[0,1]^2:\ |\mu^{[m]}(x,y)-\mu(x,y)|>\varepsilon_m\}
\end{equation}
satisfies $|E_m|\to0$ and $\partial E_m$ has zero measure. Splitting the sampled nodes into $E_m$ and its complement gives $D_{i,n}-D_{i,n}^{(m)}=R_{i,n,m}+N_{i,n,m}$, with $\|N_{i,n,m}\|\le\varepsilon_m$ and $R_{i,n,m}$ diagonal, supported only on the sampled nodes belonging to $E_m$. Because each block uses a uniform grid (up to the finitely many shifts already described), 
\begin{equation}
	\operatorname{rk}(R_{i,n,m})\le (|E_m|+o(1))d_n;
\end{equation}
hence $\{D_{i,n}^{(m)}\}_n\xrightarrow{\mathrm{a.c.s.}}\{D_{i,n}\}_n$ for every $i$. By Theorem~\ref{thm:acs-properties}, it follows that
\begin{equation}
	\{\widetilde{A}^{[m]}_{x,\textsc{glt},n}\}_n\xrightarrow{\mathrm{a.c.s.}}\{\widetilde{A}_{x,\textsc{glt},n}\}_n.
\end{equation}
Since for each fixed $m$ axioms {\bf GLT 2}, {\bf GLT 3}, and {\bf GLT 5} give the symbol $\mu^{[m]}(x,y)\hat{G}(\theta_1,\theta_2)$, axiom {\bf GLT 6} yields the limit symbol $\mu(x,y)\hat{G}(\theta_1,\theta_2)$.
Using the block multilevel GLT formalism, i.e. axiom {\bf GLT 1} -- axiom {\bf GLT 5}, combining all elements $T_{i,n}$ and using the fact that the GLT class forms an algebra, we derive the spectral symbol for the entire approximated matrix sequence. More specifically, the GLT and spectral symbol is
\begin{equation}\label{mu m hat g}
	f^{[m]}(x,y,\theta_1,\theta_2) \; = \;  \mu^{[m]}(x,y) \, \hat{G}(\theta_1,\theta_2)
\end{equation}
with
\begin{equation}\label{muVar_symbol}
	\resizebox{\linewidth}{!}{%
	$\hat{G}(\theta_1,\theta_2) =
	\left[
	\begin{array}{cc|cc|cc|cc}
		\frac{16}{3} & 0 & -\frac{4}{3} & -\frac{4}{3}e^{i\theta_2} & 0 & 0 & -\frac{4}{3}e^{i\theta_1} & -\frac{4}{3}e^{i(\theta_2+\theta_1)} \\
		0 & \frac{16}{3} & -\frac{4}{3} & -\frac{4}{3} & 0 & 0 & -\frac{4}{3}e^{i \theta_1} & -\frac{4}{3}e^{i \theta_1} \\
		\hline
		-\frac{4}{3} & -\frac{4}{3} & 4 & 0 & -\frac{4}{3} & -\frac{4}{3}e^{i\theta_2} & 0 & 2h_3 \\
		-\frac{4}{3}e^{-i\theta_2} & -\frac{4}{3} & 0 & \frac{16}{3} & -\frac{4}{3} & -\frac{4}{3} & 0 & 0 \\
		\hline
		0 & 0 & -\frac{4}{3} & -\frac{4}{3} & \frac{16}{3} & 0 & -\frac{4}{3} & -\frac{4}{3} \\
		0 & 0 & -\frac{4}{3}e^{-i\theta_2} & -\frac{4}{3} & 0 & \frac{16}{3} & -\frac{4}{3}e^{-i\theta_2} & -\frac{4}{3} \\
		\hline
		-\frac{4}{3}e^{-i\theta_1} & -\frac{4}{3}e^{-i \theta_1} & 0 & 0 & -\frac{4}{3} & -\frac{4}{3}e^{i\theta_2} & \frac{16}{3} & 0 \\
		-\frac{4}{3}e^{-i(\theta_2+\theta_1)} & -\frac{4}{3}e^{-i \theta_1} & 2\overline{h_3} & 0 & -\frac{4}{3} & -\frac{4}{3} & 0 & 4 \\
	\end{array}
	\right],$}
\end{equation}
where $h_3(\theta_1,\theta_2)=\frac{1}{6}(1+e^{i\theta_1})(1+e^{i\theta_2})$.\\
By passing to the a.c.s. limit, we derive the final spectral distribution expressed by the expected $4$-variate $8\times8$ matrix-valued function
\begin{equation}
	f(x,y,\theta_1,\theta_2) \; = \;  \mu(x,y) \, \hat{G} \, (\theta_1,\theta_2).
\end{equation}
Here the variables $(x,y)$ are the physical variables of the variable-coefficient PDE, while the variables $(\theta_1,\theta_2)$ are the Fourier variables of the operators and the $8\times 8$ dimensionality of the symbol refers to the chosen Taylor-Hood approximation. In other terms, a change in the approximation of the given PDE would affect only the dimensionality of the symbol and not the number of variables, which are intrinsic to the considered differential operator.\\
Finally, let us consider the matrix
\begin{equation}
	\mathbf{\Pi}_{4,n}\widetilde{A}_n\mathbf{\Pi}_{4,n}^{*} \; = \;
	\begin{bmatrix}
		\widetilde{A}_{x,\textsc{glt},n} & 0 \\
		0 & \widetilde{A}_{y,\textsc{glt},n}
	\end{bmatrix}
\end{equation}
with
\begin{equation}
    \mathbf{\Pi}_{4,n}=\text{diag}\bigl(\mathbf{\Pi}_{4,n}^{\prime},\mathbf{\Pi}_{4,n}^{\prime}\bigr).
\end{equation}
Since $\hat{A}_{x,n}=\hat{A}_{y,n}$, the spectral distribution of $\{\widetilde{A}_n\}_n$ is the same as the matrix sequence of its nonzero blocks, and by using Theorem \ref{thm:permutation_rect_eq} we find a permutation $\Gamma_1$ such that
\begin{equation}
	\bigl\{  \Gamma_1 \mathbf{\Pi}_{4,n}\widetilde{A}_n\mathbf{\Pi}_{4,n}^* \Gamma_1^* \bigr\}_n\; = \;  \bigl\{ A_{n,\textsc{glt}} \bigr\}_n \; \sim_{\textsc{glt}} \; F,
\end{equation}
with
\begin{equation}
	F(x,y,\theta_1,\theta_2)\; = \; \begin{bmatrix}
		f(x,y,\theta_1,\theta_2) & 0 \\
		0 & f(x,y,\theta_1,\theta_2)
	\end{bmatrix}.
\end{equation}
In preparation for the subsequent analysis, we briefly summarize the approximation and extension procedures used to derive the spectral distribution of \( A_n \) through the GLT theory. Specifically, the key transformations from $A_{\textsc{glt},n}$ to $A_n$ are:
\begin{enumerate}
	\item from $A_{\textsc{glt},n}$ to $\text{diag}(A_{x,\textsc{glt},n}, A_{y,\textsc{glt}},n)$ through a permutation matrix \(\Gamma_{1,n}\), as detailed in Theorem \ref{thm:permutation_rect_eq};
	\item from $\text{diag}(A_{x,\textsc{glt},n}, A_{y,\textsc{glt},n})$ to the permuted $\widetilde{A}_n$ via the block-diagonal permutation \(\mathbf{\Pi}_{4,n} = \mathrm{diag}(    \mathbf{\Pi}^{\prime}_{4,n}, \mathbf{\Pi}^{\prime}_{4,n})\), as described in Theorem \ref{thm:permutation_rect_eq};
	\item from $\widetilde{A}_n$ to $\hat{A}_n$ through compression by \(P_n\), as discussed in Theorem \ref{thm:hermitian-compression};
	\item from $\hat{A}_n$ to $A_n$ through a small $O(n)$-rank additive symmetric perturbation, i.e. $A_n=\hat{A}_n+E_{n}$.
\end{enumerate}
Thus, the complete transformation connecting \( A_n \) to the GLT formulation \( A_{\textsc{glt},n} \) can be succinctly expressed as:
\begin{equation}\label{A_transf}
	A_n = P_n^*\mathbf{\Pi}_{4,n}^*\Gamma_{1,n}^*A_{\textsc{glt},n}\Gamma_{1,n}\mathbf{\Pi}_{4,n}P_n+E_n.
\end{equation}

\subsection{The block structure of \texorpdfstring{$B_n$}{B\_n}}\label{sec:B}
As in the case of the matrix $\{A_n\}_n$, we begin our analysis by examining the structure of $\{B_n\}_n$. It is helpful to regard $B_n$ as a rectangular block matrix of the form
\begin{equation}
	B_n = \left[ B_{x,n} \quad B_{y,n} \right],
\end{equation}
where $B_{x,n}$ and $B_{y,n}$ are rectangular matrices of dimensions $(8n^2 - O(n)) \times (2n^2 - O(n))$. Unlike the case of $A_n$, it is important to note that $B_{x,n} \neq B_{y,n}$ and both are independent of the viscosity function $\mu(x,y)$. Consequently, using the same approximation process used for $A_n$, we will show that the singular value symbol of the resulting matrix sequences resembles a bivariate function depending only on the frequency variables $f(\theta_1,\theta_2)$ with $\theta_i \in [-\pi,\pi]$. Without loss of generality, we perform the singular value analysis on the transpose $B_n^T$.\\
Again, from the construction of $B_{n}$ one can identify a periodic block pattern along the diagonal. We have the following block structure:
\begin{displaymath}
	B_{x,n} = \left[
	\begin{array}{cc|cc|cc|c}
		L_1 & -L_4	& -L_3 	& 0 	& 0		& 0		&\\
		L_2	& 0		& -L_2 	& 0		& 0 	& 0		&\\
		L_3	& L_4 	& -L_1 	& 0 	& 0 	& 0		&\\
		\hline
		0 	& J 	& 0 	& -J 	& 0		& 0		&\\
		\hline
		0	& 0 	& L_1	& -L_4	& -L_3	& 0		&\\
		0	& 0		& L_2	& 0		& -L_2	& 0		&\\
		0	& 0		& L_3	& L_4	& -L_1	& 0		&\\
		\hline
		0 	& 0 	& 0 	& J	& 0		& -J	&\\
		\hline
		&&&&&&\ddots
	\end{array}
	\right].
\end{displaymath}
with the respective dimensions:
\begin{itemize}
	\item $L_1$: $2n \times (n+1)$;
	\item $L_2$: $(2n-1) \times (n+1)$;
	\item $L_3$: $2n \times (n+1)$;
	\item $L_4$: $2n \times n$;
	\item $J$: $(2n-1) \times n$.
\end{itemize}
Similarly, $B_y$ has a corresponding regular structure:
\begin{displaymath}
	B_{y,n} = \left[
	\begin{array}{cc|cc|cc|c}
		H_1 & H_2	& H_3 	& 0 	& 0		& 0	&\\
		0	& H_4	& 0 	& 0		& 0 	& 0	&\\
		H_3	& H_2 	& H_1 	& 0 	& 0 	& 0	&\\
		\hline
		0 	& 0 	& K 	& 0 	& 0		& 0	&\\
		\hline
		0	& 0 	& H_1	& H_2	& H_3	& 0	&\\
		0	& 0		& 0		& H_4	& 0		& 0	&\\
		0	& 0		& H_3	& H_2	& H_1	& 0	&\\
		\hline
		0 	& 0 	& 0 	& 0		& K		& 0	&\\
		\hline
		&&&&&&\ddots
	\end{array}
	\right],
\end{displaymath}
and the dimensions of each block are:
\begin{itemize}
	\item $H_1$: $2n \times (n+1)$;
	\item $H_2$: $2n \times n$;
	\item $H_3$: $2n \times (n+1)$;
	\item $H_4$: $(2n-1) \times n$;
	\item $K$: $(2n-1) \times (n+1)$.
\end{itemize}
The periodic shape of $B_{x,n}$ has the following components:
\begin{equation}
	\begin{split}
		L_1 &= \text{tridiag}\left(
		0,\
		\begin{bmatrix}
			-1/6\\
			-1/12
		\end{bmatrix},
		\begin{bmatrix}
			-1/12\\
			-1/6
		\end{bmatrix}
		\right),
		\quad
		L_2 = \text{tridiag}\left(
		0,\
		0,
		\begin{bmatrix}
			0\\
			-1/6
		\end{bmatrix}
		\right),\\
		L_3 &= \text{tridiag}\left(
		0,\
		\begin{bmatrix}
			-1/12\\
			-0
		\end{bmatrix},
		\begin{bmatrix}
			-0\\
			-1/12
		\end{bmatrix}
		\right),
		\quad
		L_4 = \text{diag}\left(
		\begin{bmatrix}
			-1/6\\
			-1/6
		\end{bmatrix}
		\right),\\
		J &= \text{diag}\left(
		\begin{bmatrix}
			-1/6\\
			-0
		\end{bmatrix}
		\right).
	\end{split}
\end{equation}
Recall that $L_2$ and $J$ are defined with one fewer row than $2n$ to match the required dimensions. Such dimension discrepancies will be fixed later.\\
On the other hand, the periodic shape of $B_{y,n}$ is composed of:
\begin{equation}
	\begin{split}
		H_1 &= \text{tridiag}\left(
		0,\
		\begin{bmatrix}
			-1/6\\
			-1/12
		\end{bmatrix},
		\begin{bmatrix}
			1/12\\
			1/6
		\end{bmatrix}
		\right),		
		\quad
		H_2 = \text{diag}\left(
		\begin{bmatrix}
			1/6\\
			-1/6
		\end{bmatrix}
		\right),\\
		H_3 &= \text{tridiag}\left(
		0,\
		\begin{bmatrix}
			-1/12\\
			0
		\end{bmatrix},
		\begin{bmatrix}
			0\\
			1/12
		\end{bmatrix}
		\right),
		\quad
		H_4 = \text{tridiag}\left(
		0,\
		\begin{bmatrix}
			0\\
			-1/6
		\end{bmatrix},
		\begin{bmatrix}
			0\\
			1/6
		\end{bmatrix}
		\right),\\
		K &= \text{tridiag}\left(
		0,\
		\begin{bmatrix}
			-1/6\\
			0
		\end{bmatrix},
		\begin{bmatrix}
			1/6\\
			0
		\end{bmatrix}
		\right).
	\end{split}
\end{equation}
Similarly, for the matrix $B_{y,n}$, its blocks $H_4$ and $K$ have dimensions that differ slightly from those of a full $2n \times n$ or $2n \times (n+1)$ block. In particular, $H_4$ is of size $(2n-1) \times n$ and $K$ is of size $(2n-1) \times (n+1)$, meaning that both are missing one row. Also note that $H_1$, $H_3$, and $K$ have one more column than $n$.

\subsection{Spectral analysis of \texorpdfstring{$B_n$}{B\_n}}
Like $A_n$, the structure of this matrix is regular and periodic, but we observe certain irregularities in block sizes. To handle these, we again perform the extension-compression trick, considering the enlarged matrices $\widetilde{B}_{x,n}$ and $\widetilde{B}_{y,n}$, obtained by adding
\begin{itemize}
	\item rows at positions \(i \equiv 4n \pmod{8n}\) and \(j \equiv 0 \pmod{8n}\);
	\item columns at positions \(i \equiv {n+1} \pmod{4n}\) and \(j \equiv 3n+1 \pmod{8n}\).
\end{itemize}
Again, by construction, there is an \(O(n)\) truncation in the last row and column of the original $B_{x,n}$ and $B_{y,n}$, preventing them from forming a proper block Toeplitz structure, but enlarging these matrices with a minimal number of additional rows and columns, we achieve a perfectly periodic block structure. By doing so, we ensure that the resulting matrix has a neat block structure. \\
We call $Q_n$ the semi-orthogonal compression that results in the described modification of the rows{,} and $R_n$ the compression to the columns applied to the entire $\widetilde{B}_n$ such that
\begin{equation}
	Q_n^* \, \widetilde{B}_{n} \, R_n =B_{n}{.}
\end{equation}
After this process, the enlarged $\widetilde{B}_n^T$ is a matrix of size $8n^2 \times 4n^2$. Thus, the hypothesis of Theorem \ref{thm:compression-equivalence} are met and spectral analysis can proceed as if the matrix were fully structured and without irregularities. By using Theorem \ref{thm:permutation_rect_eq} we consider matrices $\mathbf{\Pi}_{4,n}^{\prime}$ and $\mathbf{\Pi}_{2,n}^{\prime}$ such that
\begin{equation}
	\mathbf{\Pi}_{4,n}^{\prime,*} \, \widetilde{B}_{x,n} \, \mathbf{\Pi}_{2,n}^{\prime} \; = \; \operatorname{tridiag}(0, T_n(\hat{g}_{x,0}), T_n(\hat{g}_{x,1})^T) = B_{x,\textsc{glt},n},
\end{equation}
with
\begin{equation}
	\resizebox{\linewidth}{!}{%
	$\hat{g}_{x,0}(\theta) \; =  \;
	\begin{bmatrix}
		-\frac{1}{6} & \frac{1}{6} \\
		-\frac{1}{12} & \frac{1}{6} \\
		0 & 0 \\
		0 & 0 \\
		-\frac{1}{12} & -\frac{1}{6} \\
		0 & -\frac{1}{6} \\
		0 & -\frac{1}{6} \\
		0 & 0
	\end{bmatrix}
	+
	\begin{bmatrix}
		-\frac{1}{12} & 0 \\
		-\frac{1}{6} & 0 \\
		0 & 0 \\
		-\frac{1}{6} & 0 \\
		0 & 0 \\
		-\frac{1}{12} & 0 \\
		0 & 0 \\
		0 & 0
	\end{bmatrix} e^{-i\theta},
	\qquad
	\hat{g}_{x,1}(\theta) \; = \;
	\begin{bmatrix}
		\frac{1}{12} & 0 \\
		0 & 0 \\
		0 & 0 \\
		0 & 0 \\
		\frac{1}{6} & 0 \\
		\frac{1}{12} & 0 \\
		0 & \frac{1}{6} \\
		0 & 0
	\end{bmatrix}
	+
	\begin{bmatrix}
		0 & 0 \\
		\frac{1}{12} & 0 \\
		0 & 0 \\
		\frac{1}{6} & 0 \\
		\frac{1}{12} & 0 \\
		\frac{1}{6} & 0 \\
		0 & 0 \\
		0 & 0
	\end{bmatrix} e^{-i\theta}.$}
\end{equation}
Similarly,
\begin{equation}
    \mathbf{\Pi}_{4,n}^{\prime} \widetilde{B}_{y,n} \mathbf{\Pi}_{4,n}^{\prime} = \operatorname{tridiag}(0, T_n(\hat{g}_{y,0}), T_n(\hat{g}_{y,1})^T) = B_{y,\textsc{glt},n},
\end{equation}
with
\begin{equation}
	\resizebox{\linewidth}{!}{%
	$\hat{g}_{y,0}(\theta) =
	\begin{bmatrix}
		-\frac{1}{6} & \frac{1}{6} \\
		-\frac{1}{12} & -\frac{1}{6} \\
		0 & 0 \\
		0 & -\frac{1}{6} \\
		-\frac{1}{12} & \frac{1}{6} \\
		0 & -\frac{1}{6} \\
		0 & 0 \\
		0 & 0
	\end{bmatrix}
	+
	\begin{bmatrix}
		\frac{1}{12} & 0 \\
		\frac{1}{6} & 0 \\
		0 & 0 \\
		0 & \frac{1}{6} \\
		0 & 0 \\
		\frac{1}{12} & 0 \\
		0 & 0 \\
		0 & 0
	\end{bmatrix} e^{-i\theta},
	\qquad
	\hat{g}_{y,1}(\theta) =
	\begin{bmatrix}
		-\frac{1}{12} & 0 \\
		0 & 0 \\
		0 & 0 \\
		0 & 0 \\
		-\frac{1}{6} & 0 \\
		-\frac{1}{12} & 0 \\
		-\frac{1}{6} & 0 \\
		0 & 0
	\end{bmatrix}
	+
	\begin{bmatrix}
		0 & 0 \\
		\frac{1}{12} & 0 \\
		0 & 0 \\
		0 & 0 \\
		\frac{1}{12} & 0 \\
		\frac{1}{6} & 0 \\
		\frac{1}{6} & 0 \\
		0 & 0
	\end{bmatrix} e^{-i\theta}.$}
\end{equation}
Since $\{T_n(\hat{g}_{x,i})\}_n \sim_{\textsc{glt}} \hat{g}_{x,i}$ and $\{T_n(\hat{g}_{y,i})\}_n \sim_{\textsc{glt}} \hat{g}_{y,i}$ for $i \in \{0,1\}$, using similar argument used for $A_n$, we get $\{B_{x,n}\}_n\sim_{\sigma} G_x$ and $\{B_{y,n}\}_n\sim_{\sigma} G_y$, with
\begin{equation}\label{Bxsymbol}
	\resizebox{\linewidth}{!}{%
	$G_x(\theta) =
	\begin{bmatrix}
		-\frac{1}{6} & \frac{1}{6} \\
		-\frac{1}{12} & \frac{1}{6} \\
		0 & 0 \\
		0 & 0 \\
		-\frac{1}{12} & -\frac{1}{6} \\
		0 & -\frac{1}{6} \\
		0 & -\frac{1}{6} \\
		0 & 0
	\end{bmatrix}
	+
	\begin{bmatrix}
		-\frac{1}{12} & 0 \\
		-\frac{1}{6} & 0 \\
		0 & 0 \\
		-\frac{1}{6} & 0 \\
		0 & 0 \\
		-\frac{1}{12} & 0 \\
		0 & 0 \\
		0 & 0
	\end{bmatrix} e^{-i\theta_1}
	+
	\begin{bmatrix}
		\frac{1}{12} & 0 \\
		0 & 0 \\
		0 & 0 \\
		0 & 0 \\
		\frac{1}{6} & 0 \\
		\frac{1}{12} & 0 \\
		0 & \frac{1}{6} \\
		0 & 0
	\end{bmatrix} e^{-i\theta_2}
	+
	\begin{bmatrix}
		0 & 0 \\
		\frac{1}{12} & 0 \\
		0 & 0 \\
		\frac{1}{6} & 0 \\
		\frac{1}{12} & 0 \\
		\frac{1}{6} & 0 \\
		0 & 0 \\
		0 & 0
	\end{bmatrix} e^{-i(\theta_1 + \theta_2)},$}
\end{equation}
and
\begin{equation}\label{Bysymbol}
	\resizebox{\linewidth}{!}{%
	$G_y(\theta) =
	\begin{bmatrix}
		-\frac{1}{6} & \frac{1}{6} \\
		-\frac{1}{12} & -\frac{1}{6} \\
		0 & 0 \\
		0 & -\frac{1}{6} \\
		-\frac{1}{12} & \frac{1}{6} \\
		0 & -\frac{1}{6} \\
		0 & 0 \\
		0 & 0
	\end{bmatrix}
	+ \begin{bmatrix}
		\frac{1}{12} & 0 \\
		\frac{1}{6} & 0 \\
		0 & 0 \\
		0 & \frac{1}{6} \\
		0 & 0 \\
		\frac{1}{12} & 0 \\
		0 & 0 \\
		0 & 0
	\end{bmatrix} e^{-i\theta_1}
	+ \begin{bmatrix}
		-\frac{1}{12} & 0 \\
		0 & 0 \\
		0 & 0 \\
		0 & 0 \\
		-\frac{1}{6} & 0 \\
		-\frac{1}{12} & 0 \\
		-\frac{1}{6} & 0 \\
		0 & 0
	\end{bmatrix} e^{-i\theta_2}
	+ \begin{bmatrix}
		0 & 0 \\
		\frac{1}{12} & 0 \\
		0 & 0 \\
		0 & 0 \\
		\frac{1}{12} & 0 \\
		\frac{1}{6} & 0 \\
		\frac{1}{6} & 0 \\
		0 & 0
	\end{bmatrix} e^{-i(\theta_1 + \theta_2)}.$}
\end{equation}
Combining all the results with Theorem \ref{thm:permutation_rect_eq}, we finally obtain the singular value symbol of $B_n^T$:
\begin{equation}
	G(\theta_1, \theta_2) =
	\begin{bmatrix}
		G_x(\theta_1 , \theta_2) \\
		G_y(\theta_1, \theta_2)
	\end{bmatrix};
\end{equation}
in particular we can find permutations $\Gamma_1$ and $\Gamma_2$ such that
\begin{equation}
    \mathbf{\Pi}_{4,n}^{*} \, \widetilde{B}_{n} \, \mathbf{\Pi}_{2,n}\; = \; \Gamma_1^* \, B_{\textsc{glt},n} \, \Gamma_2,
\end{equation}
with $\bigl\{ B_{\textsc{glt},n} \bigr\}_n \sim_{\textsc{glt}} G$.\\
To conclude this section, we list the process that leads from $B_{\textsc{glt},n}$ to $B_n$:
\begin{enumerate}
	\item from $B_{\textsc{glt},n}$ to $[{B_{x,\textsc{glt},n},B_{y,\textsc{glt},n}}]$ through a permutation matrices $\Gamma_2^*$ on the left, and $\Gamma_1$ on the right, as described in Theorem \ref{thm:permutation_rect_eq};
	\item from $[B_{x,\textsc{glt},n},B_{y,\textsc{glt},n}]$ to the enlarged $\widetilde{B}_n$ through $\mathbf{\Pi_{2,n}}^*$ on the left {and} $\mathbf{\Pi_{4,n}}=\mathrm{diag}({\mathbf{\Pi}_{4,n}^\prime,\mathbf{\Pi}_{4,n}^\prime})$ on the right, as described in Theorem \ref{thm:permutation_rect_eq};
	\item from $\widetilde{B}_n$ to the compressed $B_n$ through the compressions $R_n^*$ on the left {and} $Q_n$ on the right, as described in Theorem \ref{thm:compression-equivalence}.
\end{enumerate}
Thus, the total transformation from the spectral analysis of $B_n$ to a GLT problem is given by this expression:
\begin{equation}\label{B_transf}
    B_n=R_n^*\mathbf{\Pi}_{2,n}^*\Gamma_{2,n}^*B_{\textsc{glt},n}\Gamma_{1,n}\mathbf{\Pi}_{4,n}Q_n.
\end{equation}

\section{Approximation and numerical evidences}\label{sec:precond}
In the current section we consider preliminary proposals of GLT-based approximations. For the saddle-point problems considered here, we remark that the construction of high-quality preconditioners relies on exploiting their block structure, which varies from application to application. The approximation we propose is a combination of the blocking strategy adopted in \cite{Blocking} and the well-known block diagonal preconditioner that uses the Schur complement $S_n$ \cite{benzi2005numerical} as discussed in the introduction.\\
As we know, the exact Schur complement is obtained via algebraic operations on the blocks. Hence, based on the GLT axioms in Section \ref{GLT}, the matrix sequence $\{S_n\}_n$ is expected to distribute as the symbol obtained via the very {same} operations on the symbols of the related blocks. However, compared to previous similar contributions \cite{Schur-CMAME,mazza2021matrix}, we have the additional difficulty of intricate permutation and compression matrices, which were essential for the spectral analysis of the matrix sequences associated with the blocks.\\
As a first preliminary step, the block $A_n$ is replaced by a proper matrix algebra approximation $P_{A,n}$ in such a way that, again by the GLT axioms we have
\begin{equation}
    \{A_n-P_{A,n}\}_n \sim_{\textsc{glt}} 0. 
\end{equation}
Since $A_n$ is a multilevel Toeplitz block matrix, we could apply the strategy in \cite{Blocking}, even if we are aware of the theoretical topological barriers proved in general in \cite{MultiNo,noutsos2004matrix}. In particular the Frobenius optimal circulant preconditioner for a multilevel Toeplitz matrix is not optimal. However, we adopt the blocking preconditioning substituting each Toeplitz block with its $\tau$ approximation. The same theoretical barriers apply also in the $\tau$ setting, but, as proven e.g. in \cite{MultiNo,noutsos2004matrix}, the outlier effect is much less evident in the $\tau$ case, when the order of the zero at $\theta=0$ of the GLT symbol is bounded by $2$. Indeed, let $P_{A,n}$ {be} such a block $\tau$ preconditioner of $A_n$. Then, we construct the block diagonal preconditioner of
\begin{equation}
    P_n \;=\; \begin{bmatrix}
	P_{A,n} & B_n^T \\
	B_n       & 0
    \end{bmatrix}.
\end{equation}
Therefore, if we denote by $\hat{S}_n = -B_n P_{A,n}^{-1} B_n^T$ the Schur complement of $P_n$, the block diagonal preconditioner is given by
\begin{equation}
    S_n \;=\; \begin{bmatrix}
	P_{A,n} & 0 \\
	0       & -\hat{S}_n
    \end{bmatrix}.
\end{equation}

\section{Numerical tests}\label{sec:numerical-tests}
In the current section, we present an unified overview of our numerical experiments aimed at illustrating the distributional behaviour of the various matrix sequences arising from the discretization scheme used in the previous sections. In all experiments regarding the visualizations of the spectral and singular values distributions, we fix $n=16$, i.e. matrix size equal to $4515$, and we use a dense uniform sampling of about $10^5$ points of $[0,1]^2 \times [-\pi,\pi]^2$ to evaluate the matrix-valued symbol. Subsequently, for preconditioned GMRES experiments we employ $n=8, 16, 32$ corresponding to matrix-sizes equal to $1107, 4515, 18243$, respectively.\\
Let us highlight that the visible distributional behaviour already observed for $n=16$ is excellent. In fact, the value $n=16$ is small and, by the theoretical results, the match among symbols and distributional behaviour improves as $n$ becomes larger and larger.\\
We divide the experiments concerning the symbol adherence into three groups of tests:
\paragraph{Group 1} $\mu$ is the constant function 1;

\paragraph{Group 2} $\mu$ is the continuous function
\begin{equation}
	\mu(x,y) = xy + e^{x+y};
\end{equation}

\paragraph{Group 3} $\mu$ is a piecewise continuous function given by
\begin{equation}\label{eq:piecewise_mu_1}
	\mu(x,y) =
	\begin{cases}
		\gamma  & \text{if } (x,y)\in \left[0,\frac{1}{2}\right]^2, \\[3pt]
		1+x+y & \text{otherwise}
	\end{cases}
\end{equation}
where $\gamma\in \left\{1,10,100\right\}$.\\
For the performance of the preconditioner, let us write the system in this form:
\begin{equation}
	M_n x = b\qquad\qquad\text{where}\quad M_n = \begin{bmatrix}
		A_n & B_n^T\\
		B_n & 0
	\end{bmatrix}.
\end{equation}
We divide each of the previous groups into three cases according to the right-hand side:
\paragraph{Case a} $b = 1$, the vector of all ones;

\paragraph{Case b} $b$ is a uniform sampling of $[0,1]^2$;

\paragraph{Case c} $b$ is a random sampling of $[0,1]^2$.\\

For each group and for each case, we report the number of iterations of PGMRES, using $S_n$ as preconditioner, a restart value of 20, and a tolerance of $10^{-5}$. We omit the iterations of standard GMRES, because it requires more than 1000 iterations to reach convergence.
\begin{figure}
	\centering
	\includegraphics[width=0.45\linewidth]{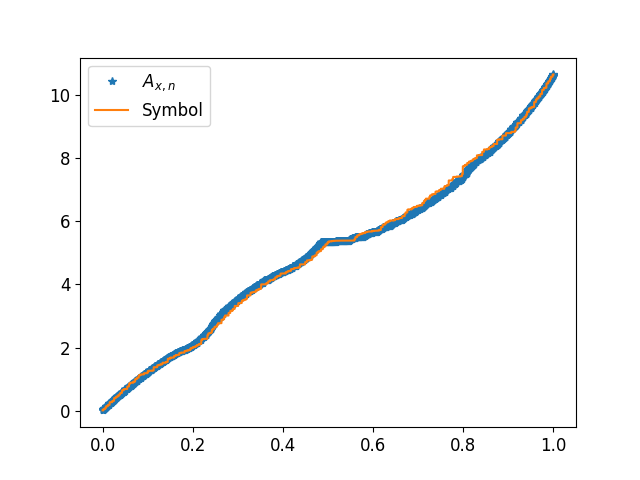}
	\includegraphics[width=0.45\linewidth]{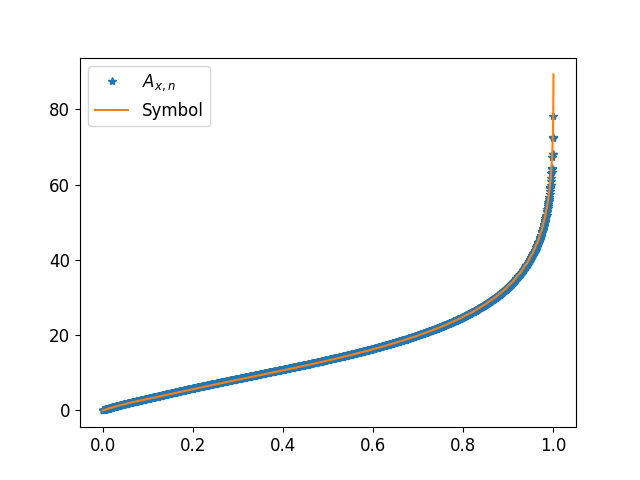}
	\hfill
	\includegraphics[width=0.32\linewidth]{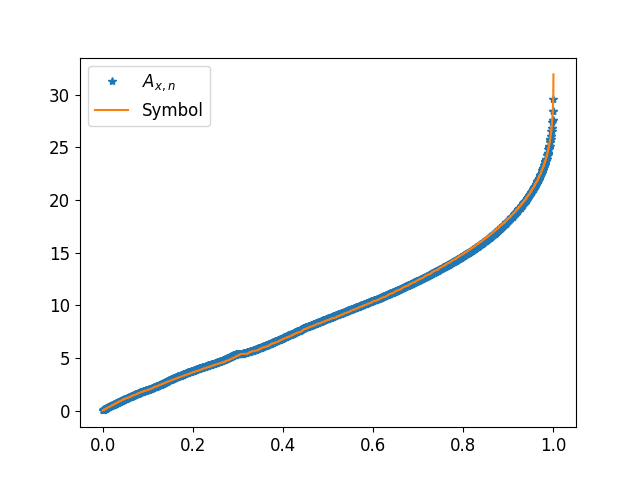}
	\includegraphics[width=0.32\linewidth]{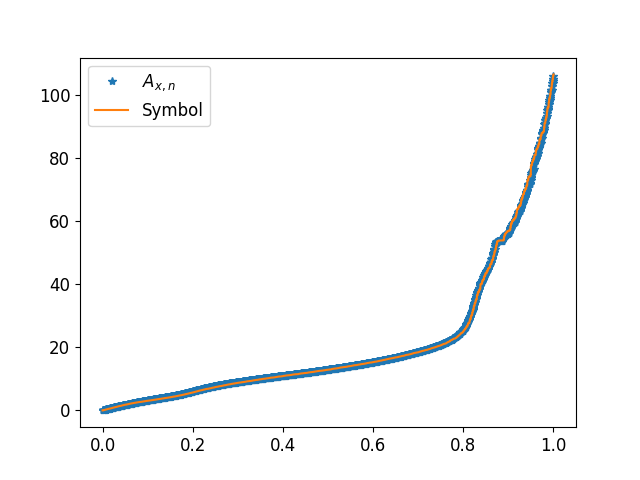}
	\includegraphics[width=0.32\linewidth]{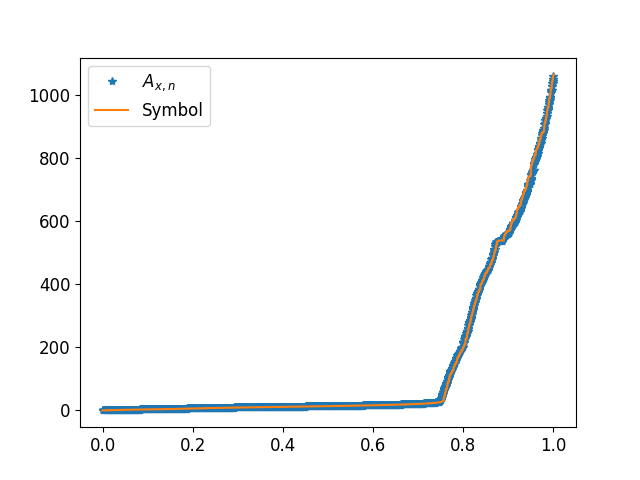}
	\caption{Adherence of the eigenvalues of $A_{x,n}$ to its symbol for $n=16$. The first line displays the plots corresponding to Group 1 (on the left) and Group 2 (on the right), while the second line shows Group 3; from left $\gamma = 1$, $\gamma = 10$, $\gamma = 100$.}
	\label{test_An123}
\end{figure}
\begin{figure}
	\centering
	\includegraphics[width=0.45\linewidth]{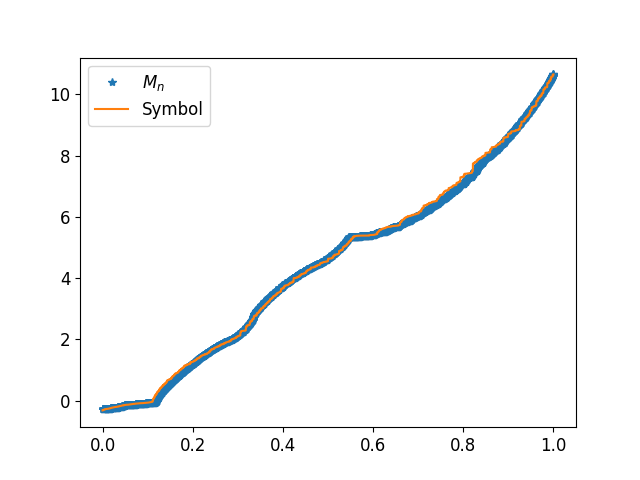}
	\includegraphics[width=0.45\linewidth]{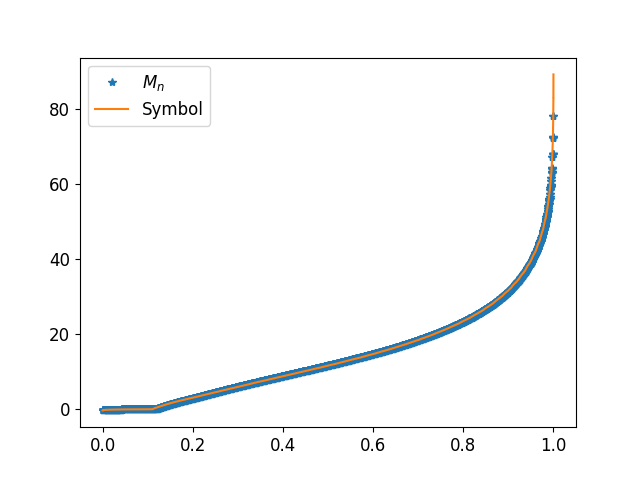}
	\hfill
	\includegraphics[width=0.32\linewidth]{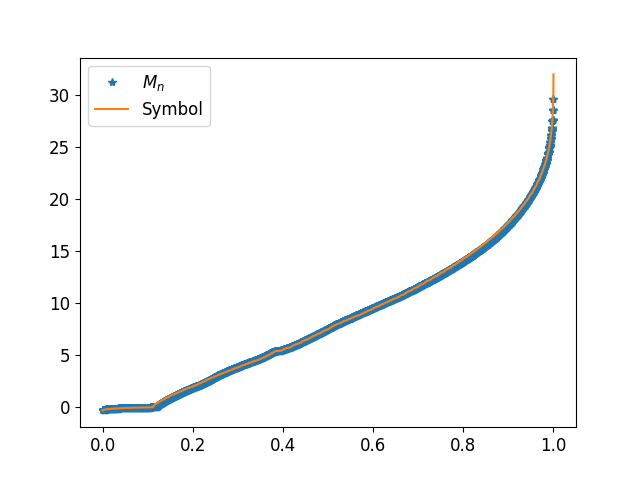}
	\includegraphics[width=0.32\linewidth]{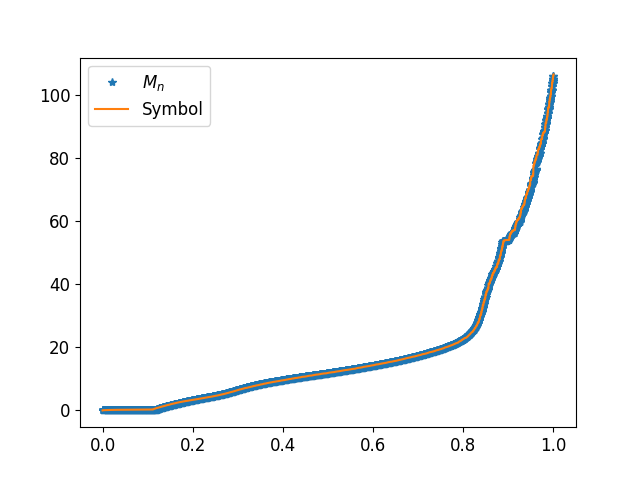}
	\includegraphics[width=0.32\linewidth]{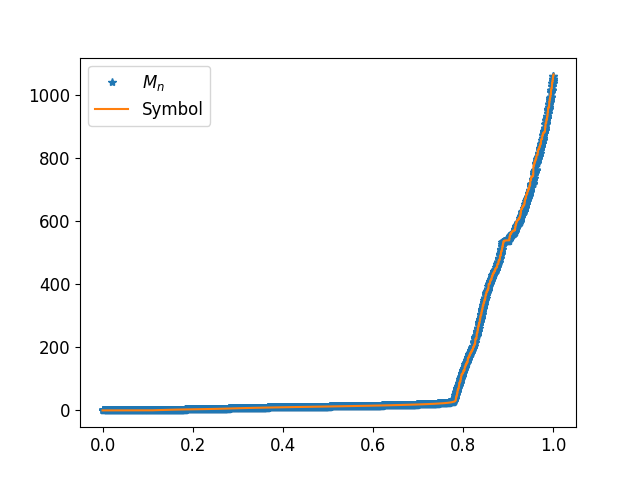}
	\caption{Adherence of the eigenvalues of $M_n$ to its symbol for $n=16$. The first line displays the plots corresponding to Group 1 (on the left) and Group 2 (on the right), while the second line shows Group 3; from left $\gamma = 1$, $\gamma = 10$, $\gamma = 100$.}
	\label{test_M123}
\end{figure}
\begin{figure}
	\centering
	\subfloat{\includegraphics[scale=.35]{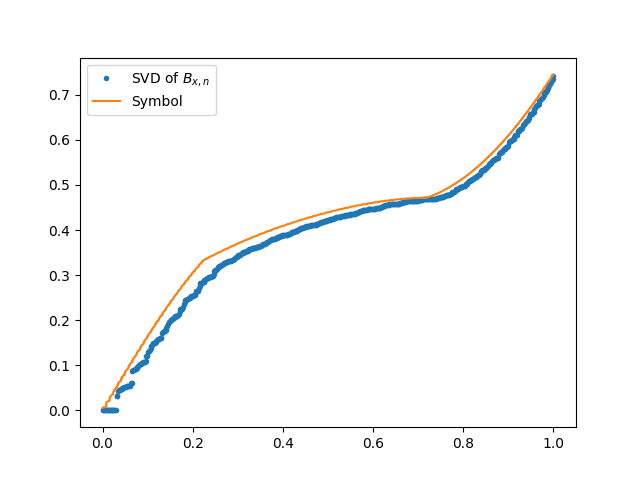}}
	\subfloat{\includegraphics[scale=.35]{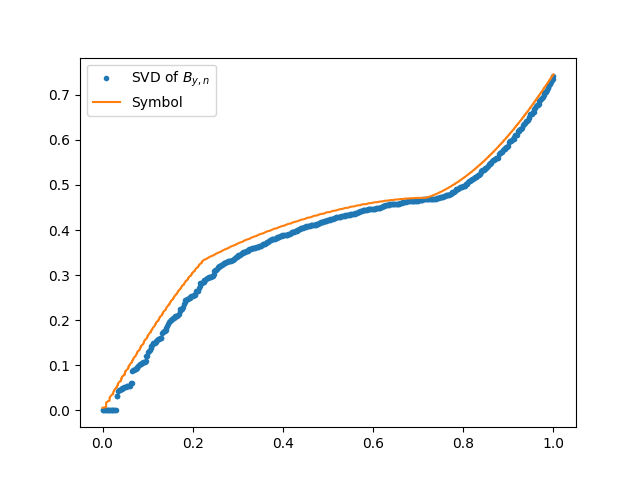}}
	\caption{Adherence of the singular values to its symbol for $n=16$. The plot on the left corresponds to the matrix sequence $\{ B_{x,n} \}_n$ and the one on the right to $\{B_{y,n}\}_n$.}
	\label{test_B}
\end{figure}
\begin{table}
	\centering
	\footnotesize
	\begin{tabular}{|c|c||c|c|c||c|c|c|}
		\toprule
		&  &  \multicolumn{3}{|c||}{Group 1}  & \multicolumn{3}{|c|}{Group 2} \\
		$n$ & $\dim(M_n)$ & Case a & Case b & Case c & Case a & Case b & Case c \\
		\midrule
		8 & 1107 & 57 & 98 & 88 & 59 & 107 & 97  \\
		16 & 4515 & 90 & 218 & 167 & 80 & 206 & 146  \\
		32 & 18243 & 154 & 625 & 444 & 118 & 554 & 407\\
		\bottomrule
	\end{tabular}
	\caption{Number of iterations required by PGMRES to achieve the tolerance of $10^{-5}$ across various values of $n$ for Group 1 and Group 2.}
	\label{gmresResult12}
\end{table}
\begin{table}
	\footnotesize
	\centering
	\begin{tabular}{|c||c|c|c||c|c|c||c|c|c|}
		\toprule
		&  \multicolumn{3}{|c||}{$\gamma=1$}  & \multicolumn{3}{|c||}{$\gamma=10$} & \multicolumn{3}{|c|}{$\gamma=100$} \\
		$n$ & Case a & Case b & Case c & Case a & Case b & Case c & Case a & Case b & Case c \\
		\midrule
		8 & 58 & 105 & 92 & 60 & 98 & 88 & 68 & 139 & 128\\
		16 & 85 & 218 & 147 & 90 & 227 & 158 & 92 & 314 & 253\\
		32 & 124 & 486 & 454 & 128 & 431 & 394 & 116 & 738 & 312\\
		\bottomrule
	\end{tabular}
	\caption{Number of iterations required by PGMRES to achieve the tolerance of $10^{-5}$ across various values of $n$ for Group 3.}
	\label{gmresResult3}
\end{table}

\paragraph{Eigenvalue distribution of \texorpdfstring{$A_n$}{A\_n}}
Figure \ref{test_An123} shows the adherence of the eigenvalue distribution with the sampling of the spectral symbol we found in (\ref{muConst_symbol}) and (\ref{muVar_symbol}).

\subparagraph{Absence of Outliers}
Throughout all tests, we do not observe eigenvalues of \(A_n\) outside the essential range of the symbol. The reason is both numerical and theoretical. Indeed, for fixed \(n\), {we} interpret $A_n= A_n(\mu)$
as a matrix-valued linear positive operator (LPO) in the sense of \cite{capizzano2000some}, due to its origin from a Galerkin approximation of a coercive operator \cite[Section 3]{capizzano2000some}. It then satisfies the monotonicity property
\begin{equation}
    A_n(\mu) \le A_n(\hat{\mu}),\quad\text{for every }\hat{\mu}\ge \mu,
\end{equation}
which, via the min-max theorem for Hermitian matrices, leads to
\begin{align*}
	\lambda_j(A_n({\rm essinf}\,\mu))
	&= [{\rm essinf}\,\mu] \,\lambda_j(A_n(1))
	\;\;\;\;\le \lambda_j\bigl(A_n(\mu)\bigr)\\
	&\le \lambda_j(A_n({\rm esssup}\,\mu))
	= [{\rm esssup}\,\mu] \,\lambda_j(A_n(1)).
\end{align*}
Since block multilevel Toeplitz matrices are also LPOs and their eigenvalues lie within the extrema of its generating function, say \((m,M)\), it follows that
\begin{equation}
    \lambda_j\bigl(A_{xx}(n)(\mu)\bigr)\in
\bigl(m\, [{\rm essinf}\,\mu],\; M\, [{\rm esssup}\,\mu]\bigr),
\end{equation}
hence no eigenvalue can emerge as an outlier beyond the essential range of the symbol.

\paragraph{Singular value distribution of \texorpdfstring{$B_x$}{B\_x} and \texorpdfstring{$B_y$}{B\_y}}
Figure~\ref{test_B} shows the adherence of the singular value distribution with the sampling of the spectral symbol we found in (\ref{Bxsymbol}) and (\ref{Bysymbol}).

\paragraph{Eigenvalues distribution of the full matrix $M_n$}
Figure \ref{test_M123} shows the adherence of the eigenvalue distribution to the sampled symbol evaluated by combining the symbols of $A_{x,n}$, $B_{x,n}$, and $B_{y,n}$ in the following way:
\begin{equation}
	\text{sym}(M_n) = \begin{bmatrix}
		\text{sym}(A_{x,n}) & 0 & \text{sym}(B_{x,n}) \\
		0 & \text{sym}(A_{x,n}) & \text{sym}(B_{y,n})\\
		\text{sym}(B_{x,n})^T & \text{sym}(B_{y,n})^T & 0
	\end{bmatrix}
\end{equation}

\begin{figure}
	\centering
	\includegraphics[width=0.45\linewidth]{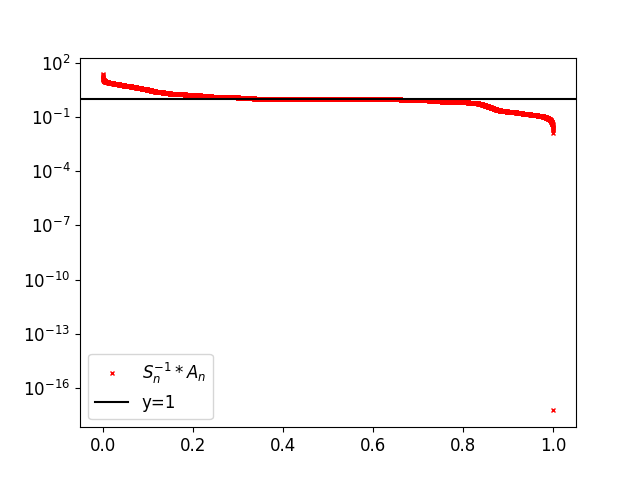}
	\includegraphics[width=0.45\linewidth]{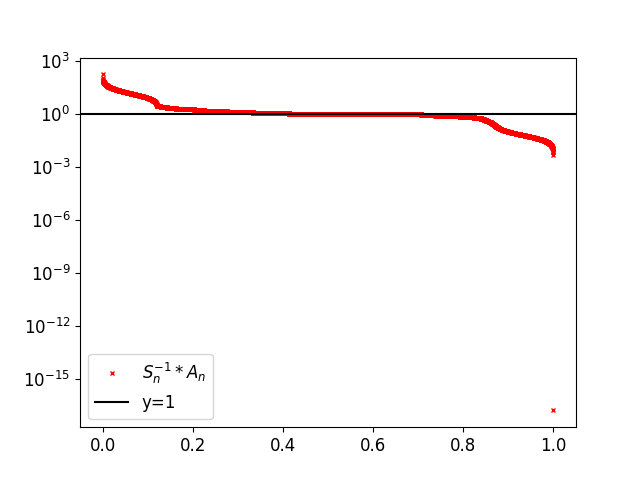}
	\hfill
	\includegraphics[width=0.32\linewidth]{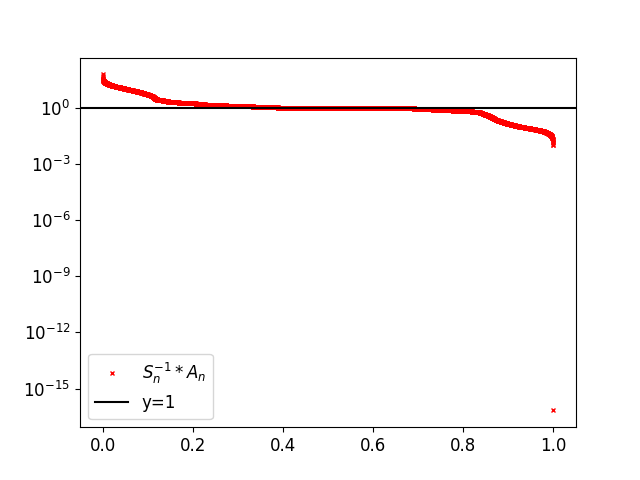}
	\includegraphics[width=0.32\linewidth]{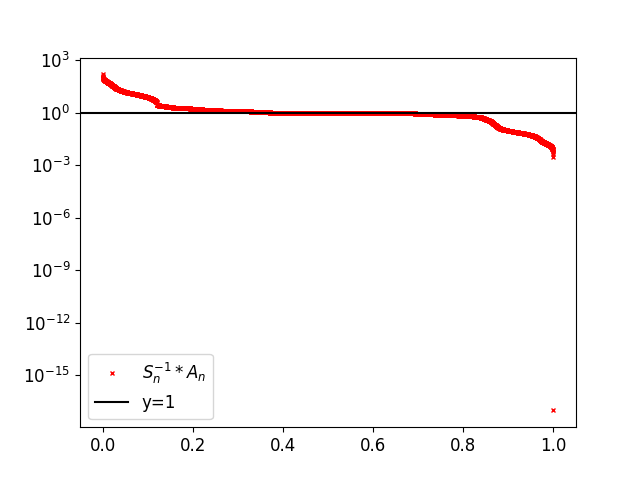}
	\includegraphics[width=0.32\linewidth]{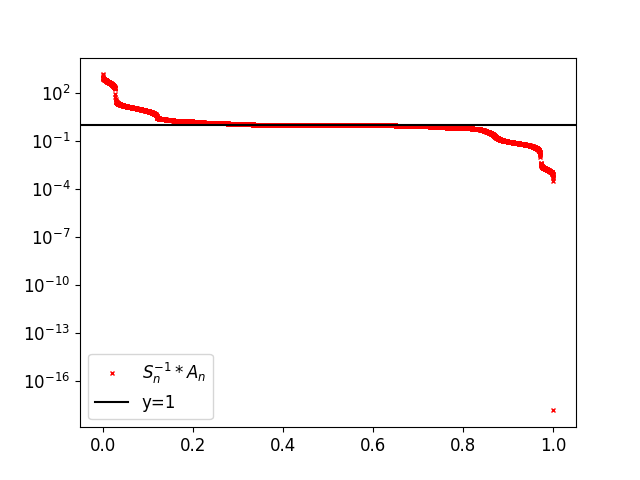}
	\caption{Cluster around 1 of the singular values of the preconditioned system. The first line displays the plots corresponding to Group 1 (on the left) and Group 2 (on the right), while the second line shows Group 3; from left $\gamma = 1$, $\gamma = 10$, $\gamma = 100$.}
	\label{fig:cluster12}
\end{figure}

\paragraph{Performance of the preliminary preconditioner}
The following experiments show how the approximation we proposed in Section \ref{sec:precond} reduces the number of iterations of GMRES, while GMRES without preconditioning does not converge within $1000$ iterations. Tables \ref{gmresResult12} and \ref{gmresResult3} show the number of iterations used by PGMRES to achieve the tolerance of $10^{-5}$, while figure \ref{fig:cluster12} displays the successful clustering of most of the singular values around 1, leaving only a small set of outliers. This clustering is the primary factor behind the rapid convergence of the preconditioned system.\\
Regarding the iteration count, we observe that the number of iterations grows linearly with $n$, i.e., it grows proportionally to the square root of the matrix size. Indeed this behaviour agrees with the theoretical results in \cite{MultiNo,noutsos2004matrix} for two dimensional problems.

\section{Conclusions}\label{sec:end}
In this work we have provided a spectral characterization of the blocks occurring in the matrix structure stemming from the approximation of a variable coefficient Stokes problem. The spectral analysis has been verified via proper numerical experiments and visualizations and it has been employed for proposing a basic GLT based preconditioning strategy. From the point of view of the practical application, these results should be regarded as very preliminary, {leading to a number of interesting open problems}, including those listed below.
\begin{itemize}
	\item Starting from application in Geophysics, such as the so-called Glacial Isostatic Adjustment (GIA) model, in \cite{Schur-CMAME} two techniques were considered on quadrilaterals: 1) the modified Taylor–Hood elements, where both displacements and pressure are discretized using bilinear basis functions but the pressure unknowns live on twice coarser mesh; 2) $\mathbb{Q}_2$-$\mathbb{Q}_1$ Taylor–Hood elements with biquadratic and bilinear basis functions for the displacements and for the pressure, correspondingly. We recall that the GIA model describes the response of the solid Earth to redistribution of mass due to alternating glaciation and deglaciation periods. While the tools employed in \cite{Schur-CMAME} are similar to those used here, the specific grids and the use of triangular elements has lead to a more complicated analysis. In fact, it is necessary to understand how to treat the two settings in a more uniform way. 
	\item The difficulties emphasized in the previous item and the fact that some parts of the analysis are based on plain use of the GLT axioms invite to put more efforts in designing procedures for the automatic computation of the GLT symbol; the idea belongs to Ratnani (see \cite[Chapter 11]{GS-I}) and it has been partially pursued in \cite{sarathkumar2024glt}, but still a lot of study has to be done.
	\item The indication for building proper preconditioners is interesting, but a lot of work remains to be performed with respect to the basic proposals given in the present work.
	\item We finally observe that the potential generalizations to the full elliptic case and to virtually any type of Galerkin approximation schemes are interesting goals for future research.
\end{itemize}

\printbibliography

\appendix
\section{Matrix entries of the stiffness matrix $A_n$}\label{app:A}
In this section we report the entries of the stiffness matrix $A_n$.

\paragraph{Zero column}
\begin{tikzpicture}[baseline=-0.1cm]
	\draw[step=0.5cm,gray,very thin] (-1,-1) grid (1,1);
	\draw[gray, very thin] (-1,0.5) -- (-0.5,1); 
	\draw[gray, very thin] (-1,0) -- (0,1); 
	\draw[gray, very thin] (-1,-0.5) -- (0.5,1);
	\draw[gray, very thin] (-1,-1) -- (1,1); 
	\draw[gray, very thin] (-0.5,-1) -- (1, 0.5);
	\draw[gray, very thin] (0,-1) -- (1, 0);
	\draw[gray, very thin] (0.5,-1) -- (1, -0.5);
	\draw[gray, very thin] (-1,-0.5) -- (-0.5, -1);
	\draw[gray, very thin] (-1,0) -- (0,-1);
	\draw[gray, very thin] (-1,0.5) -- (0.5, -1);
	\draw[gray, very thin] (-1,1) -- (1,-1);
	\draw[gray, very thin] (-0.5,1) -- (1, -0.5);
	\draw[gray, very thin] (0,1) -- (1,0);
	\draw[gray, very thin] (0.5,1) -- (1,0.5);
	\filldraw[red] (-0.875,-0.875) circle (1pt);
	\filldraw[blue] (-0.875,-0.65) circle (1pt);
	\filldraw[red] (-0.875,-0.375) circle (1pt);
	\filldraw[blue] (-0.875,-0.15) circle (1pt);
	\filldraw[red] (-0.875, 0.125) circle (1pt);
	\filldraw[blue] (-0.875, 0.35) circle (1pt);
	\filldraw[red] (-0.875, 0.625) circle (1pt);
	\filldraw[blue] (-0.875, 0.85) circle (1pt);
\end{tikzpicture}
\vspace{0.6cm}

\textcolor{red}{Even rows} $\xrightarrow{} \left| 
\begin{aligned}
	&i = 2k_1 \text{ and } i = (8n^2 -4n +1) + 2k_1 \\
	& \text{with } k_1 = 0,..,n-1   
\end{aligned} \right. $
\begin{equation*}
	A_{ij} = \left\{ \begin{aligned}
		&\frac{8}{3} (\mu_{4k_1} + \mu_{1 + 4k_1} )  \qquad &&\text{ if } j=i \\
		-&\frac{2}{3} (\mu_{4k_1} + \mu_{1 + 4k_1})\qquad &&\text{ if } j=i+2n \\
		-&\frac{4}{3} (\mu_{1 + 4k_1}) \qquad &&\text{ if } j = i+2n-1 \text{ and } \delta_0(k_1) = 0
	\end{aligned}
	\right. 
\end{equation*}

\textcolor{blue}{Odd rows} $\xrightarrow{} \left| 
\begin{aligned}
	&i = 1 + 2k_1 \text{ and } i = (8n^2 -4n +1) +1 +2k_1 \\
	& \text{with } k_1 = 0,..,n-1   
\end{aligned} \right. $
\begin{equation*}
	A_{ij} = \left\{ \begin{aligned}
		&\frac{8}{3} (\mu_{4k_1} + \mu_{2 + 4k_1})  \qquad &&\text{ if } j=i \\
		-&\frac{2}{3} (\mu_{4k_1} + \mu_{2 + 4k_1}) \qquad &&\text{ if } j=i+2n-1 \\
		-&\frac{4}{3} (\mu_{2 + 4k_1}) \qquad &&\text{ if } j = i+2n \text{ and } \delta_{n-1}(k_1) = 0
	\end{aligned}
	\right. 
\end{equation*}

\paragraph{First + n columns}
\begin{tikzpicture}[baseline=-0.1cm]
	\draw[step=0.5cm,gray,very thin] (-1,-1) grid (1,1);
	\draw[gray, very thin] (-1,0.5) -- (-0.5,1); 
	\draw[gray, very thin] (-1,0) -- (0,1); 
	\draw[gray, very thin] (-1,-0.5) -- (0.5,1);
	\draw[gray, very thin] (-1,-1) -- (1,1); 
	\draw[gray, very thin] (-0.5,-1) -- (1, 0.5);
	\draw[gray, very thin] (0,-1) -- (1, 0);
	\draw[gray, very thin] (0.5,-1) -- (1, -0.5);
	\draw[gray, very thin] (-1,-0.5) -- (-0.5, -1);
	\draw[gray, very thin] (-1,0) -- (0,-1);
	\draw[gray, very thin] (-1,0.5) -- (0.5, -1);
	\draw[gray, very thin] (-1,1) -- (1,-1);
	\draw[gray, very thin] (-0.5,1) -- (1, -0.5);
	\draw[gray, very thin] (0,1) -- (1,0);
	\draw[gray, very thin] (0.5,1) -- (1,0.5);
	\filldraw[red] (-0.75,-0.75) circle (1pt);
	\filldraw[red] (-0.25,-0.75) circle (1pt);
	\filldraw[red] (0.25,-0.75) circle (1pt);
	\filldraw[red] (0.75,-0.75) circle (1pt);
	\filldraw[cyan] (-0.75,-0.25) circle (1pt);
	\filldraw[cyan] (-0.25,-0.25) circle (1pt);
	\filldraw[cyan] (0.25,-0.25) circle (1pt);
	\filldraw[cyan] (0.75,-0.25) circle (1pt);
	\filldraw[green] (-0.75,0.25) circle (1pt);
	\filldraw[green] (-0.25,0.25) circle (1pt);
	\filldraw[green] (0.25,0.25) circle (1pt);
	\filldraw[green] (0.75,0.25) circle (1pt);
	\filldraw[blue] (-0.75,0.75) circle (1pt);
	\filldraw[blue] (-0.25,0.75) circle (1pt);
	\filldraw[blue] (0.25,0.75) circle (1pt);
	\filldraw[blue] (0.75,0.75) circle (1pt);
\end{tikzpicture}
\vspace{0.4cm}

\textcolor{red}{Bottom row} $\xrightarrow{} \left| 
\begin{aligned}
	&i = 2n+(8n -2)k_1 \text{ and }  \\
	&i = (8n^2 -4n +1) + 2n+(8n -2)k_1 \\
	& \text{with } k_1 = 0,..,n-1   
\end{aligned} \right. $\\
\begin{equation*}
	A_{ij} = \left\{ \begin{aligned}
		-&\frac{2}{3} (\mu_{4nk_1} + \mu_{1 + 4nk_1}) \qquad &&\text{ if } j=i-2n \\
		-&\frac{2}{3}  (\mu_{4nk_1} + \mu_{2 + 4nk_1}) \qquad &&\text{ if } j=i-(2n-1) \\
		&(\mu_{4nk_1} + \mu_{1 + 4nk_1} + \mu_{2 + 4nk_1} + \mu_{3 + 4nk_1}) \qquad &&\text{ if } j= i \\
		-&\frac{2}{3}  (\mu_{1 + 4nk_1} + \mu_{3 + 4nk_1}) \qquad &&\text{ if } j= i+(2n-1) \\
		-&\frac{2}{3}  (\mu_{2 + 4nk_1} + \mu_{3 + 4nk_1}) \qquad &&\text{ if } j= i+2n \\
		&\frac{1}{6} (\mu_{4nk_1} + \mu_{2 + 4nk_1}) \qquad &&\text{ if } j=i-2(2n-1) \text{ and } \delta_0(k_1) = 0 \\
		&\frac{1}{6} (\mu_{2 + 4nk_1} + \mu_{3 + 4nk_1})  \qquad &&\text{ if } j=i+4n \text{ and } \delta_{n-1}(k_1) = 0
	\end{aligned}
	\right. 
\end{equation*}

\textcolor{cyan}{Middle (odd) rows} $\xrightarrow{}
\left| 
\begin{aligned}
	&i = (2n+1)+(8n-2)k_1 + 2k_2 \text{ and }  \\
	&i = (8n^2 -4n +1) + (2n+1)+(8n-2)k_1 + 2k_2 \\
	& \text{with } k_1 = 0,..,n-1  \text{ and }  k_2=0,..,n-2
\end{aligned} \right. $
\begin{equation*}
	A_{ij} = \left\{ \begin{aligned}
		-&\frac{4}{3} (\mu_{2 + 4nk_1 + 4k_2})  \qquad &&\text{ if } j= i-2n \\
		-&\frac{4}{3} (\mu_{5 + 4nk_1 + 4k_2})  \qquad &&\text{ if } j=i-(2n-1) \\
		&\frac{8}{3} (\mu_{2 + 4nk_1 + 4k_2} + \mu_{5 + 4nk_1 + 4k_2})  \qquad &&\text{ if } j=i \\
		-&\frac{4}{3}(\mu_{2 + 4nk_1+4k_2})  \qquad &&\text{ if } j=i+(2n-1) \\
		-&\frac{4}{3} (\mu_{5 + 4nk_1+4k_2})  \qquad &&\text{ if } j=i+2n 
	\end{aligned}
	\right. 
\end{equation*}

\textcolor{green}{Middle (even) rows} $\xrightarrow{} \left| 
\begin{aligned}
	&i = 2n + (8n-2)k_1 +2k_2 \text{ and }  \\
	&i = (8n^2 -4n +1) + 2n+(8n-2)k_1+2k_2 \\
	& \text{with } k_1 = 0,..,n-1  \text{ and }  k_2=1,..,n-2
\end{aligned} \right. $
\begin{equation*}
	\resizebox{\linewidth}{!}{%
	$A_{ij} = \left\{ \begin{aligned}
		-&\frac{2}{3} (\mu_{4nk_1 + 4k_2} + \mu_{1 + 4nk_1 + 4k_2}) \qquad &&\text{ if } j=i-2n \\
		-&\frac{2}{3}  (\mu_{4nk_1 + 4k_2} + \mu_{2 + 4nk_1 + 4k_2}) \qquad &&\text{ if } j=i-(2n-1) \\
		&(\mu_{4nk_1 + 4k_2} + \mu_{1 + 4nk_1 + 4k_2} + \mu_{2 + 4nk_1 + 4k_2} + \mu_{3 + 4nk_1 + 4k_2})  \qquad &&\text{ if } j= i \\
		-&\frac{2}{3} (\mu_{1 + 4nk_1 + 4k_2} + \mu_{3 + 4nk_1 + 4k_2})  \qquad &&\text{ if } j= i+(2n-1) \\
		-&\frac{2}{3} (\mu_{2 + 4nk_1 + 4k_2} + \mu_{3 + 4nk_1 + 4k_2})  \qquad &&\text{ if } j= i+2n \\
		&\frac{1}{6}  (\mu_{4nk_1 + 4k_2} + \mu_{1 + 4nk_1 + 4k_2}) \qquad &&\text{ if } j= i-4n \text{ and } \delta_0(k_1) = 0 \\
		&\frac{1}{6}   (\mu_{4nk_1 + 4k_2} + \mu_{2 + 4nk_1 + 4k_2}) \qquad &&\text{ if } j= i-2(2n-1)\text{ and } \delta_0(k_1) = 0 \\
		&\frac{1}{6}  (\mu_{1 + 4nk_1 + 4k_2} + \mu_{3 + 4nk_1 + 4k_2}) \qquad &&\text{ if } j= i+2(2n-1) \text{ and } \delta_{n-1}(k_1) = 0 \\
		&\frac{1}{6}   (\mu_{2 + 4nk_1 + 4k_2} + \mu_{3 + 4nk_1 + 4k_2})  \qquad &&\text{ if } j=i+4n \text{ and } \delta_{n-1}(k_1) = 0                                
	\end{aligned}
	\right. $}
\end{equation*}

\textcolor{blue}{Top row} $\xrightarrow{} \left| 
\begin{aligned}
	&i = (4n-2) + (8n-2)k_1 \text{ and }  \\
	&i = (8n^2 -4n +1) + (4n-2)+(8n-2)k_1 \\
	& \text{with } k_1 = 0,..,n-1  
\end{aligned} \right. $
\begin{equation*}
	A_{ij} = \left\{ \begin{aligned}
		-&\frac{2}{3}  (\mu_{4(n-1) + 4nk_1} + \mu_{4(n-1) + 1 + 4nk_1}) \qquad&&\text{ if } j= i-2n \\
		-&\frac{2}{3} (\mu_{4(n-1) + 4nk_1} + \mu_{4(n-1) + 2 + 4nk_1}) \qquad&&\text{ if } j= i-(2n-1)\\
		&(\mu_{4(n-1) + 4nk_1} + \mu_{4(n-1) + 1 + 4nk_1} +  \\
		&\mu_{4(n-1) + 2 + 4nk_1} + \mu_{4(n-1) + 3 + 4nk_1}) \qquad&&\text{ if } j= i\\
		-&\frac{2}{3}  (\mu_{4(n-1) + 1 + 4nk_1} + \mu_{4(n-1) + 3 + 4nk_1}) \qquad&&\text{ if } j= i+(2n-1)\\
		-&\frac{2}{3}  (\mu_{4(n-1) + 2 + 4nk_1} + \mu_{4(n-1) + 3 + 4nk_1}) \qquad&&\text{ if } j= i+2n\\
		&\frac{1}{6}  (\mu_{4(n-1) + 4nk_1} + \mu_{4(n-1) + 1 + 4nk_1})  \qquad&&\text{ if } j= i-4n \text{ and } \delta_{0}(k_1) =0\\
		&\frac{1}{6}  (\mu_{4(n-1) + 1 + 4nk_1} + \mu_{4(n-1) + 3 + 4nk_1}) \qquad&&\text{ if } j=i+2(2n-1) \text{ and } \delta_{n-1}(k_1)=0              
	\end{aligned}
	\right. 
\end{equation*}

\paragraph{Second + (n-1) columns}
\begin{tikzpicture}[baseline=-0.1cm]
	\draw[step=0.5cm,gray,very thin] (-1,-1) grid (1,1);
	\draw[gray, very thin] (-1,0.5) -- (-0.5,1); 
	\draw[gray, very thin] (-1,0) -- (0,1); 
	\draw[gray, very thin] (-1,-0.5) -- (0.5,1);
	\draw[gray, very thin] (-1,-1) -- (1,1); 
	\draw[gray, very thin] (-0.5,-1) -- (1, 0.5);
	\draw[gray, very thin] (0,-1) -- (1, 0);
	\draw[gray, very thin] (0.5,-1) -- (1, -0.5);
	\draw[gray, very thin] (-1,-0.5) -- (-0.5, -1);
	\draw[gray, very thin] (-1,0) -- (0,-1);
	\draw[gray, very thin] (-1,0.5) -- (0.5, -1);
	\draw[gray, very thin] (-1,1) -- (1,-1);
	\draw[gray, very thin] (-0.5,1) -- (1, -0.5);
	\draw[gray, very thin] (0,1) -- (1,0);
	\draw[gray, very thin] (0.5,1) -- (1,0.5);
	\filldraw[red] (-0.625,-0.875) circle (1pt);
	\filldraw[red] (-0.625,-0.375) circle (1pt);
	\filldraw[red] (-0.625,0.125) circle (1pt);
	\filldraw[red] (-0.625,0.625) circle (1pt);
	\filldraw[red] (-0.125,-0.875) circle (1pt);
	\filldraw[red] (-0.125,-0.375) circle (1pt);
	\filldraw[red] (-0.125,0.125) circle (1pt);
	\filldraw[red] (-0.125,0.625) circle (1pt);
	\filldraw[red] (0.375,-0.875) circle (1pt);
	\filldraw[red] (0.375,-0.375) circle (1pt);
	\filldraw[red] (0.375,0.125) circle (1pt);
	\filldraw[red] (0.375,0.625) circle (1pt);
	\filldraw[blue] (-0.625,-0.625) circle (1pt);
	\filldraw[blue] (-0.625,-0.125) circle (1pt);
	\filldraw[blue] (-0.625,0.375) circle (1pt);
	\filldraw[blue] (-0.625,0.875) circle (1pt);
	\filldraw[blue] (-0.125,-0.625) circle (1pt);
	\filldraw[blue] (-0.125,-0.125) circle (1pt);
	\filldraw[blue] (-0.125,0.375) circle (1pt);
	\filldraw[blue] (-0.125,0.875) circle (1pt);
	\filldraw[blue] (0.375,-0.625) circle (1pt);
	\filldraw[blue] (0.375,-0.125) circle (1pt);
	\filldraw[blue] (0.375,0.375) circle (1pt);
	\filldraw[blue] (0.375,0.875) circle (1pt);
\end{tikzpicture}
\vspace{0.6cm}

\textcolor{red}{Even rows} $\xrightarrow{}
\left| 
\begin{aligned}
	&i = 4n -1 + (8n -2)k_1 + 2k_2 \text{ and }  \\
	&i = (8n^2 -4n +1) + 4n -1+(8n -2)k_1+2k_2 \\
	& \text{with } k_1 = 0,..,n-2 \text{ and } k_2 = 0,..,n-1
\end{aligned} \right. $
\begin{equation*}
	A_{ij} = \left\{ \begin{aligned}
		-&\frac{2}{3} (\mu_{1 + 4nk_1 + 4k_2} + \mu_{3 + 4nk_1 + 4k_2}) \qquad&&\text{ if } j= i-(2n-1)\\
		&\frac{8}{3} (\mu_{1 + 4nk_1 + 4k_2} + \mu_{3 + 4nk_1 + 4k_2}) \qquad&&\text{ if } j= i\\
		-&\frac{4}{3} \, \mu_{3 + 4nk_1 + 4k_2} \qquad&&\text{ if } j= i+2n\\
		-&\frac{4}{3} (\mu_{1 + 4nk_1 + 4k_2}) \qquad&&\text{ if } j= i-2n \text{ and } \delta_{0}(k_2) = 0 \\
		-&\frac{2}{3}  (\mu_{1 + 4nk_1 + 4k_2} + \mu_{3 + 4nk_1 + 4k_2}) \qquad&&\text{ if } j= i+(2n-1) \text{ and } \delta_{0}(k_2) = 0
	\end{aligned}
	\right. 
\end{equation*}

\textcolor{blue}{Odd rows} $\xrightarrow{}
\left| 
\begin{aligned}
	&i = 4n + (8n -2)k_1 + 2k_2 \text{ and }  \\
	&i = (8n^2 -4n +1) + 4n+(8n -2)k_1+2k_2 \\
	& \text{with } k_1 = 0,..,n-2 \text{ and } k_2 = 0,..,n-1
\end{aligned} \right. $
\begin{equation*}
	A_{ij} = \left\{ \begin{aligned}
		-&\frac{2}{3}  (\mu_{2 + 4nk_1 + 4k_2} + \mu_{3 + 4nk_1 + 4k_2}) \qquad &&\text{ if } j=i-2n \\
		&\frac{8}{3}  (\mu_{2 + 4nk_1 + 4k_2} + \mu_{3 + 4nk_1 + 4k_2}) \qquad &&\text{ if } j= i\\
		-&\frac{4}{3} \, \mu_{3 + 4nk_1 + 4k_2} \qquad &&\text{ if } j=i+(2n-1) \\
		-&\frac{4}{3}  (\mu_{2 + 4nk_1 + 4k_2}) \qquad &&\text{ if } j= i-(2n-1) \text{ and } \delta_{n-1}(k_2) = 0 \\
		-&\frac{2}{3}  (\mu_{2 + 4nk_1 + 4k_2} + \mu_{3 + 4nk_1 + 4k_2}) \qquad &&\text{ if } j= i+2n \text{ and } \delta_{n-1}(k_2) = 0
	\end{aligned}
	\right. 
\end{equation*}

\paragraph{Third + (n-1) columns}
\begin{tikzpicture}[baseline=-0.1cm]
	\draw[step=0.5cm,gray,very thin] (-1,-1) grid (1,1);
	\draw[gray, very thin] (-1,0.5) -- (-0.5,1); 
	\draw[gray, very thin] (-1,0) -- (0,1); 
	\draw[gray, very thin] (-1,-0.5) -- (0.5,1);
	\draw[gray, very thin] (-1,-1) -- (1,1); 
	\draw[gray, very thin] (-0.5,-1) -- (1, 0.5);
	\draw[gray, very thin] (0,-1) -- (1, 0);
	\draw[gray, very thin] (0.5,-1) -- (1, -0.5);
	\draw[gray, very thin] (-1,-0.5) -- (-0.5, -1);
	\draw[gray, very thin] (-1,0) -- (0,-1);
	\draw[gray, very thin] (-1,0.5) -- (0.5, -1);
	\draw[gray, very thin] (-1,1) -- (1,-1);
	\draw[gray, very thin] (-0.5,1) -- (1, -0.5);
	\draw[gray, very thin] (0,1) -- (1,0);
	\draw[gray, very thin] (0.5,1) -- (1,0.5);
	\filldraw[red] (-0.5,-0.75) circle (1pt);
	\filldraw[red] (-0.5,-0.25) circle (1pt);
	\filldraw[red] (-0.5,0.25) circle (1pt);
	\filldraw[red] (-0.5,0.75) circle (1pt);
	\filldraw[red] (0,-0.75) circle (1pt);
	\filldraw[red] (0,-0.25) circle (1pt);
	\filldraw[red] (0,0.25) circle (1pt);
	\filldraw[red] (0,0.75) circle (1pt);
	\filldraw[red] (0.5,-0.75) circle (1pt);
	\filldraw[red] (0.5,-0.25) circle (1pt);
	\filldraw[red] (0.5,0.25) circle (1pt);
	\filldraw[red] (0.5,0.75) circle (1pt);
	\filldraw[blue] (-0.5,-0.5) circle (1pt);
	\filldraw[blue] (-0.5,0) circle (1pt);
	\filldraw[blue] (-0.5,0.5) circle (1pt);
	\filldraw[blue] (0,-0.5) circle (1pt);
	\filldraw[blue] (0,0) circle (1pt);
	\filldraw[blue] (0,0.5) circle (1pt);
	\filldraw[blue] (0.5,-0.5) circle (1pt);
	\filldraw[blue] (0.5,0) circle (1pt);
	\filldraw[blue] (0.5,0.5) circle (1pt);
\end{tikzpicture}
\vspace{0.6cm}

\textcolor{red}{Even rows} $\xrightarrow{} \left| 
\begin{aligned}
	&i = (6n-1) + (8n-2)k_1 + 2k_2 \text{ and }  \\
	&i = (8n^2 -4n +1) + (6n-1)+(8n-2)k_1+2k_2 \\
	& \text{with } k_1 = 0,..,n-2 \text{ and } k_2 = 0,..,n-1
\end{aligned} \right. $
\begin{equation*}
	A_{ij} = \left\{ \begin{aligned}
		-&\frac{4}{3} \, \mu_{3 + 4nk_1 + 4k_2} \qquad &&\text{ if } j= i-2n \\
		-&\frac{4}{3} \, \mu_{3 + 4nk_1 + 4k_2} \qquad &&\text{ if } j= i-(2n-1) \\
		&\frac{8}{3} (\mu_{3 + 4nk_1 + 4k_2} + \mu_{4n + 4nk_1 + 4k_2}) \qquad &&\text{ if } j= i \\
		-&\frac{4}{3} \, \mu_{4n + 4nk_1 + 4k_2} \qquad &&\text{ if } j= i+(2n-1) \\
		-&\frac{4}{3} \, \mu_{4n + 4nk_1 + 4k_2} \qquad &&\text{ if } j= i+2n
	\end{aligned}
	\right. 
\end{equation*}

\textcolor{blue}{Odd rows} $\xrightarrow{} \left| 
\begin{aligned}
	&i = 6n + (8n-2)k_1 + 2k_2  \text{ and }  \\
	&i = (8n^2 -4n +1) + 6n+(8n-2)k_1+2k_2  \\
	& \text{with } k_1 = 0,..,n-2 \text{ and } k_2 = 0,..,n-2
\end{aligned} \right. $
\begin{equation*}
	\resizebox{\linewidth}{!}{%
	$A_{ij} = \left\{ \begin{aligned}
		&\frac{1}{6} (\mu_{2 + 4nk_1 + 4k_2} + \mu_{3 + 4nk_1 + 4k_2}) \qquad &&\text{ if } j= i-4n \\
		&\frac{1}{6} (\mu_{5 + 4nk_1 + 4k_2} + \mu_{7 + 4nk_1 + 4k_2}) \qquad &&\text{ if } j= i-2(2n-1) \\
		-&\frac{2}{3}  (\mu_{2 + 4nk_1 + 4k_2} + \mu_{3 + 4nk_1 + 4k_2}) \qquad &&\text{ if } j= i-2n \\
		-&\frac{2}{3}  (\mu_{5 + 4nk_1 + 4k_2} + \mu_{7 + 4nk_1 + 4k_2}) \qquad &&\text{ if } j= i-(2n-1) \\
		&\frac{1}{2} ( \mu_{2 + 4nk_1 + 4k_2} + \mu_{3 + 4nk_1 + 4k_2} + \mu_{5 + 4nk_1 + 4k_2} + \mu_{7 + 4nk_1 + 4k_2} + \mu_{4n + 4nk_1 + 4k_2} \\ &+ \mu_{4n+2 + 4nk_1 + 4k_2} + \mu_{4n+4 + 4nk_1 + 4k_2} + \mu_{4n+5 + 4nk_1 + 4k_2}) \qquad &&\text{ if } j= i \\
		-&\frac{2}{3}  (\mu_{4n + 4nk_1 + 4k_2} + \mu_{4n+2 + 4nk_1 + 4k_2}) \qquad &&\text{ if } j= i+(2n-1) \\
		-&\frac{2}{3}  (\mu_{4n+4 + 4nk_1 + 4k_2} + \mu_{4n+5 + 4nk_1 + 4k_2})  \qquad &&\text{ if } j= i+2n \\
		&\frac{1}{6} (\mu_{4n + 4nk_1 + 4k_2} + \mu_{4n+2 + 4nk_1 + 4k_2})  \qquad &&\text{ if } j= i+2(2n-1) \\
		&\frac{1}{6} (\mu_{4n+4 + 4nk_1 + 4k_2} + \mu_{4n+5 + 4nk_1 + 4k_2}) \qquad &&\text{ if } j= i+4n 
	\end{aligned}
	\right. $}
\end{equation*}

\paragraph{Fourth + (n-1) columns}
\begin{tikzpicture}[baseline=-0.1cm]
	\draw[step=0.5cm,gray,very thin] (-1,-1) grid (1,1);
	\draw[gray, very thin] (-1,0.5) -- (-0.5,1); 
	\draw[gray, very thin] (-1,0) -- (0,1); 
	\draw[gray, very thin] (-1,-0.5) -- (0.5,1);
	\draw[gray, very thin] (-1,-1) -- (1,1); 
	\draw[gray, very thin] (-0.5,-1) -- (1, 0.5);
	\draw[gray, very thin] (0,-1) -- (1, 0);
	\draw[gray, very thin] (0.5,-1) -- (1, -0.5);
	\draw[gray, very thin] (-1,-0.5) -- (-0.5, -1);
	\draw[gray, very thin] (-1,0) -- (0,-1);
	\draw[gray, very thin] (-1,0.5) -- (0.5, -1);
	\draw[gray, very thin] (-1,1) -- (1,-1);
	\draw[gray, very thin] (-0.5,1) -- (1, -0.5);
	\draw[gray, very thin] (0,1) -- (1,0);
	\draw[gray, very thin] (0.5,1) -- (1,0.5);
	\filldraw[red] (-0.375,-0.875) circle (1pt);
	\filldraw[red] (-0.375,-0.375) circle (1pt);
	\filldraw[red] (-0.375, 0.125) circle (1pt);
	\filldraw[red] (-0.375, 0.625) circle (1pt);
	\filldraw[blue] (-0.375,-0.65) circle (1pt);
	\filldraw[blue] (-0.375,-0.15) circle (1pt);
	\filldraw[blue] (-0.375, 0.35) circle (1pt);
	\filldraw[blue] (-0.375, 0.85) circle (1pt);
	\filldraw[red] (0.125,-0.875) circle (1pt);
	\filldraw[red] (0.125,-0.375) circle (1pt);
	\filldraw[red] (0.125, 0.125) circle (1pt);
	\filldraw[red] (0.125, 0.625) circle (1pt);
	\filldraw[blue] (0.125,-0.65) circle (1pt);
	\filldraw[blue] (0.125,-0.15) circle (1pt);
	\filldraw[blue] (0.125, 0.35) circle (1pt);
	\filldraw[blue] (0.125, 0.85) circle (1pt);
	\filldraw[red] (0.625,-0.875) circle (1pt);
	\filldraw[red] (0.625,-0.375) circle (1pt);
	\filldraw[red] (0.625, 0.125) circle (1pt);
	\filldraw[red] (0.625, 0.625) circle (1pt);
	\filldraw[blue] (0.625,-0.65) circle (1pt);
	\filldraw[blue] (0.625,-0.15) circle (1pt);
	\filldraw[blue] (0.625, 0.35) circle (1pt);
	\filldraw[blue] (0.625, 0.85) circle (1pt);
\end{tikzpicture}
\vspace{0.6cm}

\textcolor{red}{Even rows} $\xrightarrow{} \left| 
\begin{aligned}
	&i = (8n-2)k_1 + 2k_2 \text{ and }  \\
	&i = (8n^2 -4n +1) + (8n-2)k_1+2k_2 \\
	& \text{with } k_1 = 1,..,n-1 \text{ and } k_2 = 0,..,n-1
\end{aligned} \right. $

\begin{equation*}
	A_{ij} = \left\{ \begin{aligned}
		-&\frac{4}{3}  \mu_{4nk_1 + 4k_2} \qquad &&\text{ if } j=i-(2n-1) \\
		&\frac{8}{3}  (\mu_{4nk_1 + 4k_2} + \mu_{1 + 4nk_1 + 4k_2}) \qquad &&\text{ if } j=i \\
		-&\frac{2}{3}  (\mu_{4nk_1 + 4k_2} + \mu_{1 + 4nk_1 + 4k_2}) \qquad &&\text{ if } j=i+2n \\
		-&\frac{2}{3}   (\mu_{4nk_1 + 4k_2} + \mu_{1 + 4nk_1 + 4k_2})  \qquad &&\text{ if } j= i-2n \text{ and } \delta_0(k_2)=0 \\ 
		-&\frac{4}{3}  (\mu_{1 + 4nk_1 + 4k_2}) \qquad &&\text{ if } j= i+(2n-1) \text{ and } \delta_0(k_2)=0
	\end{aligned}
	\right. 
\end{equation*}

\textcolor{blue}{Odd rows} $\xrightarrow{} \left| 
\begin{aligned}
	&i = 1 + (8n-2)k_1 + 2k_2 \text{ and }  \\
	&i = (8n^2 -4n +1) + 1+(8n-2)k_1+2k_2 \\
	& \text{with } k_1 = 1,..,n-1 \text{ and } k_2 = 0,..,n-1
\end{aligned} \right. $
\begin{equation*}
	A_{ij} = \left\{ \begin{aligned}
		-&\frac{4}{3} \mu_{4nk_1 + 4k_2} \qquad &&\text{ if } j= i-2n \\
		&\frac{8}{3} (\mu_{4nk_1 + 4k_2} + \mu_{2 + 4nk_1 + 4k_2}) \qquad &&\text{ if } j= i \\
		-&\frac{2}{3} (\mu_{4nk_1 + 4k_2} + \mu_{2 + 4nk_1 + 4k_2}) \qquad &&\text{ if } j= i+(2n-1) \\
		-&\frac{2}{3} (\mu_{4nk_1 + 4k_2} + \mu_{2 + 4nk_1 + 4k_2}) \qquad &&\text{ if } j= i-(2n-1) \text{ and } \delta_{n-1}(k_2)=0 \\
		-&\frac{4}{3} (\mu_{2 + 4nk_1 + 4k_2}) \qquad &&\text{ if } j= i+2n \text{ and } \delta_{n-1}(k_2)=0        
	\end{aligned}
	\right. 
\end{equation*}

\paragraph{Last column}
\begin{tikzpicture}[baseline=-0.1cm]
	\draw[step=0.5cm,gray,very thin] (-1,-1) grid (1,1);
	\draw[gray, very thin] (-1,0.5) -- (-0.5,1); 
	\draw[gray, very thin] (-1,0) -- (0,1); 
	\draw[gray, very thin] (-1,-0.5) -- (0.5,1);
	\draw[gray, very thin] (-1,-1) -- (1,1); 
	\draw[gray, very thin] (-0.5,-1) -- (1, 0.5);
	\draw[gray, very thin] (0,-1) -- (1, 0);
	\draw[gray, very thin] (0.5,-1) -- (1, -0.5);
	\draw[gray, very thin] (-1,-0.5) -- (-0.5, -1);
	\draw[gray, very thin] (-1,0) -- (0,-1);
	\draw[gray, very thin] (-1,0.5) -- (0.5, -1);
	\draw[gray, very thin] (-1,1) -- (1,-1);
	\draw[gray, very thin] (-0.5,1) -- (1, -0.5);
	\draw[gray, very thin] (0,1) -- (1,0);
	\draw[gray, very thin] (0.5,1) -- (1,0.5);
	\filldraw[red] (0.875,-0.875) circle (1pt);
	\filldraw[red] (0.875,-0.375) circle (1pt);
	\filldraw[red] (0.875,0.125) circle (1pt);
	\filldraw[red] (0.875,0.625) circle (1pt);
	\filldraw[blue] (0.875,-0.625) circle (1pt);
	\filldraw[blue] (0.875,-0.125) circle (1pt);
	\filldraw[blue] (0.875,0.375) circle (1pt);
	\filldraw[blue] (0.875,0.875) circle (1pt);
\end{tikzpicture}
\vspace{0.6cm}

\textcolor{red}{Even rows} $\xrightarrow{} \left| 
\begin{aligned}
	&i = 8n^2-6n+1 + 2k_1\text{ and }  \\
	&i = (8n^2 -4n +1) + 8n^2-6n+1+2k_1 \\
	& \text{with } k_1 = 0,..,n-1 
\end{aligned} \right. $
\begin{equation*}
	A_{ij} = \left\{ \begin{aligned}
		-&\frac{2}{3} (\mu_{4n(n-1) + 1 + 4k_1} + \mu_{4n(n-1) + 3 + 4k_1}) \qquad &&\text{ if } j=i-(2n-1) \\
		&\frac{8}{3} (\mu_{4n(n-1) + 1 + 4k_1} + \mu_{4n(n-1) + 3 + 4k_1}) \qquad &&\text{ if } j= i \\
		-&\frac{4}{3} \mu_{4n(n-1) + 1 + 4k_1} \qquad &&\text{ if } j= i-2n \text{ and } \delta_{0}(k_1)=0
	\end{aligned}
	\right. 
\end{equation*}

\textcolor{blue}{Odd rows} $\xrightarrow{} \left| 
\begin{aligned}
	&i = 8n^2-6n+2 + 2k_1 \text{ and }  \\
	&i = (8n^2 -4n +1) + 8n^2-6n+2+2k_1 \\
	& \text{with } k_1 = 0,..,n-1 
\end{aligned} \right. $
\begin{equation*}
	A_{ij} = \left\{ \begin{aligned}
		-&\frac{2}{3} (\mu_{4n(n-1) + 2 + 4k_1} + \mu_{4n(n-1) + 3 + 4k_1}) \qquad &&\text{ if } j= i-2n \\
		&\frac{8}{3} (\mu_{4n(n-1) + 2 + 4k_1} + \mu_{4n(n-1) + 3 + 4k_1}) \qquad &&\text{ if } j= i \\
		-&\frac{4}{3}  \mu_{4n(n-1) + 2 + 4k_1} \qquad &&\text{ if } j= i-(2n-1) \text{ and } \delta_{n-1}(k_1)=0
	\end{aligned}
	\right. 
\end{equation*}

\newpage
\section{Matrix entries of $B_n^T$}\label{app:B}
In this section we report the matrix entries of $B_{x,n}$ and $B_{y,n}$.

\subsection{Matrix $B_{x,n}$}

\paragraph{Zero column + n columns}
\begin{tikzpicture}[baseline=-0.1cm]
	\draw[step=0.5cm,gray,very thin] (-1,-1) grid (1,1);
	\draw[gray, very thin] (-1,0.5) -- (-0.5,1); 
	\draw[gray, very thin] (-1,0) -- (0,1); 
	\draw[gray, very thin] (-1,-0.5) -- (0.5,1);
	\draw[gray, very thin] (-1,-1) -- (1,1); 
	\draw[gray, very thin] (-0.5,-1) -- (1, 0.5);
	\draw[gray, very thin] (0,-1) -- (1, 0);
	\draw[gray, very thin] (0.5,-1) -- (1, -0.5);
	\draw[gray, very thin] (-1,-0.5) -- (-0.5, -1);
	\draw[gray, very thin] (-1,0) -- (0,-1);
	\draw[gray, very thin] (-1,0.5) -- (0.5, -1);
	\draw[gray, very thin] (-1,1) -- (1,-1);
	\draw[gray, very thin] (-0.5,1) -- (1, -0.5);
	\draw[gray, very thin] (0,1) -- (1,0);
	\draw[gray, very thin] (0.5,1) -- (1,0.5);
	\filldraw[red] (-0.875,-0.875) circle (1pt);
	\filldraw[red] (-0.875,-0.375) circle (1pt);
	\filldraw[red] (-0.875, 0.125) circle (1pt);
	\filldraw[red] (-0.875, 0.625) circle (1pt);
	\filldraw[red] (-0.375,-0.875) circle (1pt);
	\filldraw[red] (-0.375,-0.375) circle (1pt);
	\filldraw[red] (-0.375, 0.125) circle (1pt);
	\filldraw[red] (-0.375, 0.625) circle (1pt);
	\filldraw[red] (0.125,-0.875) circle (1pt);
	\filldraw[red] (0.125,-0.375) circle (1pt);
	\filldraw[red] (0.125, 0.125) circle (1pt);
	\filldraw[red] (0.125, 0.625) circle (1pt);
	\filldraw[red] (0.625,-0.875) circle (1pt);
	\filldraw[red] (0.625,-0.375) circle (1pt);
	\filldraw[red] (0.625, 0.125) circle (1pt);
	\filldraw[red] (0.625, 0.625) circle (1pt);
	\filldraw[blue] (-0.875,-0.65) circle (1pt);
	\filldraw[blue] (-0.875,-0.15) circle (1pt);
	\filldraw[blue] (-0.875, 0.35) circle (1pt);
	\filldraw[blue] (-0.875, 0.85) circle (1pt);
	\filldraw[blue] (-0.375,-0.65) circle (1pt);
	\filldraw[blue] (-0.375,-0.15) circle (1pt);
	\filldraw[blue] (-0.375, 0.35) circle (1pt);
	\filldraw[blue] (-0.375, 0.85) circle (1pt);
	\filldraw[blue] (0.125,-0.65) circle (1pt);
	\filldraw[blue] (0.125,-0.15) circle (1pt);
	\filldraw[blue] (0.125, 0.35) circle (1pt);
	\filldraw[blue] (0.125, 0.85) circle (1pt);
	\filldraw[blue] (0.625,-0.65) circle (1pt);
	\filldraw[blue] (0.625,-0.15) circle (1pt);
	\filldraw[blue] (0.625, 0.35) circle (1pt);
	\filldraw[blue] (0.625, 0.85) circle (1pt);
\end{tikzpicture}
\vspace{0.6cm}

\textcolor{red}{Even rows} $\xrightarrow{}
\left| 
\begin{aligned}
	&i = (8n - 2)k_1 + 2k_2 \\
	& \text{with } k_1 = 0,..,n-1 \text{ and } k_2 = 0,..,n-1
\end{aligned} \right. $
\begin{equation*}
	A_{ij} = \left\{ \begin{aligned}
		-&\frac{1}{6} \qquad &&\text{ if } j= (2n+1)k_1 + k_2 \\
		-&\frac{1}{12} \qquad &&\text{ if } j= (2n+1)k_1+1 +k_2 \\
		&\frac{1}{6} \qquad &&\text{ if } j= (2n+1)k_1+(1+n)+k_2 \\
		&\frac{1}{12} \qquad &&\text{ if } j= (2n+1)k_1+(1+2n)+k_2  
	\end{aligned}
	\right. 
\end{equation*}

\textcolor{blue}{Odd rows} $\xrightarrow{}
\left| 
\begin{aligned}
	&i = (8n - 2)k_1 +1 + 2k_2 \\
	& \text{with } k_1 = 0,..,n-1 \text{ and } k_2 = 0,..,n-1
\end{aligned} \right. $
\begin{equation*}
	A_{ij} = \left\{ \begin{aligned}
		-&\frac{1}{12} \qquad &&\text{ if } j= (2n+1)k_1 + k_2 \\
		-&\frac{1}{6} \qquad &&\text{ if } j= (2n+1)k_1 + 1 + k_2 \\
		&\frac{1}{6} \qquad &&\text{ if } j= (2n+1)k_1 + (1+n) + k_2\\
		&\frac{1}{12} \qquad &&\text{ if } j=(2n+1)k_1 + (2+2n) + k_2                      
	\end{aligned}
	\right. 
\end{equation*}

\paragraph{First column + n columns}
\begin{tikzpicture}[baseline=-0.1cm]
	\draw[step=0.5cm,gray,very thin] (-1,-1) grid (1,1);
	\draw[gray, very thin] (-1,0.5) -- (-0.5,1); 
	\draw[gray, very thin] (-1,0) -- (0,1); 
	\draw[gray, very thin] (-1,-0.5) -- (0.5,1);
	\draw[gray, very thin] (-1,-1) -- (1,1); 
	\draw[gray, very thin] (-0.5,-1) -- (1, 0.5);
	\draw[gray, very thin] (0,-1) -- (1, 0);
	\draw[gray, very thin] (0.5,-1) -- (1, -0.5);
	\draw[gray, very thin] (-1,-0.5) -- (-0.5, -1);
	\draw[gray, very thin] (-1,0) -- (0,-1);
	\draw[gray, very thin] (-1,0.5) -- (0.5, -1);
	\draw[gray, very thin] (-1,1) -- (1,-1);
	\draw[gray, very thin] (-0.5,1) -- (1, -0.5);
	\draw[gray, very thin] (0,1) -- (1,0);
	\draw[gray, very thin] (0.5,1) -- (1,0.5);
	\filldraw[red] (-0.75,-0.75) circle (1pt);
	\filldraw[red] (-0.75,0.25) circle (1pt);
	\filldraw[red] (-0.25,-0.75) circle (1pt);
	\filldraw[red] (-0.25,0.25) circle (1pt);
	\filldraw[red] (0.25,-0.75) circle (1pt);
	\filldraw[red] (0.25,0.25) circle (1pt);
	\filldraw[red] (0.75,-0.75) circle (1pt);
	\filldraw[red] (0.75,0.25) circle (1pt);
	\filldraw[blue] (-0.75,-0.25) circle (1pt);
	\filldraw[blue] (-0.75,0.75) circle (1pt);
	\filldraw[blue] (-0.25,-0.25) circle (1pt);
	\filldraw[blue] (-0.25,0.75) circle (1pt);
	\filldraw[blue] (0.25,-0.25) circle (1pt);
	\filldraw[blue] (0.25,0.75) circle (1pt);
	\filldraw[blue] (0.75,-0.25) circle (1pt);
	\filldraw[blue] (0.75,0.75) circle (1pt);
\end{tikzpicture} 
\vspace{0.6cm}

\textcolor{red}{Even rows} $\xrightarrow{}$ zeros

\textcolor{blue}{Odd rows} $\xrightarrow{}
\left| 
\begin{aligned}
	&i = (2n + 1) + (8n - 2)k_1 + 2k_2\\
	& \text{with } k_1 = 0,..,n-1 \text{ and } k_2 = 0,..,n-2
\end{aligned} \right. $
\begin{equation*}
	A_{ij} = \left\{ \begin{aligned}
		-&\frac{1}{6} \qquad &&\text{ if } j=(2n+1)k_1 + 1 + k_2 \\
		&\frac{1}{6} \qquad &&\text{ if } j= (2n+1)k_1 + (2n+2) + k_2                      
	\end{aligned}
	\right. 
\end{equation*}

\paragraph{Second column + n columns}
\begin{tikzpicture}[baseline=-0.1cm]
	\draw[step=0.5cm,gray,very thin] (-1,-1) grid (1,1);
	\draw[gray, very thin] (-1,0.5) -- (-0.5,1); 
	\draw[gray, very thin] (-1,0) -- (0,1); 
	\draw[gray, very thin] (-1,-0.5) -- (0.5,1);
	\draw[gray, very thin] (-1,-1) -- (1,1); 
	\draw[gray, very thin] (-0.5,-1) -- (1, 0.5);
	\draw[gray, very thin] (0,-1) -- (1, 0);
	\draw[gray, very thin] (0.5,-1) -- (1, -0.5);
	\draw[gray, very thin] (-1,-0.5) -- (-0.5, -1);
	\draw[gray, very thin] (-1,0) -- (0,-1);
	\draw[gray, very thin] (-1,0.5) -- (0.5, -1);
	\draw[gray, very thin] (-1,1) -- (1,-1);
	\draw[gray, very thin] (-0.5,1) -- (1, -0.5);
	\draw[gray, very thin] (0,1) -- (1,0);
	\draw[gray, very thin] (0.5,1) -- (1,0.5);
	\filldraw[red] (-0.625,-0.875) circle (1pt);
	\filldraw[red] (-0.625,-0.375) circle (1pt);
	\filldraw[red] (-0.625,0.125) circle (1pt);
	\filldraw[red] (-0.625,0.625) circle (1pt);
	\filldraw[red] (-0.125,-0.875) circle (1pt);
	\filldraw[red] (-0.125,-0.375) circle (1pt);
	\filldraw[red] (-0.125,0.125) circle (1pt);
	\filldraw[red] (-0.125,0.625) circle (1pt);
	\filldraw[red] (0.375,-0.875) circle (1pt);
	\filldraw[red] (0.375,-0.375) circle (1pt);
	\filldraw[red] (0.375,0.125) circle (1pt);
	\filldraw[red] (0.375,0.625) circle (1pt);
	\filldraw[red] (0.875,-0.875) circle (1pt);
	\filldraw[red] (0.875,-0.375) circle (1pt);
	\filldraw[red] (0.875,0.125) circle (1pt);
	\filldraw[red] (0.875,0.625) circle (1pt);
	\filldraw[blue] (-0.625,-0.625) circle (1pt);
	\filldraw[blue] (-0.625,-0.125) circle (1pt);
	\filldraw[blue] (-0.625,0.375) circle (1pt);
	\filldraw[blue] (-0.625,0.875) circle (1pt);
	\filldraw[blue] (-0.125,-0.625) circle (1pt);
	\filldraw[blue] (-0.125,-0.125) circle (1pt);
	\filldraw[blue] (-0.125,0.375) circle (1pt);
	\filldraw[blue] (-0.125,0.875) circle (1pt);
	\filldraw[blue] (0.375,-0.625) circle (1pt);
	\filldraw[blue] (0.375,-0.125) circle (1pt);
	\filldraw[blue] (0.375,0.375) circle (1pt);
	\filldraw[blue] (0.375,0.875) circle (1pt);
	\filldraw[blue] (0.875,-0.625) circle (1pt);
	\filldraw[blue] (0.875,-0.125) circle (1pt);
	\filldraw[blue] (0.875,0.375) circle (1pt);
	\filldraw[blue] (0.875,0.875) circle (1pt);
\end{tikzpicture}
\vspace{0.6cm}

\textcolor{red}{Even rows} $\xrightarrow{}
\left| 
\begin{aligned}
	&i = (8n-2)k_1 + (4n-1) + 2k_2\\
	& \text{with } k_1 = 0,..,n-1 \text{ and } k_2 = 0,..,n-1
\end{aligned} \right. $
\begin{equation*}
	A_{ij} = \left\{ \begin{aligned}
		-&\frac{1}{12} \qquad &&\text{ if } j= (2n+1)k_1 + k_2 \\
		-&\frac{1}{6} \qquad &&\text{ if } j= (2n+1)k_1 + (n+1) + k_2 \\
		&\frac{1}{6} \qquad &&\text{ if } j=(2n+1)k_1 + (2n+1) + k_2 \\
		&\frac{1}{12} \qquad &&\text{ if } j= (2n+1)k_1 + (2n+2) + k_2                      
	\end{aligned}
	\right. 
\end{equation*}

\textcolor{blue}{Odd rows} $\xrightarrow{}
\left| 
\begin{aligned}
	&i = (8n-2)k_1 + 4n + 2k_2\\
	& \text{with } k_1 = 0,..,n-1 \text{ and } k_2 = 0,..,n-1
\end{aligned} \right. $
\begin{equation*}
	A_{ij} = \left\{ \begin{aligned}
		-&\frac{1}{12} \qquad &&\text{ if } j= (2n+1)k_1 + 1 + k_2 \\
		-&\frac{1}{6} \qquad &&\text{ if } j= (2n+1)k_1 + (n+1) + k_2 \\
		&\frac{1}{12} \qquad &&\text{ if } j= (2n+1)k_1 + (2n+1) + k_2 \\
		&\frac{1}{6} \qquad &&\text{ if } j= (2n+1)k_1 + (2n+2) + k_2                      
	\end{aligned}
	\right. 
\end{equation*}

\paragraph{First column + (n-1) columns}
\begin{tikzpicture}[baseline=-0.1cm]
	\draw[step=0.5cm,gray,very thin] (-1,-1) grid (1,1);
	\draw[gray, very thin] (-1,0.5) -- (-0.5,1); 
	\draw[gray, very thin] (-1,0) -- (0,1); 
	\draw[gray, very thin] (-1,-0.5) -- (0.5,1);
	\draw[gray, very thin] (-1,-1) -- (1,1); 
	\draw[gray, very thin] (-0.5,-1) -- (1, 0.5);
	\draw[gray, very thin] (0,-1) -- (1, 0);
	\draw[gray, very thin] (0.5,-1) -- (1, -0.5);
	\draw[gray, very thin] (-1,-0.5) -- (-0.5, -1);
	\draw[gray, very thin] (-1,0) -- (0,-1);
	\draw[gray, very thin] (-1,0.5) -- (0.5, -1);
	\draw[gray, very thin] (-1,1) -- (1,-1);
	\draw[gray, very thin] (-0.5,1) -- (1, -0.5);
	\draw[gray, very thin] (0,1) -- (1,0);
	\draw[gray, very thin] (0.5,1) -- (1,0.5);
	\filldraw[red] (-0.5,-0.75) circle (1pt);
	\filldraw[red] (-0.5,-0.25) circle (1pt);
	\filldraw[red] (-0.5,0.25) circle (1pt);
	\filldraw[red] (-0.5,0.75) circle (1pt);
	\filldraw[red] (0,-0.75) circle (1pt);
	\filldraw[red] (0,-0.25) circle (1pt);
	\filldraw[red] (0,0.25) circle (1pt);
	\filldraw[red] (0,0.75) circle (1pt);
	\filldraw[red] (0.5,-0.75) circle (1pt);
	\filldraw[red] (0.5,-0.25) circle (1pt);
	\filldraw[red] (0.5,0.25) circle (1pt);
	\filldraw[red] (0.5,0.75) circle (1pt);
	\filldraw[blue] (-0.5,-0.5) circle (1pt);
	\filldraw[blue] (-0.5,0) circle (1pt);
	\filldraw[blue] (-0.5,0.5) circle (1pt);
	\filldraw[blue] (0,-0.5) circle (1pt);
	\filldraw[blue] (0,0) circle (1pt);
	\filldraw[blue] (0,0.5) circle (1pt);
	\filldraw[blue] (0.5,-0.5) circle (1pt);
	\filldraw[blue] (0.5,0) circle (1pt);
	\filldraw[blue] (0.5,0.5) circle (1pt);
\end{tikzpicture}
\vspace{0.6cm}

\textcolor{red}{Even rows} $\xrightarrow{}
\left| 
\begin{aligned}
	&i = (8n - 2)k_1 + (6n-1) + 2k_2\\
	& \text{with } k_1 = 0,..,n-2 \text{ and } k_2 = 0,..,n-1
\end{aligned} \right. $
\begin{equation*}
	A_{ij} = \left\{ \begin{aligned}
		-&\frac{1}{6} \qquad &&\text{ if } j=(2n+1)k_1 + (n+1) + k_2 \\
		&\frac{1}{6} \qquad &&\text{ if } j= (2n+1)k_1 + (3n+2) + k_2                     
	\end{aligned}
	\right. 
\end{equation*}

\textcolor{blue}{Odd rows} $\xrightarrow{}$ zeros
\newpage

\subsection{Matrix $B_{y,n}$}

\paragraph{Zero column + n columns}
\begin{tikzpicture}[baseline=-0.1cm]
	\draw[step=0.5cm,gray,very thin] (-1,-1) grid (1,1);
	\draw[gray, very thin] (-1,0.5) -- (-0.5,1); 
	\draw[gray, very thin] (-1,0) -- (0,1); 
	\draw[gray, very thin] (-1,-0.5) -- (0.5,1);
	\draw[gray, very thin] (-1,-1) -- (1,1); 
	\draw[gray, very thin] (-0.5,-1) -- (1, 0.5);
	\draw[gray, very thin] (0,-1) -- (1, 0);
	\draw[gray, very thin] (0.5,-1) -- (1, -0.5);
	\draw[gray, very thin] (-1,-0.5) -- (-0.5, -1);
	\draw[gray, very thin] (-1,0) -- (0,-1);
	\draw[gray, very thin] (-1,0.5) -- (0.5, -1);
	\draw[gray, very thin] (-1,1) -- (1,-1);
	\draw[gray, very thin] (-0.5,1) -- (1, -0.5);
	\draw[gray, very thin] (0,1) -- (1,0);
	\draw[gray, very thin] (0.5,1) -- (1,0.5);
	\filldraw[red] (-0.875,-0.875) circle (1pt);
	\filldraw[red] (-0.875,-0.375) circle (1pt);
	\filldraw[red] (-0.875, 0.125) circle (1pt);
	\filldraw[red] (-0.875, 0.625) circle (1pt);
	\filldraw[red] (-0.375,-0.875) circle (1pt);
	\filldraw[red] (-0.375,-0.375) circle (1pt);
	\filldraw[red] (-0.375, 0.125) circle (1pt);
	\filldraw[red] (-0.375, 0.625) circle (1pt);
	\filldraw[red] (0.125,-0.875) circle (1pt);
	\filldraw[red] (0.125,-0.375) circle (1pt);
	\filldraw[red] (0.125, 0.125) circle (1pt);
	\filldraw[red] (0.125, 0.625) circle (1pt);
	\filldraw[red] (0.625,-0.875) circle (1pt);
	\filldraw[red] (0.625,-0.375) circle (1pt);
	\filldraw[red] (0.625, 0.125) circle (1pt);
	\filldraw[red] (0.625, 0.625) circle (1pt);
	\filldraw[blue] (-0.875,-0.65) circle (1pt);
	\filldraw[blue] (-0.875,-0.15) circle (1pt);
	\filldraw[blue] (-0.875, 0.35) circle (1pt);
	\filldraw[blue] (-0.875, 0.85) circle (1pt);
	\filldraw[blue] (-0.375,-0.65) circle (1pt);
	\filldraw[blue] (-0.375,-0.15) circle (1pt);
	\filldraw[blue] (-0.375, 0.35) circle (1pt);
	\filldraw[blue] (-0.375, 0.85) circle (1pt);
	\filldraw[blue] (0.125,-0.65) circle (1pt);
	\filldraw[blue] (0.125,-0.15) circle (1pt);
	\filldraw[blue] (0.125, 0.35) circle (1pt);
	\filldraw[blue] (0.125, 0.85) circle (1pt);
	\filldraw[blue] (0.625,-0.65) circle (1pt);
	\filldraw[blue] (0.625,-0.15) circle (1pt);
	\filldraw[blue] (0.625, 0.35) circle (1pt);
	\filldraw[blue] (0.625, 0.85) circle (1pt);
\end{tikzpicture}
\vspace{0.6cm}

\textcolor{red}{Even rows} $\xrightarrow{}
\left| 
\begin{aligned}
	&i = (8n^2-4n+1) + (8n - 2)k_1 + 2k_2 \\
	& \text{with } k_1 = 0,..,n-1 \text{ and } k_2 = 0,..,n-1
\end{aligned} \right. $
\vspace{0.2cm}

\begin{equation*}
	A_{ij} = \left\{ \begin{aligned}
		-&\frac{1}{6} \qquad &&\text{ if } j= (2n+1)k_1 + k_2 \\
		&\frac{1}{12} \qquad &&\text{ if } j= (2n+1)k_1+1 +k_2 \\
		&\frac{1}{6} \qquad &&\text{ if } j= (2n+1)k_1+(1+n)+k_2 \\
		-&\frac{1}{12} \qquad &&\text{ if } j= (2n+1)k_1+(1+2n)+k_2  
	\end{aligned}
	\right. 
\end{equation*}

\textcolor{blue}{Odd rows} $\xrightarrow{}
\left| 
\begin{aligned}
	&i = (8n^2-4n+1) + (8n - 2)k_1 +1 + 2k_2 \\
	& \text{with } k_1 = 0,..,n-1 \text{ and } k_2 = 0,..,n-1
\end{aligned} \right. $
\begin{equation*}
	A_{ij} = \left\{ \begin{aligned}
		-&\frac{1}{12} \qquad &&\text{ if } j= (2n+1)k_1 + k_2 \\
		&\frac{1}{6} \qquad &&\text{ if } j= (2n+1)k_1 + 1 + k_2 \\
		-&\frac{1}{6} \qquad &&\text{ if } j= (2n+1)k_1 + (1+n) + k_2\\
		&\frac{1}{12} \qquad &&\text{ if } j=(2n+1)k_1 + (2+2n) + k_2                      
	\end{aligned}
	\right. 
\end{equation*}
\vspace{0.5cm}

\paragraph{First column + n columns}
\begin{tikzpicture}[baseline=-0.1cm]
	\draw[step=0.5cm,gray,very thin] (-1,-1) grid (1,1);
	\draw[gray, very thin] (-1,0.5) -- (-0.5,1); 
	\draw[gray, very thin] (-1,0) -- (0,1); 
	\draw[gray, very thin] (-1,-0.5) -- (0.5,1);
	\draw[gray, very thin] (-1,-1) -- (1,1); 
	\draw[gray, very thin] (-0.5,-1) -- (1, 0.5);
	\draw[gray, very thin] (0,-1) -- (1, 0);
	\draw[gray, very thin] (0.5,-1) -- (1, -0.5);
	\draw[gray, very thin] (-1,-0.5) -- (-0.5, -1);
	\draw[gray, very thin] (-1,0) -- (0,-1);
	\draw[gray, very thin] (-1,0.5) -- (0.5, -1);
	\draw[gray, very thin] (-1,1) -- (1,-1);
	\draw[gray, very thin] (-0.5,1) -- (1, -0.5);
	\draw[gray, very thin] (0,1) -- (1,0);
	\draw[gray, very thin] (0.5,1) -- (1,0.5);
	\filldraw[red] (-0.75,-0.75) circle (1pt);
	\filldraw[red] (-0.75,0.25) circle (1pt);
	\filldraw[red] (-0.25,-0.75) circle (1pt);
	\filldraw[red] (-0.25,0.25) circle (1pt);
	\filldraw[red] (0.25,-0.75) circle (1pt);
	\filldraw[red] (0.25,0.25) circle (1pt);
	\filldraw[red] (0.75,-0.75) circle (1pt);
	\filldraw[red] (0.75,0.25) circle (1pt);
	\filldraw[blue] (-0.75,-0.25) circle (1pt);
	\filldraw[blue] (-0.75,0.75) circle (1pt);
	\filldraw[blue] (-0.25,-0.25) circle (1pt);
	\filldraw[blue] (-0.25,0.75) circle (1pt);
	\filldraw[blue] (0.25,-0.25) circle (1pt);
	\filldraw[blue] (0.25,0.75) circle (1pt);
	\filldraw[blue] (0.75,-0.25) circle (1pt);
	\filldraw[blue] (0.75,0.75) circle (1pt);
\end{tikzpicture}
\vspace{0.6cm}

\textcolor{red}{Even rows} $\xrightarrow{}$ zeros\\

\textcolor{blue}{Odd rows} $\xrightarrow{} \left| 
\begin{aligned}
	&i = (8n^2-4n+1) + (2n + 1) + (8n - 2)k_1 + 2k_2\\
	& \text{with } k_1 = 0,..,n-1 \text{ and } k_2 = 0,..,n-2
\end{aligned} \right. $
\begin{equation*}
	A_{ij} = \left\{ \begin{aligned}
		-&\frac{1}{6} \qquad &&\text{ if } j=(2n+1)k_1 + (n+1) + k_2 \\
		&\frac{1}{6} \qquad &&\text{ if } j= (2n+1)k_1 + (n+2) + k_2                      
	\end{aligned}
	\right. 
\end{equation*}
\newpage

\paragraph{Second column + n columns}
\begin{tikzpicture}[baseline=-0.1cm]
	\draw[step=0.5cm,gray,very thin] (-1,-1) grid (1,1);
	\draw[gray, very thin] (-1,0.5) -- (-0.5,1); 
	\draw[gray, very thin] (-1,0) -- (0,1); 
	\draw[gray, very thin] (-1,-0.5) -- (0.5,1);
	\draw[gray, very thin] (-1,-1) -- (1,1); 
	\draw[gray, very thin] (-0.5,-1) -- (1, 0.5);
	\draw[gray, very thin] (0,-1) -- (1, 0);
	\draw[gray, very thin] (0.5,-1) -- (1, -0.5);
	\draw[gray, very thin] (-1,-0.5) -- (-0.5, -1);
	\draw[gray, very thin] (-1,0) -- (0,-1);
	\draw[gray, very thin] (-1,0.5) -- (0.5, -1);
	\draw[gray, very thin] (-1,1) -- (1,-1);
	\draw[gray, very thin] (-0.5,1) -- (1, -0.5);
	\draw[gray, very thin] (0,1) -- (1,0);
	\draw[gray, very thin] (0.5,1) -- (1,0.5);
	\filldraw[red] (-0.625,-0.875) circle (1pt);
	\filldraw[red] (-0.625,-0.375) circle (1pt);
	\filldraw[red] (-0.625,0.125) circle (1pt);
	\filldraw[red] (-0.625,0.625) circle (1pt);
	\filldraw[red] (-0.125,-0.875) circle (1pt);
	\filldraw[red] (-0.125,-0.375) circle (1pt);
	\filldraw[red] (-0.125,0.125) circle (1pt);
	\filldraw[red] (-0.125,0.625) circle (1pt);
	\filldraw[red] (0.375,-0.875) circle (1pt);
	\filldraw[red] (0.375,-0.375) circle (1pt);
	\filldraw[red] (0.375,0.125) circle (1pt);
	\filldraw[red] (0.375,0.625) circle (1pt);
	\filldraw[red] (0.875,-0.875) circle (1pt);
	\filldraw[red] (0.875,-0.375) circle (1pt);
	\filldraw[red] (0.875,0.125) circle (1pt);
	\filldraw[red] (0.875,0.625) circle (1pt);
	\filldraw[blue] (-0.625,-0.625) circle (1pt);
	\filldraw[blue] (-0.625,-0.125) circle (1pt);
	\filldraw[blue] (-0.625,0.375) circle (1pt);
	\filldraw[blue] (-0.625,0.875) circle (1pt);
	\filldraw[blue] (-0.125,-0.625) circle (1pt);
	\filldraw[blue] (-0.125,-0.125) circle (1pt);
	\filldraw[blue] (-0.125,0.375) circle (1pt);
	\filldraw[blue] (-0.125,0.875) circle (1pt);
	\filldraw[blue] (0.375,-0.625) circle (1pt);
	\filldraw[blue] (0.375,-0.125) circle (1pt);
	\filldraw[blue] (0.375,0.375) circle (1pt);
	\filldraw[blue] (0.375,0.875) circle (1pt);
	\filldraw[blue] (0.875,-0.625) circle (1pt);
	\filldraw[blue] (0.875,-0.125) circle (1pt);
	\filldraw[blue] (0.875,0.375) circle (1pt);
	\filldraw[blue] (0.875,0.875) circle (1pt);
\end{tikzpicture}
\vspace{0.6cm}

\textcolor{red}{Even rows} $\xrightarrow{} \left| 
\begin{aligned}
	&i = (8n^2-4n+1) + (8n-2)k_1 + (4n-1) + 2k_2\\
	& \text{with } k_1 = 0,..,n-1 \text{ and } k_2 = 0,..,n-1
\end{aligned} \right. $
\begin{equation*}
	A_{ij} = \left\{ \begin{aligned}
		-&\frac{1}{12} \qquad &&\text{ if } j= (2n+1)k_1 + k_2 \\
		&\frac{1}{6} \qquad &&\text{ if } j= (2n+1)k_1 + (n+1) + k_2 \\
		-&\frac{1}{6} \qquad &&\text{ if } j=(2n+1)k_1 + (2n+1) + k_2 \\
		&\frac{1}{12} \qquad &&\text{ if } j= (2n+1)k_1 + (2n+2) + k_2                      
	\end{aligned}
	\right. 
\end{equation*}

\textcolor{blue}{Odd rows} $\xrightarrow{} \left| 
\begin{aligned}
	&i = (8n^2-4n+1)+ (8n-2)k_1 + 4n + 2k_2\\
	& \text{with } k_1 = 0,..,n-1 \text{ and } k_2 = 0,..,n-1
\end{aligned} \right. $
\begin{equation*}
	A_{ij} = \left\{ \begin{aligned}
		&\frac{1}{12} \qquad &&\text{ if } j= (2n+1)k_1 + 1 + k_2 \\
		-&\frac{1}{6} \qquad &&\text{ if } j= (2n+1)k_1 + (n+1) + k_2 \\
		-&\frac{1}{12} \qquad &&\text{ if } j= (2n+1)k_1 + (2n+1) + k_2 \\
		&\frac{1}{6} \qquad &&\text{ if } j= (2n+1)k_1 + (2n+2) + k_2                      
	\end{aligned}
	\right. 
\end{equation*}
\vspace{0.5cm}

\paragraph{First column + (n-1) columns}
\begin{tikzpicture}[baseline=-0.1cm]
	\draw[step=0.5cm,gray,very thin] (-1,-1) grid (1,1);
	\draw[gray, very thin] (-1,0.5) -- (-0.5,1); 
	\draw[gray, very thin] (-1,0) -- (0,1); 
	\draw[gray, very thin] (-1,-0.5) -- (0.5,1);
	\draw[gray, very thin] (-1,-1) -- (1,1); 
	\draw[gray, very thin] (-0.5,-1) -- (1, 0.5);
	\draw[gray, very thin] (0,-1) -- (1, 0);
	\draw[gray, very thin] (0.5,-1) -- (1, -0.5);
	\draw[gray, very thin] (-1,-0.5) -- (-0.5, -1);
	\draw[gray, very thin] (-1,0) -- (0,-1);
	\draw[gray, very thin] (-1,0.5) -- (0.5, -1);
	\draw[gray, very thin] (-1,1) -- (1,-1);
	\draw[gray, very thin] (-0.5,1) -- (1, -0.5);
	\draw[gray, very thin] (0,1) -- (1,0);
	\draw[gray, very thin] (0.5,1) -- (1,0.5);
	\filldraw[red] (-0.5,-0.75) circle (1pt);
	\filldraw[red] (-0.5,-0.25) circle (1pt);
	\filldraw[red] (-0.5,0.25) circle (1pt);
	\filldraw[red] (-0.5,0.75) circle (1pt);
	\filldraw[red] (0,-0.75) circle (1pt);
	\filldraw[red] (0,-0.25) circle (1pt);
	\filldraw[red] (0,0.25) circle (1pt);
	\filldraw[red] (0,0.75) circle (1pt);
	\filldraw[red] (0.5,-0.75) circle (1pt);
	\filldraw[red] (0.5,-0.25) circle (1pt);
	\filldraw[red] (0.5,0.25) circle (1pt);
	\filldraw[red] (0.5,0.75) circle (1pt);
	\filldraw[blue] (-0.5,-0.5) circle (1pt);
	\filldraw[blue] (-0.5,0) circle (1pt);
	\filldraw[blue] (-0.5,0.5) circle (1pt);
	\filldraw[blue] (0,-0.5) circle (1pt);
	\filldraw[blue] (0,0) circle (1pt);
	\filldraw[blue] (0,0.5) circle (1pt);
	\filldraw[blue] (0.5,-0.5) circle (1pt);
	\filldraw[blue] (0.5,0) circle (1pt);
	\filldraw[blue] (0.5,0.5) circle (1pt);
\end{tikzpicture}
\vspace{0.6cm}

\textcolor{red}{Even rows} $\xrightarrow{} \left| 
\begin{aligned}
	&i = (8n^2-4n+1) + (8n - 2)k_1 + (6n-1) + 2k_2\\
	& \text{with } k_1 = 0,..,n-2 \text{ and } k_2 = 0,..,n-1
\end{aligned} \right. $
\begin{equation*}
	A_{ij} = \left\{ \begin{aligned}
		-&\frac{1}{6} \qquad &&\text{ if } j=(2n+1)k_1 + (2n+1) + k_2 \\
		&\frac{1}{6} \qquad &&\text{ if } j= (2n+1)k_1 + (2n+2) + k_2                     
	\end{aligned}
	\right. 
\end{equation*}

\textcolor{blue}{Odd rows} $\xrightarrow{}$ zeros

\end{document}